\documentclass{article}
\usepackage{color}

\usepackage{amsmath}

\begin{document}

\title{Well posedness of an integrodifferential kinetic model of Fokker-Planck type for angiogenesis}
\author{Ana Carpio \thanks{Departamento de Matematica Aplicada, Universidad Complutense, 28040 Madrid, ana\_carpio@mat.ucm.es}, 
Gema Duro \thanks{Departamento de An\'alisis Econ\'omico: Econom\'{\i}a Cuantitativa, Universidad Aut\'onoma de Madrid, 28049 Madrid, Spain}}

\date{May 16, 2015}

\maketitle

{\bf Abstract.}
Tumor induced angiogenesis processes including the effect of stochastic motion 
and branching of blood vessels can be described coupling a (nonlocal in time) integrodifferential  kinetic equation of Fokker-Planck type with a diffusion equation 
for the tumor induced angiogenic factor. The chemotactic force field depends on the flux of blood vessels through the angiogenic factor. We develop an existence and uniqueness theory for this system under natural assumptions on the initial data. 
The proof combines the construction of fundamental  solutions for  associated linearized problems with comparison principles, sharp estimates of the velocity integrals and compactness results for this type of  kinetic and parabolic operators.
\\


\section{Introduction}
\label{sec:intro}

Angiogenesis (growth of blood vessels) is fundamental for tissue development and repair.  Numerous inflammatory, immune, ischaemic and malignant diseases are fostered by angiogenic disorders \cite{angiogenesis1}. In particular, 
angiogenesis supports tumor spread. Many efforts have been devoted to antiangiogenic therapies to neutralize tumor growth \cite{angiogenesis2}. Current investigations try also to control the formation of new vessels to regenerate damaged tissues and to prevent retinopathies in premature children. To those purposes, it is essential to develop adequate mathematical models of the process, that must be continuously updated to incorporate new experimental knowledge.
There are many models addressing partial aspects of angiogenesis dynamics, see references \cite{angiogenesis3,angiogenesis6,angiogenesis5,angiogenesis4} for instance.
Being able to reproduce the stochastic nature of the vessel branching process seems to be a key novel feature identified in recent experiments. We consider here a deterministic integrodifferential system suitable to describe the development of the stochastic vessel network, as shown by recent numerical studies \cite{capasso}.  We aim to formulate an existence and uniqueness theory for this type of models, that may serve as a tool for a rigorous derivation of the mean field system from the initial stochastic equations and as a basis for numerical schemes.

We consider the following equations for the evolution of the density of vessel tips $p$ in response to the tumor angiogenic factor released by cells $c$: 
\begin{eqnarray} \frac{\partial}{\partial t} p(t,\mathbf{x},\mathbf{v})&=&
 \alpha(c(t,\mathbf{x})) \rho_{}(\mathbf{v})  p(t,\mathbf{x},\mathbf{v})- \gamma p(t,\mathbf{x},\mathbf{v}) \int_0^t \hskip -2mm d\,s \int_{I\!\! R^{N}} 
 \hskip -4mm d{\mathbf v}' p(s,\mathbf{x},\mathbf{v}')  \nonumber\\
 && - \mathbf{v}\cdot \nabla_\mathbf{x}   p(t,\mathbf{x},\mathbf{v}) + k \nabla_\mathbf{v} \cdot (\mathbf{v} p(t,\mathbf{x},\mathbf{v}))
\nonumber\\
&& - \nabla_\mathbf{v} \cdot \left[\mathbf{F}\left(c(t,\mathbf{x})\right)p(t,\mathbf{x},\mathbf{v})  \right]\!+ \sigma
\Delta_\mathbf{v} p(t,\mathbf{x},\mathbf{v}), \label{eq:p}  \\
\frac{\partial}{\partial t}c(t,\mathbf{x}) &=& d \Delta_{\mathbf x} c(t,\mathbf{x}) - \eta c(t,\mathbf{x}) { j(t,\mathbf{x})} \label{eq:c}, \\
p(0,\mathbf{x},\mathbf{v}) &=& p_0(\mathbf{x},\mathbf{v}), \quad
c(0,\mathbf{x}) =c_0(\mathbf{x}), \label{eq:pc0}
\end{eqnarray}
where
\begin{eqnarray}
\alpha(c)=\alpha_1\frac{\frac{c}{c_R}}{1+\frac{c}{c_R}} \geq 0, \quad
 {\mathbf F}(c)= \frac{d_1}{(1+\gamma_1c)^{q_1}}\nabla_{\mathbf x} c,
\label{eq:alphaF} \\
{ j(t,\mathbf{x})}= \int_{I\!\! R^{N}} { {|\mathbf{v}|} 
\over  1+ e^{\delta (|{\mathbf v}|^2 - v_{max}^2) }} p(t,\mathbf{x},\mathbf{v})\, 
d \mathbf{v},
\quad \tilde{p}(t,\mathbf{x})= \int_{I\!\! R^{N}}  p(t,\mathbf{x},\mathbf{v})
\, d \mathbf{v}, \label{eq:intpintvp}
\end{eqnarray}
for ${\mathbf x}, {\mathbf v} \in I\!\! R^{N} \times I\!\! R^{N}$, 
$N  =  2,3,$ $t  \in   [0, \infty).$
The constants $\gamma$, $k$, $\sigma$, $d$, $\eta$, $\alpha_1$, $c_R$, $d_1$, $\gamma_1$, $q_1$, $\delta$ and $v_{max}$ are positive.  In the original model \cite{capasso}, $\rho(\mathbf{v})=\delta_{\mathbf v_0}$ is a Dirac measure, approximated by gaussians for numerical purposes.
Here, we take $\rho(\mathbf{v})$ to be smooth and bounded positive functions
decaying at infinity. The density of vessels $p$ decays at infinity. The tumor angiogenic factor $c$ decays also, except for a finite  region in $(x_2,x_3)$, or $x_2$, depending on the dimension, for which it behaves like  a positive constant as  $x_1 \rightarrow \infty$. 
In dimension two, these models may be adapted to investigate retinopathy associated angiogenesis problems. The retina is a two dimensional membrane. For tumor related angiogenesis, two dimensional studies are a simplification that may guide three dimensional analyses, where blood vessel tip behavior is more complex.  Only two dimensional numerical tests  have been performed so far with this model \cite{capasso}.  Blood vessel behavior in three dimensions might lead to slight variations in the nonlinear terms.

In this model, the source term  $\alpha(c) \rho \, p$
represents vessel tip branching. Vessel death (anastomosis) is described by the integral sink  $- \gamma p \int_0^t \tilde{p}$.  The Fokker-Planck operator represents vessel extension.   As discussed in reference \cite{angiogenesis3}, the chemotatic force ${\bf F}(c)$ may be taken to depend on the  flux of blood vessels {$j$}, or on the marginal tip density $\tilde p$ if the latter replaces  {$j$} in the diffusion equation (\ref{eq:c})  for the tumor angiogenic factor (TAF). If the sink term in \eqref{eq:c} depends on $\tilde p$, it represents consumption of TAF due to all the cells in the network. This is similar to the consumption term appearing in previous reaction-diffusion models  \cite{anderson}. If TAF consumption is mostly due to the  additional endothelial cells that produce vessel extensions, then the sink term in \eqref{eq:c} is proportional to the velocity of the tips and it may be described through an integral ${\mathbf j}= \int \mathbf{v}  g(|\mathbf{v}|)  \,p\, d \mathbf{v}$, $g\geq 0$. A possible choice could be $g(|\mathbf{v}|)=1$. However, the velocity of cells is limited, and usually quite small. In practice, numerical simulations introduce velocity cut-offs, that may be represented through Fermi-Dirac distributions, as done in  definition (\ref{eq:intpintvp}) for {$j$.
The magnitude entering the sink term in equation (\ref{eq:c})  must be a scalar.
The choice $|\mathbf j|$ (euclidean norm of $\mathbf j$) was explored numerically in reference \cite{capasso}. However, from the modeling point of view it might introduce artifacts when $|\mathbf j|$ vanishes,  which might 
happen under particular, though unusual, symmetry conditions. The cells
consume TAF even in such situations, therefore the sink term in equation (\ref{eq:c})  cannot vanish.
The choice $j$ made in definition (\ref{eq:intpintvp}) avoids such artifacts.
From the mathematical point of view, the presence of the density $p$ 
inside a square root  may introduce lipschitzianity and uniqueness problems when $|\mathbf j|$ vanishes. This situation is avoided in the simulations
performed in finite spatial domains in reference \cite{capasso}. Working
in unbounded spatial domains, $|\mathbf j|$ vanishes at infinity, and
we must work with $j$ to be able to prove uniqueness results. Existence
proofs follow the same lines in both cases, choosing either $|\mathbf j|$ or $j$.}

Equation (\ref{eq:p}) reminds of Vlasov-Poisson-Fokker-Planck problems. There are, however, two important differences. First, the force field ${\mathbf F}$ is not computed from $\tilde{p}$ solving a Poisson equation. It involves the gradient of solutions of heat equations whose sources depend on the flux of blood vessels 
{$j(p)$}. Second, a nonlocal in time integrodifferential sink involving $\int_0^t \tilde{p}$ is present. This source is negative, therefore it may interfere with  the expected positive sign of $p$ unless handled properly.

During the last 30 years, intensive efforts have been devoted to establish global existence, uniqueness and regularity results for Vlasov-Poisson-Fokker-Planck problems. 
A pioneering paper by Degond \cite{degond} proved existence and uniqueness of solutions  in the whole space. Such solutions were global in dimensions $n<3$ but local  in dimension $n=3$. He developed a weak existence theory for the linear problem, including a maximum principle and bounds independent of the force field, and put forward an iterative scheme to construct solutions of the nonlinear problem, exploiting velocity decay to obtain $L^{\infty}$ estimates and to control the velocity moments by means of interpolation inequalities. Victory and O'Dwyer  \cite{victoryclassical} formulated  a theory of fundamental solutions for linear Fokker-Planck equations in the whole space, that they exploited to construct classical  solutions. They were able to adapt  methods for parabolic equations \cite{friedman}, overcoming the difficulties created by the degenerate diffusion.
Existence of global solutions in three dimensions for some classes of initial data followed \cite{victoryweak,rein}. The work on kinetic equations and renormalized solutions by Diperna, Lions, Perthame  and other authors promoted the appearance of global existence and smoothness analyses under milder assumptions on the initial data. Novel compactness results for transport operators and studies of the propagation of moments were essential \cite{diperna, perthame}.  A series of papers by Bouchut and coworkers established existence of global smooth solutions in three dimensions and studied their long  time behavior, see \cite{bouchut,bouchut2} and references therein. More detailed analysis of long time asymptotics in the whole space followed, see \cite{carpio} and references therein. The study of Fokker-Planck type problems in bounded spatial domains, that we do not consider here, poses additional problems due to the interaction of the transport operator with the boundaries \cite{chen}.

The knowledge accumulated for Vlasov-Poisson-Fokker-Planck systems will 
serve as background to build our theory. 
The key idea to develop an existence theory for our angiogenesis model (\ref{eq:p})-(\ref{eq:intpintvp})  is to include the integrodifferential term in the reference linear operator. We will design an iterative scheme in which the velocity integrals $\tilde{p}$ and {$j$} are frozen from one step to the next. This allows us to construct solutions of the linearized Fokker-Plank and heat problems using fundamental solutions. It also ensures nonnegativity  of the solutions, that is crucial to obtain preliminary uniform estimates.
Proving compactness requires sharp estimates on the force field ${\mathbf F}(c)$  and on the anastomosis sink in terms of $p$. This will be achieved expressing the angiogenic factor $c$ as a convolution with heat kernels and estimating the velocity integrals of $p$ exploiting the decay in the velocity space. 
Specific compactness results for this type of kinetic transport operators will be essential  to extract  convergent subsequences. Passing to the limit 
in the equations satisfied by them, we obtain a solution of the original problem. We construct  global in time solutions in dimensions two and three.
Stability bounds in terms of the norms of the initial data follow. Uniqueness ensues from coupled Gronwall type inequalities. These solutions are strong in the sense that the force field $\mathbf F$ is a bounded function. 

The numerical simulations in \cite{capasso} start from Gaussian initial data. We will consider here milder assumptions.  In the sequel,  $L^{\infty}$ will denote the space of measurable and bounded functions, $L^{1}$ the space of measurable and integrable functions, and $L^{q}$ the space of measurable functions $p$ for which $|p|^q$ is integrable. Subindices will make explicit the involved variables. The mathematical structure and the physical interpretation of the system suggest  some natural hypotheses on the data:
\begin{itemize} 
\item $p_0 \geq 0$, $c_0 \geq 0$, so that $p$ and $c$ are both nonnegative.
\item $c_0 \in L^{\infty}_{\mathbf x}$, $\nabla_{\mathbf x} c_0 \in L^{\infty}_{\mathbf x}$, so that $c$, $\nabla_{\mathbf x} c$ and ${\mathbf F}$ may be bounded functions.
\item $p_0\in L^1_{{\mathbf x}{\mathbf v}}$,  so that we can define the mass 
$\int\int p \, d{\mathbf x}d{\mathbf v}$. 
\item $p_0 \in L^{\infty}_{\mathbf x}L^1_{\mathbf v}$, so that 
$\tilde{p}(t,{\mathbf x}) =\int p \, d{\mathbf v}$ may be a bounded function and we can define fundamental solutions for $p$.
\item  $|{\mathbf v}| g(|{\mathbf v}|) p_0  \in  L^{\infty}_{\mathbf x}L^1_{\mathbf v}$, 
so that  { $j(t,{\mathbf x}) =\int |{\mathbf v}| g(|{\mathbf v}|) \, p \, d{\mathbf v}$} may be a bounded function and  we can define fundamental solutions for $c$.
\end{itemize}
 When { $j$} is given by formula (\ref{eq:intpintvp}), the
technical hypotheses
\begin{itemize}
\item  $(1+|{\mathbf v}|^2)^{\beta/2} p_0 \in L^{\infty}_{{\mathbf x}{\mathbf v}} 
\cap L^1_{{\mathbf x}{\mathbf v}}$,  $\beta > N$,  
and $\nabla_{\mathbf x} c_0 \in L^2_{\mathbf x}$,
\item $\nabla_{\mathbf v} p_0 \in L^{\infty}_{\mathbf x}L^1_{\mathbf v}$,
\end{itemize}
will be needed to establish existence, and to ensure uniqueness of the
solutions, respectively.

The paper is organized as follows. Section 2 recalls basic results on linear Fokker-Planck problems: fundamental solutions, integral expressions, bounds in terms of data and sources and comparison principles. Section 3 extends these results to 
linear problems including an additional linear term in $p$. Section 4 establishes 
$L^1$ and $L^{\infty}$ estimates on the velocity integrals 
{ $j(p)$} and $\tilde{p}$ for the nonlinear problem using velocity decay.  Section 5 introduces the iterative scheme and proves  the global existence and uniqueness result. Studying long time asymptotics remains an open problem.
Section 6 proves a local existence  and uniqueness result for chemotactic forces depending on fluxes without velocity cut-offs. Global existence is an open issue 
in this case. Existence of solutions
when the coefficient $\rho(\mathbf v)$ is a Dirac measure remains
an open question too, since we lack estimates to pass to the limit
in the regularized problems that approximate the delta function
by gaussians.

\section{Initial value problems for linear Fokker-Planck
operators}
\label{sec:linearvfp}

The initial value problem (IVP) for linear Fokker-Planck
operators:
\begin{eqnarray} 
\frac{\partial}{\partial t} p(t,\mathbf{x},\mathbf{v}) \!+\! 
\mathbf{v} \!\cdot\! \nabla_\mathbf{x}   p(t,\mathbf{x},\mathbf{v}) 
\!+\!  \nabla_\mathbf{v} \!\cdot\! [({\mathbf F}(t,\mathbf{x}) \!-\! k\mathbf{v} ) p(t,\mathbf{x},\mathbf{v}) ] \!-\! \sigma
\Delta_\mathbf{v} p(t,\mathbf{x},\mathbf{v}) \nonumber \\
= f(t,\mathbf{x},\mathbf{v}),  \label{eq:linp} \\
p(0,\mathbf{x},\mathbf{v})=p_0(\mathbf{x},\mathbf{v}),  
\label{eq:linp0}
\end{eqnarray}
in the whole space has been studied in detail, see \cite{degond,victoryclassical, victoryweak, rein, carpio} and references therein. 
We are particularly  interested here in existence theories constructed using fundamental solutions, that will be the basis of our later analysis of the integrodifferential angiogenesis model. Fundamental solutions are known in stochastic settings as transition probability densities. They furnish a bridge between deterministic and stochastic formulations. Fundamental solutions have the advantage of providing a flexible framework to establish positivity results, maximum principles and estimates on spatiotemporal decay.

Let us set $Q_T= [0,T] \times I\!\!R^{N}  \times I\!\!R^{N}$. A function $\Gamma_{\mathbf F}(t, {\mathbf x}, {\mathbf v};  \tau, {\boldsymbol \xi}, {\boldsymbol \nu})$ defined for $(t, {\mathbf x}, {\mathbf v})  \in Q_T$, 
$(\tau, {\boldsymbol \xi}, {\boldsymbol \nu}) \in Q_T$, $t>\tau,$  is a fundamental
solution of the initial value problem (\ref{eq:linp})-(\ref{eq:linp0}) when:
\begin{itemize}
\item[i)] For $(\tau, {\boldsymbol \xi}, {\boldsymbol \nu}) \in Q_T$ fixed, it satisfies the equation (\ref{eq:linp}) with $f=0$ as a function of $({\mathbf x}, {\mathbf v}) \in I\!\!R^{N}  \times I\!\!R^{N}$, $t \geq \tau$,
\item[ii)] For every continuous and bounded function $p_0({\mathbf x}, {\mathbf v})$
\[
{\rm lim}_{t \rightarrow \tau} \int_{I\!\!R^{N}}\int_{I\!\!R^{N}}
\Gamma_{\mathbf F}(t, {\mathbf x}, {\mathbf v};  \tau, {\boldsymbol \xi}, {\boldsymbol \nu}) p_0({\boldsymbol \xi}, {\boldsymbol \nu}) d{\boldsymbol \xi} d{\boldsymbol \nu}
= p_0({\mathbf x}, {\mathbf v}).
\]
\end{itemize}

When ${\mathbf F}=0$, the field free fundamental solution of equation (\ref{eq:linp}) has a known explicit form \cite{chandrasehkar,victoryclassical}:
\begin{eqnarray}
G(t, {\mathbf x}, {\mathbf v};  \tau, {\boldsymbol \xi}, {\boldsymbol \nu})=
\left[ {k \, e^{k(t-\tau)} \over 4\pi \sigma \sqrt{
{e^{2k(t-\tau)}-1 \over 2 k}(t-\tau) - {(e^{k(t-\tau)}-1)^2 \over k^2}
} }  \right]^{N} \nonumber \\
exp\left\{ -{(k{\mathbf x} \!-\! k{\boldsymbol \xi}
+{\mathbf v}\!-\!{\boldsymbol \nu})^2 \over 4 \sigma(t-\tau)}
- { \big| {e^{k(t-\tau)}-1 \over (t-\tau)} ({\mathbf x}\!-\!{\boldsymbol \xi} 
+ { {\mathbf v}-{\boldsymbol \nu} \over k}) + 
({\boldsymbol \nu} \!-\! {\mathbf v} e^{k(t-\tau)}) \big|^2 \over
4 \sigma \left[ {e^{2k(t-\tau)}-1 \over 2k} - {(e^{k(t-\tau)}-1)^2 
\over k^2 (t-\tau)}
\right]
}
\right\}. 
\label{sf0}
\end{eqnarray}

Some relevant properties are summarized in references \cite{rein,bouchut,bouchut2,bouchut3}. We recall below some properties to be used throughout the paper. \\

{\bf Lemma 2.1.} {\it The fundamental solution $G(t, {\mathbf x}, {\mathbf v};  
\tau, {\boldsymbol \xi}, {\boldsymbol \nu})$
of the field free initial value problem is positive and all its relevant derivatives  exist in a classical sense. 
For $(t, {\mathbf x}, {\mathbf v}), (\tau, {\boldsymbol \xi}, {\boldsymbol \nu}) 
\in Q_{\infty}$, $t > \tau$, it satisfies the following relations:
\begin{itemize}
\item[(G1)] $\displaystyle \int_{I\!\!R^{N}}\int_{I\!\!R^{N}} G(t, {\mathbf x}, {\mathbf v};  
\tau, {\boldsymbol \xi}, {\boldsymbol \nu})  d{\boldsymbol \xi} d{\boldsymbol \nu} = e^{{N}k(t-\tau)}$, \\
$\displaystyle  \int_{I\!\!R^{N}}\int_{I\!\!R^{N}} G(t, {\mathbf x}, {\mathbf v};  
\tau, {\boldsymbol \xi}, {\boldsymbol \nu})  d{\mathbf x} d{\mathbf v} = 1.$\\
\item[(G2)] For $t > \tau' > \tau \geq 0$, 
\[
\int_{I\!\!R^{N}} \int_{I\!\!R^{N}}
G(t, {\mathbf x}, {\mathbf v};   \tau', {\boldsymbol \xi}', {\boldsymbol \nu}') 
G( \tau', {\boldsymbol \xi}', {\boldsymbol \nu}';   \tau, {\boldsymbol \xi}, {\boldsymbol \nu}) d{\boldsymbol \xi}' d{\boldsymbol \nu}'
=  
G(t, {\mathbf x}, {\mathbf v};    \tau, {\boldsymbol \xi}, {\boldsymbol \nu}). \]
\item[(G3)] When $z_i=x_i$ or $\xi_i$, $i=1,\ldots,N$,
\[ |\partial_{z_i} G(t, {\mathbf x}, {\mathbf v};  
\tau, {\boldsymbol \xi}, {\boldsymbol \nu})| \leq M
{\displaystyle G(t, {\mathbf x}/2, {\mathbf v}/2;  \tau, {\boldsymbol \xi}/2, {\boldsymbol \nu}/2)
\over (t-\tau)^{3/2}}.\]
\item[(G4)] When $z_i=v_i$ or $\nu_i$, $i=1,\ldots,N$,
\[ |\partial_{z_i} G(t, {\mathbf x}, {\mathbf v};  
\tau, {\boldsymbol \xi}, {\boldsymbol \nu})| \leq M
{\displaystyle G(t, {\mathbf x}/2, {\mathbf v}/2;  \tau, {\boldsymbol \xi}/2, {\boldsymbol \nu}/2)
\over (t-\tau)^{1/2}}.\]
\item[(G5)] $ |\Delta_{\mathbf v} G(t, {\mathbf x}, {\mathbf v};  
\tau, {\boldsymbol \xi}, {\boldsymbol \nu})| \leq M
{\displaystyle G(t, {\mathbf x}/2, {\mathbf v}/2;  \tau, {\boldsymbol \xi}/2, {\boldsymbol \nu}/2)
\over (t-\tau)}.$
\item[(G6)] $ |\partial_{v_i} \partial_{x_j}
G(t, {\mathbf x}, {\mathbf v};  
\tau, {\boldsymbol \xi}, {\boldsymbol \nu})| \leq M
{\displaystyle G(t, {\mathbf x}/2, {\mathbf v}/2;  \tau, {\boldsymbol \xi}/2, {\boldsymbol \nu}/2) \over (t-\tau)^2},$ $\; i,j=1, \ldots, N$.
\end{itemize}
}

%
%
%
%
%
%
%

When a field ${\mathbf F}$ is added, fundamental solutions can be constructed as a perturbation of $G$. This is done in reference \cite{victoryclassical} for smooth bounded ${\mathbf F}$. We recall below a collection of useful estimates. \\

{\bf Lemma 2.2.} {\it Assume that ${\mathbf F}(t,{\mathbf x})$ is continuous and bounded with respect to $t \in [0,T]$ and continuously differentiable with respect to ${\mathbf x}$, with bounded derivatives in $I\!\!R^{N}  \times I\!\!R^{N}.$ 
Then, there exists a fundamental solution $\Gamma_{\mathbf F}(t, {\mathbf x}, {\mathbf v};  \tau, {\boldsymbol \xi}, {\boldsymbol \nu})$ of the initial value problem (\ref{eq:linp})-(\ref{eq:linp0}) in $Q_T$ whose relevant derivatives exist in a classical sense.
When $(t, {\mathbf x}, {\mathbf v}), (\tau, {\boldsymbol \xi}, {\boldsymbol \nu}) \in Q_T$, $t > \tau$, it satisfies the following identities and bounds:
\begin{itemize}
\item[(F1)] $\displaystyle \int_{I\!\!R^{N}} \! \int_{I\!\!R^{N}} 
\hskip -3mm \Gamma_{\mathbf F}(t, {\mathbf x}, {\mathbf v};  
\tau, {\boldsymbol \xi}, {\boldsymbol \nu})  d{\boldsymbol \xi} d{\boldsymbol \nu} = e^{{N}k(t-\tau)}$, \\
$\displaystyle  \int_{I\!\!R^{N}} \! \int_{I\!\!R^{N}} 
\hskip -3mm \Gamma_{\mathbf F}(t, {\mathbf x}, {\mathbf v};  
\tau, {\boldsymbol \xi}, {\boldsymbol \nu})  d{\mathbf x} d{\mathbf v} = 1.$

\item[(F2)] For $t > \tau' > \tau \geq 0$, 
\[
\int_{I\!\!R^{N}} \! \int_{I\!\!R^{N}} \hskip -3mm
\Gamma_{{\mathbf F}}(t, {\mathbf x}, {\mathbf v};   \tau', {\boldsymbol \xi}', {\boldsymbol \nu}') 
\Gamma_{{\mathbf F}}( \tau', {\boldsymbol \xi}', {\boldsymbol \nu}';  
 \tau, {\boldsymbol \xi}, {\boldsymbol \nu}) d{\boldsymbol \xi}' d{\boldsymbol \nu}'
=  
\Gamma_{{\mathbf F}}(t, {\mathbf x}, {\mathbf v};    \tau, {\boldsymbol \xi}, {\boldsymbol \nu}). \]

\item[(F3)] $\Gamma_{{\mathbf F}}$ is a solution of the integral equations:
\begin{eqnarray*}
\Gamma_{{\mathbf F}}(t, {\mathbf x}, {\mathbf v};   \tau, {\boldsymbol \xi}, 
{\boldsymbol \nu}) =
G(t, {\mathbf x}, {\mathbf v};   \tau, {\boldsymbol \xi}, {\boldsymbol \nu}) + 
\\
 \int_{\tau}^t \! \int_{I\!\!R^{N}} \! \int_{I\!\!R^{N}} \partial_{{\boldsymbol \nu}'}
G(t, {\mathbf x}, {\mathbf v};   \tau', {\boldsymbol \xi}', {\boldsymbol \nu}')
{\mathbf F}(\tau',{\boldsymbol \xi}')
\Gamma_{{\mathbf F}}(\tau', {\boldsymbol \xi}', {\boldsymbol \nu}';   \tau, {\boldsymbol \xi}, {\boldsymbol \nu})  d{\boldsymbol \xi}' d{\boldsymbol \nu}' 
d \tau' \\
= G(t, {\mathbf x}, {\mathbf v};   \tau, {\boldsymbol \xi}, {\boldsymbol \nu}) -\\
\int_{\tau}^t \! \int_{I\!\!R^{N}} \! \int_{I\!\!R^{N}}
G(t, {\mathbf x}, {\mathbf v};   \tau', {\boldsymbol \xi}', {\boldsymbol \nu}')
{\mathbf F}(\tau',{\boldsymbol \xi}')
\partial_{{\boldsymbol \nu}'} \Gamma_{{\mathbf F}}(\tau', {\boldsymbol \xi}', {\boldsymbol \nu}';   \tau, {\boldsymbol \xi}, {\boldsymbol \nu})  d{\boldsymbol \xi}' d{\boldsymbol \nu}' d \tau'.
\end{eqnarray*}

\item[(F4)] $0 \leq \Gamma_{\mathbf F}(t, {\mathbf x}, {\mathbf v};  
\tau, {\boldsymbol \xi}, {\boldsymbol \nu}) \leq M(T,\|\mathbf F\|_{\infty})
\, G(t, {\mathbf x}/2, {\mathbf v}/2;  \tau, {\boldsymbol \xi}/2, {\boldsymbol \nu}/2).$

\item[(F5)] $ |\partial_{v_i} \Gamma_{\mathbf F}(t, {\mathbf x}, {\mathbf v};  
\tau, {\boldsymbol \xi}, {\boldsymbol \nu})| \leq M(T,\|\mathbf F\|_{\infty})
{\displaystyle G(t, {\mathbf x}/2, {\mathbf v}/2;  \tau, {\boldsymbol \xi}/2, {\boldsymbol \nu}/2)
\over (t-\tau)^{1/2}}.$

\item[(F6)] $ |\partial_{x_i} \Gamma_{\mathbf F}(t, {\mathbf x}, {\mathbf v};  
\tau, {\boldsymbol \xi}, {\boldsymbol \nu})| \leq M(T,\|\mathbf F\|_{\infty},
\|\nabla_{\mathbf x}\mathbf F\|_{\infty})
{\displaystyle G(t, {\mathbf x}/2, {\mathbf v}/2;  \tau, {\boldsymbol \xi}/2, {\boldsymbol \nu}/2)
\over (t-\tau)^{3/2}}.$

\item[(F7)] $ |\Delta_{\mathbf v} \Gamma_{\mathbf F}(t, {\mathbf x}, {\mathbf v};  
\tau, {\boldsymbol \xi}, {\boldsymbol \nu})| \leq M(T,\|\mathbf F\|_{\infty},
\|\nabla_{\mathbf x}\mathbf F\|_{\infty})
{\displaystyle G(t, {\mathbf x}/2, {\mathbf v}/2;  \tau, {\boldsymbol \xi}/2, {\boldsymbol \nu}/2)
\over (t-\tau)}.$
\end{itemize}
}

Bounds on higher order derivatives and derivatives with respect to the dual variables ${\boldsymbol \xi}, {\boldsymbol \nu}$ similar to those for $G$ hold too. Once the fundamental solutions are known, the solution of the initial value problem (\ref{eq:linp})-(\ref{eq:linp0}) can be explicitly constructed as volume integrals involving the fundamental solution and the data \cite{victoryclassical}.\\

{\bf Lemma 2.3.} {\it Assume that $p_0$ and $f$ are bounded, continuous, locally H\"older continuous in ${\mathbf x}$, uniformly in $t$ and ${\mathbf v}$ with exponent bigger than $2/3$, and integrable in ${\mathbf x},{\mathbf v}$ for $t\in [0,T]$. Then, the initial value problem (\ref{eq:linp})-(\ref{eq:linp0}) has a unique solution
\begin{eqnarray}
p(t,{\mathbf x},{\mathbf v})= \int_{I\!\!R^{N}}\int_{I\!\!R^{N}}
\Gamma_{\mathbf F}(t, {\mathbf x}, {\mathbf v};  0, {\boldsymbol \xi}, {\boldsymbol \nu})  
p_0({\boldsymbol \xi}, {\boldsymbol \nu}) d{\boldsymbol \xi} d{\boldsymbol \nu} + \nonumber \\
\int_0^t \int_{I\!\!R^{N}}\int_{I\!\!R^{N}}
\Gamma_{\mathbf F}(t, {\mathbf x}, {\mathbf v};  \tau, {\boldsymbol \xi}, {\boldsymbol \nu})  
f(\tau, {\boldsymbol \xi}, {\boldsymbol \nu}) d{\boldsymbol \xi} d{\boldsymbol \nu} d \tau.
\label{solint}
\end{eqnarray}
The solution $p$  is continuous in $Q_T$, continuously differentiable in $(0,T] \times I\!\!R^{N} \times I\!\!R^{N}$ (once with respect to ${\mathbf x}$ and twice with respect to ${\mathbf v}$), with bounded derivatives. Moreover, for $t\in [0,T]$
\begin{eqnarray}
\|p(t)\|_1 \leq  \|p_0\|_1 + \int_0^t  \|f(\tau)\|_1 d\tau, \label{int1}\\
\|p(t)\|_{\infty} \leq e^{Nkt}\|p_0\|_{\infty} + \int_0^t 
e^{Nk(t-\tau)} \|f(\tau)\|_{\infty} d\tau. 
\label{intinf}
\end{eqnarray}
}

These results have been extended to measurable bounded fields 
${\mathbf F}(t,{\mathbf x})$ and integrable or bounded data in reference \cite{carpio}
(see also Proposition 3.1 in the next section, and reference \cite{gema}). The time dependence on $T$ of the constants appearing in bounds of the fundamental solution can be removed when $\|{\mathbf F}(t)\|_{\infty}$ decays in time \cite{carpio}.  This provides global  in time estimates of the spatial and temporal behavior of the fundamental solutions and the solutions of the initial value problem \cite{carpio}. \\

{\bf Lemma 2.4.} {\it If ${\mathbf F}(t,{\mathbf x})$ is measurable and bounded, a generalized fundamental solution $\Gamma_{\mathbf F}(t, {\mathbf x}, 
{\mathbf v};  \tau, {\boldsymbol \xi}, {\boldsymbol \nu})$ of the initial value problem (\ref{eq:linp})-(\ref{eq:linp0}) exists. It satisfies the identities and estimates (F1)-(F5).  Moreover, if $\|{\mathbf F}(t)\|_{\infty}$ decays as  $C_{\alpha} (1+t)^{-\alpha-1/2}$ for $\alpha\geq 0$ (and $C_{\alpha}$ is small enough in case $\alpha=0$) the bounds hold uniformly for any positive $t$.
The integral expression (\ref{solint}) is the unique weak solution of the initial value problem (\ref{eq:linp})-(\ref{eq:linp0}) when 
$p_0 \in L^\infty \cap L^1(I\!\!R^{N} \times I\!\!R^{N})$ and
$f \in L^\infty(0,T;L^\infty \cap L^1(I\!\!R^{N} \times I\!\!R^{N})).$ It
satisfies bounds (\ref{int1})-(\ref{intinf}).
For all $\phi \in C_c^2([0,T) \times I\!\!R^{N} \times I\!\!R^{N})$, we have   
\begin{eqnarray}
\int_0^T \hskip -3mm \int_{I\!\!R^{N} \times I\!\!R^{N}} \hskip -8mm
p(t, {\mathbf x}, {\mathbf v}) \Big[
\frac{\partial}{\partial t} \phi(t,\mathbf{x},\mathbf{v}) \!+\! 
\mathbf{v} \!\cdot\! \nabla_\mathbf{x}   \phi(t,\mathbf{x},\mathbf{v}) 
\!+\!   [({\mathbf F}(t,\mathbf{x}) \!-\! k\mathbf{v} ] \!\cdot\!  \nabla_\mathbf{v} 
\phi(t,\mathbf{x},\mathbf{v})   \nonumber \\  \!+ \sigma
\Delta_\mathbf{v} \phi(t,\mathbf{x},\mathbf{v}) \Big] 
d{\mathbf x}  d{\mathbf v} dt 
+ \int_{I\!\!R^{N} \times I\!\!R^{N}} \phi(0,\mathbf{x},\mathbf{v}) 
p_0(\mathbf{x},\mathbf{v}) d{\mathbf x}  d{\mathbf v} =0. \label{weak}
\end{eqnarray}
}

The integral expression for the solutions of the initial value problem together 
with the nonnegativity of the fundamental solutions yield a straightforward 
comparison principle: \\

{\bf Lemma 2.5.} {\it Let $p^{(1)}$, $p^{(2)}$ be solutions of the initial value problem (\ref{eq:linp})-(\ref{eq:linp0}) with data $p^{(1)}_0\leq p^{(2)}_0$
and $f^{(1)} \leq f^{(2)}$ in $Q_T.$ Then, $p^{(1)}\leq p^{(2)}$ 
in $Q_T$.} \\

To prove strong convergences later on we will need adequate upper solutions 
provided by the following Lemma. \\

{\bf Lemma 2.6.} {\it Assume that ${\mathbf F}(t,{\mathbf x})$ is measurable and bounded in $(0,T)\times I\!\!R^{N}.$ Let $U$ be a solution of 
\begin{eqnarray}
\frac{\partial}{\partial t} U \!+\! 
\mathbf{v} \!\cdot\! \nabla_\mathbf{x}   U
\!+\!  \nabla_\mathbf{v} \!\cdot\! [({\mathbf F} \!-\! k\mathbf{v} ) U) ] \!-\! \sigma
\Delta_\mathbf{v} U  \!=\! \alpha_1 \rho_1 U
\quad \mbox{in} \; Q_T,
\label{upper1}
\end{eqnarray}
such that $U(0)=p_0$, for positive constants $\sigma, k, \alpha_1, \rho_1$.
Then, $U \leq {\cal P} = e^{\alpha_1 \rho_1 t} P$, $P$ being a solution of
\begin{eqnarray}
\frac{\partial}{\partial t} P \!+\! 
\mathbf{v} \!\cdot\! \nabla_\mathbf{x}   P
\!-\!  k \nabla_\mathbf{v} \!\cdot\! (\mathbf{v} P)  \!-\! 4 \sigma
\Delta_\mathbf{v} P  \!= 0
\quad \mbox{in} \; Q_T,
\label{upper2}
\end{eqnarray}
taking the initial value $P(0)=2^{2N} M(T, \|{\mathbf F}\|_{\infty}) p_0$,
where $M(T, \|{\mathbf F}\|_{\infty})$ is the constant appearing in inequality 
(F4) in Lemma 2.2.
}\\

{\bf Proof.}
We write $U= e^{\alpha_1 \rho_1 t} {\cal U}$, where ${\cal U}$ is the solution of
\begin{eqnarray} 
\frac{\partial}{\partial t} {\cal U}  \!+\! 
\mathbf{v} \!\cdot\! \nabla_\mathbf{x}   {\cal U}
\!+\!  \nabla_\mathbf{v} \!\cdot\! [({\mathbf F}
\!-\! k\mathbf{v} ) {\cal U} ]   
\!-\! \sigma \Delta_\mathbf{v} {\cal U} = 0.  
\nonumber 
\end{eqnarray}
${\cal U}$ admits the integral expression:
\begin{eqnarray}
{\cal U}(t,{\mathbf x},{\mathbf v})= \int_{I\!\!R^{N} \times I\!\!R^{N}}
\hskip -8mm
\Gamma_{\mathbf F}(t, {\mathbf x}, {\mathbf v};  
0, {\boldsymbol \xi}, {\boldsymbol \nu})  
p_0({\boldsymbol \xi}, {\boldsymbol \nu}) d{\boldsymbol \xi} d{\boldsymbol \nu}.
\nonumber 
\end{eqnarray}
Thanks to estimate (F4) in Lemma 2.2: 
\begin{eqnarray}
{\cal U}(t,{\mathbf x},{\mathbf v}) \!\leq\!  
M(T, \|{\mathbf F}\|_{\infty})  \!\! \int_{I\!\!R^{N} \times I\!\!R^{N}} 
\hskip -7mm 
G(t, {{\mathbf x}\over 2}, {{\mathbf v}\over 2};  0, {{\boldsymbol \xi}\over 2}, {{\boldsymbol \nu}\over 2}) p_0({\boldsymbol \xi}, {\boldsymbol \nu}) 
d{\boldsymbol \xi} d{\boldsymbol \nu}
= P(t,{\mathbf x},{\mathbf v}).
\nonumber 
\end{eqnarray}




\section{Fokker-Planck problems with additional potentials}
\label{sec:potential}

The nonlinear angiogenesis model (\ref{eq:p})-({\ref{eq:intpintvp}}) includes an integrodifferential term in the right hand side that becomes negative when the density $p$ is positive. If we keep it in the source term it becomes hard to guarantee that $p$ remains nonnegative. This suggests considering linear Fokker-Planck problems with an additional lower order term:
\begin{eqnarray} 
\frac{\partial}{\partial t} p(t,\mathbf{x},\mathbf{v}) \!+\! 
\mathbf{v} \!\cdot\! \nabla_\mathbf{x}   p(t,\mathbf{x},\mathbf{v}) 
\!+\!  \nabla_\mathbf{v} \!\cdot\! [({\mathbf F}(t,\mathbf{x}) \!-\! k\mathbf{v} ) 
p(t,\mathbf{x},\mathbf{v}) ] \!-\! \sigma
\Delta_\mathbf{v} p(t,\mathbf{x},\mathbf{v}) \nonumber \\
+ a(t,\mathbf{x},\mathbf{v}) p(t,\mathbf{x},\mathbf{v}) = f(t,\mathbf{x},\mathbf{v}),  
\label{eq:linpa} \\
p(0,\mathbf{x},\mathbf{v})=p_0(\mathbf{x},\mathbf{v}), 
\label{eq:linp0a}
\end{eqnarray}
for $t \in [0,T] $ and $(\mathbf{x},\mathbf{v}) \in  I\!\!R^{N}  \times I\!\!R^{N}$.

When $a=0$, fundamental solutions $\Gamma_{\mathbf F}$ are constructed in reference \cite{victoryclassical} solving the integral equation (F3) in Lemma 2.2:
\begin{eqnarray}
\Gamma_{\mathbf F}(t, {\mathbf x}, {\mathbf v};  \tau,  {\boldsymbol \xi}, {\boldsymbol \nu}) =
G(t, {\mathbf x}, {\mathbf v};  \tau,  {\boldsymbol \xi}, {\boldsymbol \nu}) 
+ \nonumber \\
\int_{\tau}^t \int_{I\!\!R^{N} \times I\!\!R^{N}}
\nabla_{\boldsymbol \nu'} G(t, {\mathbf x}, {\mathbf v};   \tau',  {\boldsymbol \xi}', {\boldsymbol \nu}') 
{\mathbf F}(\tau',  {\boldsymbol \xi}')
\Gamma_{\mathbf F}( \tau',  {\boldsymbol \xi}', {\boldsymbol \nu}';   \tau,  {\boldsymbol \xi}, {\boldsymbol \nu}) 
d{\boldsymbol \xi}' d{\boldsymbol \nu}' d\tau'.
\label{eq:integralF}
\end{eqnarray}

For $a \neq 0$, we correct $\Gamma_{\mathbf F}$ adding the new contribution from $a$. We must solve the integral equation:
\begin{eqnarray}
\Gamma_{{\mathbf F},a}(t, {\mathbf x}, {\mathbf v};  \tau,  {\boldsymbol \xi}, {\boldsymbol \nu}) 
=
\Gamma_{{\mathbf F}}(t, {\mathbf x}, {\mathbf v};  \tau,  {\boldsymbol \xi}, {\boldsymbol \nu})  \nonumber \\
- \int_{\tau}^t \int_{I\!\!R^{N} \times I\!\!R^{N}} \hskip -6mm
\Gamma_{\mathbf F}(t, {\mathbf x}, {\mathbf v};   \tau',  {\boldsymbol \xi}', {\boldsymbol \nu}') a(\tau', {\boldsymbol \xi}', {\boldsymbol \nu}')
\Gamma_{{\mathbf F},a}( \tau',  {\boldsymbol \xi}', {\boldsymbol \nu}';   \tau,  {\boldsymbol \xi}, {\boldsymbol \nu}) 
d{\boldsymbol \xi}' d{\boldsymbol \nu}' d\tau'.
\label{eq:integralFa}
\end{eqnarray}

A fundamental solution is then constructed following standard techniques for uniformly parabolic operators \cite{victoryclassical, friedman, ilin, gema}. We prove the following result:\\

{\bf Proposition 3.1.} {\it If ${\mathbf F}(t,{\mathbf x})$ and $a(t,{\mathbf x}, \mathbf v)$ are measurable and bounded, a generalized fundamental solution $\Gamma_{\mathbf F,a}(t, {\mathbf x}, {\mathbf v};  \tau, {\boldsymbol \xi}, {\boldsymbol \nu})$ of the initial value problem (\ref{eq:linpa})-(\ref{eq:linp0a}) exists. It satisfies the    estimates (F4)-(F5) with coefficients depending on $T$, $\|{\mathbf F}\|_{\infty}$
and $\|a\|_{\infty}$.  

If $p_0 \in L^q(I\!\!R^{N} \times I\!\!R^{N})$
and $f \in L^{\infty}(0,T;L^q(I\!\!R^{N} \times I\!\!R^{N}))$,  $\forall q \in [1, \infty]$,
the integral expression (\ref{solint}), with $\Gamma_{\mathbf F}$  replaced
by $\Gamma_{\mathbf F,a}$, is a weak solution of the initial value 
problem (\ref{eq:linpa})-(\ref{eq:linp0a}). 
This solution $p$ of (\ref{eq:linpa})-(\ref{eq:linp0a}) belongs to
$L^{\infty}(0,T;L^q(I\!\!R^{N} \times I\!\!R^{N}))$ for all $1\leq q \leq \infty$
and is unique.
For nonnegative $p_0$ and $f$, the solution $p$ is also nonnegative.
} \\
 
 {\bf Proof.}
The integral equation (\ref{eq:integralFa}) is solved by a successive approximation procedure:
\begin{eqnarray}
\Gamma_{{\mathbf F},a}^{(0)}(t, {\mathbf x}, {\mathbf v};   \tau,  {\boldsymbol \xi}, {\boldsymbol \nu})=
\Gamma_{{\mathbf F}}(t, {\mathbf x}, {\mathbf v};   \tau, {\boldsymbol \xi}, {\boldsymbol \nu}),
\nonumber \\
\Gamma_{{\mathbf F},a}^{(\ell+1)}(t, {\mathbf x}, {\mathbf v};  \tau,  {\boldsymbol \xi}, {\boldsymbol \nu}) 
= \Gamma_{{\mathbf F}}(t, {\mathbf x}, {\mathbf v};  \tau,  {\boldsymbol \xi}, {\boldsymbol \nu})  \nonumber \\
- \int_{\tau}^t \int_{I\!\!R^{N} \times I\!\!R^{N}} \hskip -4mm
\Gamma_{{\mathbf F}}(t, {\mathbf x}, {\mathbf v};   \tau',  {\boldsymbol \xi}', {\boldsymbol \nu}') 
a(\tau', {\boldsymbol \xi}', {\boldsymbol \nu}')
\Gamma_{{\mathbf F},a}^{(\ell)}( \tau',  {\boldsymbol \xi}', {\boldsymbol \nu}';   \tau,  {\boldsymbol \xi}, {\boldsymbol \nu}) 
d{\boldsymbol \xi}' d{\boldsymbol \nu}' d\tau'.
\label{iterationFa}
\end{eqnarray}
The limiting function is represented as the sum of a telescopic series:
\begin{eqnarray}
\Gamma_{{\mathbf F},a}(t, {\mathbf x}, {\mathbf v};  \tau,  {\boldsymbol \xi}, {\boldsymbol \nu}) =
\Gamma_{{\mathbf F}}(t, {\mathbf x}, {\mathbf v};  \tau,  {\boldsymbol \xi}, {\boldsymbol \nu})
+ \nonumber \\ \sum_{\ell=0}^{\infty} 
\left( \Gamma_{{\mathbf F},a}^{(\ell+1)}(t, {\mathbf x}, {\mathbf v};  \tau,  {\boldsymbol \xi}, {\boldsymbol \nu})
-\Gamma_{{\mathbf F},a}^{(\ell)}(t, {\mathbf x}, {\mathbf v};  \tau,  {\boldsymbol \xi}, {\boldsymbol \nu})\right).
\label{telescopicFa}
\end{eqnarray}

We prove below an estimate of the form:
\begin{eqnarray}
|\Gamma_{{\mathbf F},a}^{(\ell+1)}(t, {\mathbf x}, {\mathbf v};  \tau,  {\boldsymbol \xi}, {\boldsymbol \nu}) -\Gamma_{{\mathbf F},a}^{(\ell)}(t, {\mathbf x}, {\mathbf v};  \tau, {\boldsymbol \xi}, {\boldsymbol \nu})| \leq  \nonumber \\
\| a \|_{\infty}^{\ell+1} {T^{\ell+1} \over (\ell+1)!} \Gamma_{{\mathbf F}}(t, {\mathbf x}, {\mathbf v}; \tau,  {\boldsymbol \xi}, {\boldsymbol \nu}).
\label{telescopicele}
\end{eqnarray}

Thanks to identity (F2) in Lemma 2.2 we find for $\ell=0$:
\begin{eqnarray}
|\Gamma_{{\mathbf F},a}^{(1)}(t, {\mathbf x}, {\mathbf v};  \tau, {\boldsymbol \xi}, {\boldsymbol \nu})
-\Gamma_{{\mathbf F},a}^{(0)}(t, {\mathbf x}, {\mathbf v};  \tau, {\boldsymbol \xi}, {\boldsymbol \nu})|  \leq \nonumber \\
\| a \|_{\infty} \int_{\tau}^t \int_{I\!\!R^{N} \times I\!\!R^{N}}
\Gamma_{{\mathbf F}}(t, {\mathbf x}, {\mathbf v};   \tau',  {\boldsymbol \xi}', {\boldsymbol \nu}') 
\Gamma_{{\mathbf F}}( \tau',  {\boldsymbol \xi}', {\boldsymbol \nu}';   \tau,  {\boldsymbol \xi}, {\boldsymbol \nu}) d{\boldsymbol \xi}' d{\boldsymbol \nu}' d\tau'
\leq  \nonumber\\
  \| a \|_{\infty} (t-\tau)  
\Gamma_{{\mathbf F}}(t, {\mathbf x}, {\mathbf v};   \tau,  {\boldsymbol \xi}, {\boldsymbol \nu}).  \nonumber 
\end{eqnarray}

Assuming that estimate (\ref{telescopicele}) holds up to $\ell$, we check
that it also holds for $\ell+1$:
\begin{eqnarray}
|\Gamma_{{\mathbf F},a}^{(\ell+1)}(t, {\mathbf x}, {\mathbf v};  \tau, {\boldsymbol \xi}, {\boldsymbol \nu})
-\Gamma_{{\mathbf F},a}^{(\ell)}(t, {\mathbf x}, {\mathbf v};  \tau, {\boldsymbol \xi}, {\boldsymbol \nu})| \nonumber \\
\leq \| a \|_{\infty} \int_{\tau}^t \int_{I\!\!R^{N} \times I\!\!R^{N}}
\hskip -8mm
\Gamma_{{\mathbf F}}(t, {\mathbf x}, {\mathbf v};   \tau', {\boldsymbol \xi}', {\boldsymbol \nu}') 
\big|\Gamma^{(\ell)}_{{\mathbf F},a}( \tau', {\boldsymbol \xi}', {\boldsymbol \nu}';   \tau, {\boldsymbol \xi}, {\boldsymbol \nu}) - \nonumber \\
\Gamma^{(\ell-1)}_{{\mathbf F},a}( \tau', {\boldsymbol \xi}', {\boldsymbol \nu}';   \tau, {\boldsymbol \xi}, {\boldsymbol \nu}) \big| d{\boldsymbol \xi}' d{\boldsymbol \nu}' d\tau' \nonumber \\
\leq \| a \|_{\infty}^{\ell+1}  
\int_{\tau}^t {(\tau'-\tau)^{\ell} \over \ell !} \int_{I\!\!R^{N} \times I\!\!R^{N}} 
\hskip -8mm
\Gamma_{{\mathbf F}}(t, {\mathbf x}, {\mathbf v}; \tau', {\boldsymbol \xi'}, {\boldsymbol \nu}')  \Gamma_{{\mathbf F}}(\tau', {\boldsymbol \xi}', {\boldsymbol \nu}';   \tau, {\boldsymbol \xi}, {\boldsymbol \nu}) 
d{\boldsymbol \xi}' d{\boldsymbol \nu}' d\tau' \nonumber \\
\leq \| a \|_{\infty}^{\ell+1} {(t-\tau)^{\ell+1} \over (\ell+1)!}
\Gamma_{{\mathbf F}}(t, {\mathbf x}, {\mathbf v}; \tau, {\boldsymbol \xi}, {\boldsymbol \nu}), \nonumber
\end{eqnarray}
using again identity (F2). This proves estimate (\ref{telescopicele})  for all $\ell$ by induction.

The upper bound $ \| a \|_{\infty}^{\ell+1} {T^{\ell+1} \over (\ell+1)!}$ tends to zero as $\ell$ tends to infinity, which ensures convergence of the iterative scheme. The series $ \sum_{\ell=0}^{\infty} \| a \|_{\infty}^{\ell+1} {T^{\ell+1} \over (\ell+1)!}$ converges, which proves the series expansion (\ref{telescopicFa}). Inserting the bounds (\ref{telescopicele}) in the series, we see that the fundamental solution satisfies:
\begin{eqnarray}
|\Gamma_{{\mathbf F},a}(t, {\mathbf x}, {\mathbf v};  \tau, {\boldsymbol \xi}, {\boldsymbol \nu})| \leq 
\left[ 1 +\sum_{\ell=0}^{\infty} \| a \|_{\infty}^{\ell+1} {T^{\ell+1} \over (\ell+1)!} \right] 
\Gamma_{{\mathbf F}}(t, {\mathbf x}, {\mathbf v}; 
\tau, {\boldsymbol \xi}, {\boldsymbol \nu}).
\label{boundFa}
\end{eqnarray}
Using $(F4)$ from Lemma 2.2, we obtain a similar upper bound in terms of the field free fundamental solution where the constant $M$ now depends  also on $\|a\|_{\infty}$.

Bounds on the derivatives with respect to ${\mathbf v}$ are obtained from a differentiated version of the integral equation, as in references \cite{victoryclassical,gema}.
Bounds on other derivatives 
require higher regularity on $a$ and ${\mathbf F}$. 

The integral expression for the solutions of the initial value problem follows as in reference \cite{victoryclassical} when $a$ and ${\mathbf F}$ are smooth enough.
When they are just bounded functions, we approximate the coefficients by regularized sequences as in references \cite{carpio,gema} and pass to the limit in the approximated problems, combining uniform bounds provided by the integral
expressions and energy inequalities (see Lemma 4.5)
with compactness results for kinetic models (see Lemma 5.2).
The integral expression, combined with the bounds of the fundamental solution in terms of the fundamental solution of the
free field problem (${\mathbf F}$, $a=0$), yields the $L^q$ estimates on the solution of the initial value problem. Uniqueness follows from energy inequalities
\cite{bouchut2}.

The nonnegativity of the fundamental solution can be established in different ways. We may argue that it is the limit of fundamental solutions  for uniformly parabolic problems with an additional $ \delta \Delta_{\mathbf x} p$ term that vanishes as $ \delta \rightarrow 0$, or extend the maximum principle arguments used in reference \cite{friedman} to prove positivity of
fundamental solutions for uniformly parabolic problems, as commented in reference \cite{victoryclassical}. 
This is valid regardless the sign of $a$.   


Let us remark for completeness that nonnegativity of solutions for the initial value problem can be proved without resorting to fundamental solutions provided the data satisfy $p_0\geq 0, f \geq 0$, together with $p_0 \in L^2(I\!\!R^{N} \times I\!\!R^{N})$ and $ f \in L^2([0,T] \times I\!\!R^{N} \times I\!\!R^{N})$, to ensure that $p \in C([0,T],L^2(I\!\!R^{N} \times I\!\!R^{N}))$, $\nabla_{\mathbf v} p \in L^2([0,T] \times I\!\!R^{N} \times I\!\!R^{N})$ and energy identities hold \cite {bouchut2}.   Making the change of variables $v=v'e^{bt}$ and setting $p=e^{ct} q$, $c=Nk$, we find for $q$ the equation:
\begin{eqnarray*}
\frac{\partial}{\partial t} q +
v' e^{bt}  \!\cdot\! \nabla_\mathbf{x}   q
\!+\!  e^{-bt}  {\mathbf F} \!\cdot\!   \nabla_\mathbf{v'}   q
\!-\!  k v'  \nabla_\mathbf{v'} q
\!-\! \sigma e^{-2bt} \Delta_\mathbf{v'} q 
+ a q  = e^{-Nkt} f.
\end{eqnarray*}
Multiplying by $q^-={\rm Max}(-q,0)$ and integrating:
\begin{eqnarray*}
-{1\over 2}\|q^-(t) \|_2^2 
- \int_0^t \sigma e^{-2bs} \| \nabla_\mathbf{v'} q^-(s) \|_2^2 ds
- \int_0^t  \int a(s) q^-(s)^2 ds \\
- k \int_0^t  \|q^-\|_2^2 ds =  \int_0^t \int e^{-Nkt}f q^- ds \geq 0.
\end{eqnarray*}
Therefore $q^-=0$ and the solution $p$ is nonnegative, provided $a\geq 0$. When
$a$ is bounded, a similar result is achieved choosing $c$ large enough in the
change of variables.



%
%

\section{Decay in the velocity space}
\label{sec:estimatesv}

The  nonlinear angiogenesis problem involves the  velocity integrals of the density $\tilde{p}(t,\mathbf{x})= \int_{I\!\! R^{N}} p(t,\mathbf{x},\mathbf{v}) \, d\mathbf{v} $ and  ${j(t,\mathbf{x})}= \int_{I\!\! R^{N}} { |\mathbf{v}|} 
 g(|\mathbf{v}|)  p(t,\mathbf{x},\mathbf{v})\,  d \mathbf{v}$. Unlike Vlasov-Poisson-Fokker-Planck problems, here the force field ${\mathbf F}$ depends on {$j$}, not on $\tilde{p}$.
We will use an iterative scheme to establish existence of solutions for the nonlinear problem. Constructing the iterates and extracting convergent subsequences will require $L^{\infty}$ bounds on $\tilde{p}$ and 
{$j$}. The necessary bounds will follow from  estimates on the velocity decay of $p$ and its derivatives with respect to ${\mathbf v}$.
Let us first state the basic result relating the norms ${\mathbf F}$ and 
{$j$}. \\

{\bf Lemma 4.1.} {\it Let us consider the field ${\mathbf F}$ defined in (\ref{eq:alphaF}), $c$ being a solution of the parabolic problem (\ref{eq:c})-(\ref{eq:pc0}).  
We assume that $c_0\geq 0$, $c_0 \in L^{\infty}(I\!\!R^N)$, $\nabla_{\mathbf x} c_0 \in L^r(I\!\!R^N)$ and ${ j} \in L^{\infty}(0,T; L^q(I\!\!R^N))$, for some $q$, $r$ fulfilling $1 \leq q\leq r \leq \infty$. We furthermore assume that, either ${ j}$ is a bounded function or $c\geq 0$. Then, the $L^{r}_{\mathbf x}$ norm of  ${\mathbf F}$ satisfies:
\begin{eqnarray}
\| {\mathbf F} (t) \|_{L^r_{\mathbf x}} \leq 
d_1 \|\nabla_{\mathbf x} c_0 \|_{L^r_{\mathbf x}}
+ d_1\eta C_{N,r,q} \| c_0\|_{L^\infty_{\mathbf x}} t^{{1\over 2} - {N \over 2}( {1\over q} -{1\over r}) } \| { j} \|_{L^{\infty}(0,t;L^q_{\mathbf x})},
\label{boundFr}
\end{eqnarray}
for $t \in [0,T]$, provided $1/N > (1/q-1/r)$. 
In particular, the $L^{\infty}$ norm of  ${\mathbf F}$ satisfies:
\begin{eqnarray}
\| {\mathbf F} (t) \|_{L^\infty_{\mathbf x}} \leq d_1 \|\nabla_{\mathbf x} c_0 \|_{L^\infty_{\mathbf x}}
+ d_1\eta C_{N,q} \| c_0\|_{L^\infty_{\mathbf x}} t^{{1\over 2} - {N\over 2q}} 
\| { j} \|_{L^{\infty}(0,t;L^q_{\mathbf x})},
\label{boundFinf}
\end{eqnarray}
for $q>N,$ $t \in [0,T]$. 
} \\

{\bf Proof.} 
For ${ j}$ bounded, problem (\ref{eq:c})-(\ref{eq:pc0}) has a unique solution  $c$ that can be expressed in terms of the positive fundamental solution of the linear parabolic operator \cite{aronson,gema,kusuoka} and $c \geq 0$ if $c_0\geq 0$.
Since $c\geq 0$, classical maximum principles for heat equations \cite{friedman} imply that $c$ is bounded from above by the solution of the heat equation obtained setting $\eta=0$ in (\ref{eq:c}), keeping the same initial datum $c_0$.
Using $L^{\infty}$ estimates for heat equations, we conclude that $c \leq \|c_0\|_{L^{\infty}_{\mathbf x}}$. To control the norms of ${\mathbf F}$, we need to estimate the norms of $\nabla_{\mathbf x} c$.

Differentiating the integral expression for $c$:
\begin{eqnarray}
c(t) = K(t)*_{\mathbf  x}  c_0 \!-\! \eta \int_0^t \hskip -2.5mm 
K(t-s)*_{\mathbf  x} c(s) { j}(s) ds, \label{intc}
\end{eqnarray}
we find:
\begin{eqnarray}
\nabla_{\mathbf  x} c(t) = K(t)*_{\mathbf  x} \nabla_{\mathbf  x} c_0 
\!-\! \eta \int_0^t \hskip -2.5mm \nabla_{\mathbf  x} 
K(t-s)*_{\mathbf  x} c(s) { j}(s) ds, \label{intdc}
\end{eqnarray}
where $K$ is the heat kernel. Combining classical estimates on $L^r$ norms of convolutions \cite{brezis}:
\begin{eqnarray}
\|a * b\|_{r} \leq \|a\|_{q'} \|b\|_q, \quad 1/q'+1/q = 1+ 1/r, \; 1 \leq q,q',r \leq \infty,
\label{convolution}
\end{eqnarray}
and the known decay of $L^{q'}$ norms of derivatives of heat kernels \cite{giga}:
\begin{eqnarray}
\|\nabla_{\mathbf  x} K (t)\|_{q'} \leq C_{N,q'} t^{-N/2(1-1/q')-1/2}, 
\quad t>0, 
\label{decayheat}
\end{eqnarray}
we obtain:  
\[
\| {\mathbf F} (t) \|_{L^r_{\mathbf x}} \leq 
d_1 \|\nabla_{\mathbf x} c_0 \|_{L^r_{\mathbf x}}
+  d_1 \eta C_{N,q'} \| c_0\|_{L^\infty_{\mathbf x}} 
\int_0^t (t-s)^{-{1 \over 2} - {N\over 2}(1 - {1\over q'})}
 \| { j}(s) \|_{L^q_{\mathbf x}} ds,
\]
with the restriction $1+{1/r} ={1/q}+{1/q'}$,  provided
$1/2 -  N/2 (1 - 1/q')= 1/2-N/2 (1/q-1/r)>0$, that  is, 
$1/N > (1/q-1/r)$. Estimate (\ref{boundFr}) follows, with $C_{N,r,q}>0$.
When $r=\infty$, this condition becomes $q>N$.
Estimate (\ref{boundFinf}) follows, with $C_{N,q}>0$. 
\\

For general choices of the weight $g$, we can obtain estimates 
on the $L^q_{\mathbf x}$  norms of $\mathbf j$, ${ j}$  and 
$\tilde p$ if, in addition to the $L^{1}_{\mathbf x \mathbf v}$ 
and the $L^\infty_{\mathbf x \mathbf v}$ norms of the density $p$, we control 
the norms of velocity moments $|\mathbf v|^\beta p$, as shown by the next 
Lemma.\\

{\bf Lemma 4.2.} {\it Set ${\mathbf j}(t,{\mathbf x})= \int_{I\!\! R^{N}} 
{\mathbf v}  g(|\mathbf{v}|) p(t,{\mathbf x},{\mathbf v}) d{\mathbf v}$, with $p\geq 0$, $g \geq 0$.  
The following inequalities hold:
\begin{eqnarray}
|{\mathbf j}(t,{\mathbf x})|= |\int_{I \!\! R^{N}} {\mathbf v} 
g(|\mathbf{v}|) p(t,{\mathbf x},{\mathbf v}) 
d {\mathbf v}| \leq {N}\int_{I \!\! R^{N}} |{\mathbf v}| 
g(|\mathbf{v}|) p(t,{\mathbf x},{\mathbf v}) d {\mathbf v},
\label{boundj} \\
 \| |\mathbf v|^{\ell} p \|_{L^1_{\mathbf x \mathbf v}}  \leq 
\| p \|_{L^1_{\mathbf x \mathbf v}}^{1- {\ell \over \beta}} \;
\| |\mathbf v|^\beta p \|_{L^1_{\mathbf x \mathbf v}}^{{\ell \over \beta}}, 
\quad \beta > \ell >0, \label{boundmlp1} \\
\| \int_{I \!\! R^N} |\mathbf v|^\ell p \, d {\mathbf v} 
\|_{L^{N+\beta \over N+\ell}_{\mathbf x}}  \leq C_{N,\beta,\ell} \,
\| p \|_{L^\infty_{\mathbf x \mathbf v}}^{\beta-\ell \over N+\beta} \;
\| |\mathbf v|^{\beta} p \|_{L^1_{\mathbf x \mathbf v}}^{N+\ell \over N+\beta}, 
\quad \beta > \ell >0,  \label{boundmlpinf} \\
\|\int_{I \!\! R^{N}} |{\mathbf v}| p d{\mathbf v} \|_{L^\infty_{\mathbf x}} \leq  
C_{N,\beta} \|p\|_{L^\infty_{\mathbf x \mathbf v}}^{1-{N+1\over \beta}}
\| (1+|{\mathbf v}|^2)^{\beta/2} p\|_{L^\infty_{\mathbf x \mathbf v}}^{N+1\over\beta}, 
\quad \beta>N+1,
\label{boundvinf} \\
\| \int_{I \!\! R^{N}}  p d {\mathbf v} \|_{L^\infty_{\mathbf x}}
\leq  C_{N,\beta}  \| p \|_{L^\infty_{\mathbf x \mathbf v}}^{1-{N \over \beta}}
\|(1+|{\mathbf  v}|^2)^{\beta/2} p \|_{L^\infty_{\mathbf x \mathbf v}}^{N \over\beta}, \quad \beta>N,
\label{boundinf} \\
\|(1+|{\mathbf  v}|^2)^{(\beta-1)/2} p \|_{L^\infty_{\mathbf x \mathbf v}} 
\leq C_\beta \|p\|_{L^\infty_{\mathbf x \mathbf v}}^{1/\beta} 
\|(1+|{\mathbf  v}|^2)^{\beta/2} p \|_{L^\infty_{\mathbf x \mathbf v}}^{1-1/\beta}, \quad \beta>1, 
\label{boundinterp}
\end{eqnarray}
provided the involved integrals and norms are finite. 
} \\

{\bf Proof.} 
Inequality (\ref{boundj}) ensues from 
$\left(\sum_{i=1}^N \! v_i^2\right)^{1\over 2}
\! \! \leq  \! \sum_{i=1}^N \! |v_i|
\! \leq \! N \! \left(\sum_{i=1}^N \! |v_i|^2 \right)^{1\over 2}$, 
thanks to Young's  inequality  
$2 a b \leq {1\over 2} a^2 + {1\over 2} b^2$, $a,b \geq 0$.

Inequality (\ref{boundmlp1}) follows applying  H\"older's inequality \cite{brezis}
to $|\mathbf v|^{\ell} p^{\ell\over \beta} $ and $ p^{(\beta - \ell) \over \beta}$
with $q={\beta \over \ell}$, $q'={\beta\over \beta-\ell}$:
\[
\int_{I \!\! R^{N}\times I \!\! R^{N}} \hskip -8mm|\mathbf v|^{\ell} 
p d {\mathbf v} d\mathbf x
= \int_{I \!\! R^{N}\times I \!\! R^{N}} \hskip -8mm |\mathbf v|^{\ell} p^{\ell\over \beta} 
p^{\beta - \ell \over \beta} d\mathbf v d\mathbf x
\leq  
\||\mathbf v|^{\ell} p^{\ell\over \beta}  \|_{L^{\beta \over \ell}_{\mathbf x \mathbf v}} \;
\| p^{\beta - \ell \over \beta} \|_{L^{\beta\over \beta-\ell}_{\mathbf x \mathbf v}}. 
\]

Inequality (\ref{boundmlpinf}) is proven in reference \cite{bouchut} 
(pp. 246) when $\ell=0$. We prove the general case in a similar way.
Using ${|{\mathbf v}|^{\beta-\ell} \over R^{\beta-\ell}} \geq 1$ for 
$|{\mathbf v}|> R$, we get:
\begin{eqnarray*}
\int_{I \!\! R^N} \hskip -2mm |{\mathbf v}|^{\ell} p(t,{\mathbf x},{\mathbf v})
d {\mathbf v} \leq \|p \|_{\infty} \int_{|{\mathbf v}| \leq R} \hskip -4mm
|{\mathbf v}|^{\ell} d {\mathbf v} + {1 \over R^{\beta-\ell}} \int_{|{\mathbf v}|> R}
\hskip -4mm |{\mathbf v}|^{\beta} p \,d {\mathbf v}  \\ \leq 
C_{N,\ell}  \|p \|_{\infty} R^{N+\ell} + {1 \over R^{\beta-\ell}}  \int_{I \!\! R^{N}} 
\hskip -2mm |{\mathbf v}|^{\beta} p \,d {\mathbf v}, 
\end{eqnarray*}
where $C_{N,\ell}={\rm meas(S_{N-1})} (N+\ell)^{-1}$, ${\rm meas(S_{N-1})}$
being the measure of the unit sphere.
Choosing $R$ to minimize the right hand side, we find:
\[
\int_{I \!\! R^N} \hskip -2mm |{\mathbf v}|^\ell p(t,{\mathbf x},{\mathbf v})
d {\mathbf v}  \leq C_{N,\beta,\ell}
\| p \|_{L^\infty_{\mathbf x \mathbf v}}^{\beta-\ell \over N+\beta} \;
\left(\int_{I \!\! R^N} \hskip -2mm 
|\mathbf v|^{\beta} p \, d {\mathbf v} \right)^{N+\ell \over N+\beta}. 
\]
Estimate (\ref{boundmlpinf}) follows taking the 
$L^{N+\beta\over N+\ell}_{\mathbf x}$ norm.

To prove (\ref{boundvinf}) we observe that:
\begin{eqnarray}
\int_{I \!\! R^{N}} \hskip -4mm |{\mathbf v}| |p(t,{\mathbf x},{\mathbf v})| 
d {\mathbf v} \leq \int_{|{\mathbf v}|\leq R} \hskip -6mm 
|{\mathbf v}| d {\mathbf v} \| p \|_{L^\infty_{\mathbf x \mathbf v}}
\!+\! \int_{|{\mathbf v}|>R}  
{  d {\mathbf v} \over (1+|{\mathbf  v}|^2)^{\delta \over 2}}  
\| (1+|{\mathbf  v}|^2)^{\delta+1 \over 2} p \|_{L^\infty_{\mathbf x \mathbf v}} 
\nonumber \\
\leq  C_{N,\beta} \left[  R^{N+1} 
\| p \|_{L^\infty_{\mathbf x \mathbf v}} + R^{N-\delta}
 \| (1+|{\mathbf  v}|^2)^{\delta+1 \over 2} p \|_{L^\infty_{\mathbf x \mathbf v}} \right]. 
\nonumber 
\end{eqnarray}
Optimizing with respect to $R$, we obtain (\ref{boundvinf}) for $\beta = \delta+1>N+1$.

Inequalities (\ref{boundinf}) and (\ref{boundinterp}) are proven in Ref.  
\cite{degond}, see Lemma B.1 therein.\\

With the choice $g(|\mathbf v|) = 
 [1+ e^{\delta(|\mathbf v|^2-v_{max}^2)}  ]^{-1} $
performed in definition (\ref{eq:intpintvp}), immediate estimates of the
$L^q_{\mathbf x}$ norms of the flux { $j$} in terms of  
norms of $p$ follow. \\

{\bf Lemma 4.3.} {\it For any  $p \geq 0$, the norms 
$\| { j} \|_{L^q_{\mathbf x}}$, $1 \leq q
\leq \infty,$ of the flux defined in equation (\ref{eq:intpintvp}) are bounded
in terms of the $L^1$ and $L^\infty$ norms of the weight, 
$\||\mathbf v|g\|_{L^1_{\mathbf v}}$, $\||\mathbf v|g\|_{L^\infty_{\mathbf v}}$, 
and the function $p$, $\|p \|_{L^1_{\mathbf x \mathbf v}}$,
$\|p \|_{L^\infty_{\mathbf x \mathbf v}}$.
}\\

{\bf Proof.} The weight $|\mathbf v| g(|\mathbf v|) = |\mathbf v| 
[1+ e^{\delta(|\mathbf v|^2-v_{max}^2)}  ]^{-1}$ is bounded and integrable. 
We find:
\begin{eqnarray}
\|{ j} \|_{L^1_{\mathbf x}} \leq 
 \||\mathbf v|g\|_{L^\infty_{\mathbf v}}
\|p\|_{L^1_{\mathbf x \mathbf v}}, \label{j1} \\
\| { j}\|_{L^\infty_{\mathbf x}} \leq 
\||\mathbf v|g\|_{L^1_{\mathbf v}}
\|p \|_{L^\infty_{\mathbf x \mathbf v}}. \label{jinfty} 
\end{eqnarray}
By interpolation \cite{brezis},  all the intermediate $L^q_{\mathbf x}$ norms
are also bounded using those magnitudes. \\

Let us now obtain $L^{\infty}_{\mathbf x}$ estimates for the marginal density $\tilde p$ appearing in the anastomosis term. We adapt a strategy introduced in \cite{degond} (see Lemma 3.1 therein), that exploits decay in the velocity space to obtain spatial bounds of the moments.  \\

{\bf Proposition 4.4.} {\it Let $p \geq 0$ be a solution of 
\begin{eqnarray} 
\frac{\partial}{\partial t} p(t,\mathbf{x},\mathbf{v}) \!+\! 
\mathbf{v} \!\cdot\! \nabla_\mathbf{x}   p(t,\mathbf{x},\mathbf{v}) 
\!+\!  \nabla_\mathbf{v} \!\cdot\! [({\mathbf F}(t,\mathbf{x}) \!-\! k\mathbf{v} ) 
p(t,\mathbf{x},\mathbf{v}) ] \!-\! \sigma
\Delta_\mathbf{v} p(t,\mathbf{x},\mathbf{v}) \nonumber \\
+ a(t,\mathbf{x},\mathbf{v}) p(t,\mathbf{x},\mathbf{v}) = 0,  
\label{eq:linpaj} \\
p(0,\mathbf{x},\mathbf{v})=p_0(\mathbf{x},\mathbf{v}),
\label{eq:linp0aj}
\end{eqnarray}
under the hypotheses:
\begin{itemize}
\item [(i)] $a \in L^{\infty}((0,T)\times I\!\!R^{N} \times I\!\!R^{N})$, 
\, ${\bf F} \in L^{\infty}((0,T)\times I\!\!R^{N})$,
\item[(ii)] $ (1+|{\mathbf  v}|^2)^{\beta/2} p_0 \in L^\infty(I\!\!R^{N} \times I\!\!R^{N})$, \; $\beta >N$,
\item[(iii)]
$(1+|{\mathbf v}|^2)^{\beta/2} p_0 \in
L^1(I\!\!R^{N} \times I\!\!R^{N})$,  \; $p_0\geq 0.$
\end{itemize}
Then,  
$\| (1+|{\mathbf v}|^{2})^{\beta/2} p \|_{L^\infty(0,T;L^{\infty}_{\mathbf x \mathbf v})
}$
and $\| p\|_{L^\infty(0,T;L^\infty_{\mathbf x}L^1_{\mathbf v})}$ 
are bounded by constants depending on  the parameters
$\sigma$,  $k$, $T$,  $\beta$, $N$, and the norms
$\| (1+|{\mathbf v}|^{2})^{\beta/2} p_0 \|_{L^{\infty}_{\mathbf x \mathbf v
}}$,   $\|{\bf F}\|_{L^\infty(0,T;L^{\infty}_{\mathbf x})}$,
 $\|a^-\|_{L^\infty(0,T;L^{\infty}_{\mathbf x \mathbf v})}$. If $\beta > 
N+ 1$, $\| |\mathbf v|p\|_{L^\infty(0,T;L^\infty_{\mathbf x}L^1_{\mathbf v})}$ 
is similarly bounded.

Assume further:
\begin{itemize}
\item[(iv)]  ${\bf F}(c(t,{\mathbf x}))$ is given by 
(\ref{eq:alphaF}), $c$, $\nabla_{\mathbf x} c$ are given by (\ref{intc})-(\ref{intdc}) 
with $c_0 \in W^{1,\infty}_{\mathbf x}$, $c\geq 0$,
$c$ is coupled to (\ref{eq:linpaj}) by (\ref{eq:intpintvp}), and
$a= \gamma \int_0^t \tilde p \, ds -\alpha(c) \rho$.
\end{itemize}
Then, similar bounds hold with constants depending on the  parameters 
$\sigma$, $k$, $d$, $\eta$,  
$d_1$, $\alpha_1$, $T$, $\beta$, $N$, and the norms
$\| (1+|{\mathbf v}|^{2})^{\beta/2} p_0 \|_{L^{\infty}_{\mathbf x \mathbf v
}}$, $\|c_0\|_{W^{1,\infty}_{\mathbf x}}$, 
$\|\rho \|_{L^{\infty}_{\mathbf v}}$, 
$\||\mathbf v| g\|_{L^\infty_{\mathbf v}\cap L^1_{\mathbf v}}$,}
$\|p\|_{L^\infty(0,T;L^{\infty}_{\mathbf x \mathbf v}\cap
L^{1}_{\mathbf x \mathbf v})}$.\\


{\bf Proof.} Using expressions in terms of fundamental solutions
$(1+|\mathbf v|^2)^{\beta/2} p \in C([0,T]; L^1_{\mathbf x} \cap L^{\infty}_{\mathbf v})$
\cite{bouchut2,bouchut3}.
We set $Y(t,{\mathbf x},{\mathbf v})=(1+|{\mathbf v}|^2)^{\beta/2} 
p(t,{\mathbf x},{\mathbf v})$.
Multiplying  equation (\ref{eq:linpaj}) by $(1+|{\mathbf v}|^2)^{\beta/2}$, 
$\beta>0$, we get:
\begin{eqnarray}
{\partial \over \partial t} Y \!+\! {\mathbf v} \nabla_{\mathbf x} Y \!+\! \Big( {\mathbf F} 
+ 2 \sigma \beta {{\mathbf v} \over 1 + |{\mathbf v}|^2} \!-\! k {\mathbf v} \Big) 
\nabla_{\mathbf v} Y \!-\! \Delta_{\mathbf v} Y \!=\! (N k\!-\!a) Y 
\!+\! R_1 \!+\! R_2 \!+\! R_3,
\nonumber 
\end{eqnarray}
where
\begin{eqnarray}
R_1 = 
{\beta} (1+|{\mathbf v}|^2)^{\beta/2 -1}  {\mathbf F} \cdot {\mathbf v} p,  \nonumber \\
R_2 = -k {\beta/2} (1+|{\mathbf v}|^2)^{\beta/2 -1} 2 |{\mathbf  v}|^2 p 
= -k \beta {|{\mathbf  v}|^2 \over (1+|{\mathbf v}|^2)} Y,  \nonumber \\
R_3= \sigma \beta (\beta +2)   {|{\mathbf v}|^2 \over (1+|{\mathbf v}|^2)^2} Y - 
N \sigma \beta {1 \over 1+|{\mathbf v}|^2 } Y.
\nonumber
\end{eqnarray}
Thanks to the $L^{\infty}$ estimate in reference \cite{degond} for this type of operators (see Proposition A.3 therein) we get:
\begin{eqnarray}
\| Y(t) \|_{L^\infty_{\mathbf x \mathbf v}} 
\!\!\leq \!\! \| Y(0) \|_{L^\infty_{\mathbf x \mathbf v}} 
\!\!+\!\!\! \int_0^t \hskip -2mm\Big( \!\!
(N k\!+\!\|a^-\|_{\infty} \!) \| Y \|_{L^\infty_{\mathbf x \mathbf v}} 
\!\!+\!\! \|R_1 \|_{L^\infty_{\mathbf x \mathbf v}} 
\!\!+\!\! \|R_2 \|_{L^\infty_{\mathbf x \mathbf v}} 
\!\!+\!\! \|R_3 \|_{L^\infty_{\mathbf x \mathbf v}} \!\! \Big) ds. \nonumber
\end{eqnarray}
The factors ${|{\mathbf v}|^\varepsilon \over 1+|{\mathbf v}|^2} \leq 1$, for
$0 \leq \varepsilon \leq 2$. Therefore,
\begin{eqnarray}
 \|R_1\|_{L^\infty_{\mathbf x \mathbf v}} \leq  
 \beta \|(1+|{\mathbf v}|^2)^{\beta/2 -1}  {\mathbf F} \cdot {\mathbf v} p 
 \|_{L^\infty_{\mathbf x \mathbf v}}, \nonumber \\
 \|R_2\|_{L^\infty_{\mathbf x \mathbf v}} \leq   
  k \beta \|Y \|_{L^\infty_{\mathbf x \mathbf v}}, \nonumber \\
 \|R_3\|_{L^\infty_{\mathbf x \mathbf v}} \leq   \sigma \beta (\beta + 2 + N) 
 \|Y \|_{L^\infty_{\mathbf x \mathbf v}}. \nonumber
\end{eqnarray}
The key term to be bounded is 
$\|(1+|{\mathbf v}|^2)^{\beta/2 -1}  {\mathbf F} \cdot {\mathbf v} p \|_{L^\infty_{\mathbf x \mathbf v}}$.

Under assumptions (i), (ii) and (iii), we may resort to interpolation inequality (\ref{boundinterp}) or just set:
\begin{eqnarray*}
\|(1+|{\mathbf v}|^2)^{\beta/2 -1}  
{\mathbf F} \cdot {\mathbf v} p \|_{L^\infty_{\mathbf x \mathbf v}}
\leq   { N |{\mathbf v}| \over 1+|{\mathbf v}|^2}
\|{\bf F}\|_{\infty} \|Y\|_{L^\infty_{\mathbf x \mathbf v}}
\leq N\|{\bf F}\|_{\infty} \|Y\|_{L^\infty_{\mathbf x \mathbf v}}.
\label{boundF0}
\end{eqnarray*}
Taking 
$A= (N\|{\bf F}\|_{\infty}+k) \beta 
+ \sigma \beta (\beta + 2 + N) + Nk + \|a^-\|_{\infty}$ 
and  $B = \| Y(0) \|_{\infty}$, 
Gronwall's inequality implies
\begin{eqnarray}
\| Y(t) \|_{L^\infty_{\mathbf x \mathbf v}} \leq B e^{A t}, \quad t \in [0,T].
\nonumber 
\end{eqnarray}

When (iv) holds, we are in fact dealing with the nonlinear system.
Lemma 4.1 yields estimate (\ref{boundFinf}) for $\|\mathbf F \|_{\infty}$ in terms 
of  $\| { j} \|_{L^q_{\mathbf x}}$ for $q>N,$
$t \in [0,T]$.   Lemma 4.3 shows that these norms can be estimated
from the norms of the weight $|\mathbf v|g(\mathbf v)$ and the solution  $p$.
Therefore, 
\begin{eqnarray*}
\|Ê\mathbf F\|_{\infty} \leq C(d_1,  \eta, d,  N, \|c_0\|_{W^{1,\infty}},  
\||\mathbf v| g\|_{L^\infty_{\mathbf v}\cap L^1_{\mathbf v}},
\|p\|_{L^\infty(0,T;L^{\infty}_{\mathbf x \mathbf v}\cap
L^{1}_{\mathbf x \mathbf v})}).
\end{eqnarray*}
In this case, $\|a^-Ê\|_{\infty}\leq  \alpha_1 \|\rho\|_{L^\infty_{\mathbf v}}$. 
Therefore, the estimate on $Y$ and its dependence
on parameters and norms is a consequence of the previous statement
for linear problems with general coefficients $\mathbf F $ and $a^{-}$.

Once the velocity decay has been established, the $L^{\infty}$ bounds on 
$ \int_{I\!\! R^{N}} p d {\mathbf v}$ follow from inequality
(\ref{boundinf}) in Lemma 4.2. When $\beta > N+1$, the estimate
on $ \int_{I\!\! R^{N}} |\mathbf v| p d {\mathbf v}$ follows from
inequality (\ref{boundvinf}).
\\

Let us finally recall a Lemma from \cite{bouchut2} (see Proposition 2.2 therein), useful to estimate derivatives with respect to ${\mathbf v}$. \\

{\bf Lemma 4.5.} {\it Let ${\mathbf F}(t,{\mathbf x}) \in L^{\infty}((0,T)\times I\!\! 
R^{N})$, 
$a(t,{\mathbf x}) \in L^{\infty}((0,T)\times I\!\! R^N)$, \ 
$p_0 \in L^2(I\!\! R^N \times I\!\! R^N)$ 
and $f \in L^2((0,T) \times I\!\! R^N \times I\!\! R^N))$ for $T>0$. Set $\beta
=2k + 2\| a^-\|_{\infty}+ 1$, $a^-$ being the negative  part of $a$.
Then, the solution  $p \in C([0,T]; L^2(I\!\! R^N \times I\!\! R^N))$ of 
(\ref{eq:linpa})-(\ref{eq:linp0a})  satisfies:
\begin{itemize}
\item For any $t\in [0,T]$
\begin{eqnarray}
\int_{I\!\!R^{N} \times I\!\!R^{N}} \hskip -9mm 
p^2 (t) \, d\mathbf x \, d \mathbf v = 
\int_{I\!\!R^{N} \times I\!\!R^{N}} \hskip -9mm 
p^2 (0) \, d\mathbf x \, d \mathbf v  
+ \!\int_0^t \!\! \int_{I\!\!R^{N} \times I\!\!R^{N}} 
\hskip -9mm (2k+2a) p^2 
 \, ds \, d\mathbf x \, d \mathbf v
\nonumber \\
\!- 2 \sigma \!\! \int_0^t \!\! \int_{I\!\!R^{N} \times I\!\!R^{N}}
\hskip -9mm |\nabla_{\mathbf v} p|^2 \, ds \, d\mathbf x \, d \mathbf v
+  2 \int_0^t \!\! \int_{I\!\!R^{N} \times I\!\!R^{N}} 
\hskip -9mm p f \, ds \, d\mathbf x \, d \mathbf v ,
\label{energyidentity} \\
\int_{I\!\!R^{N} \times I\!\!R^{N}} \hskip -9mm 
p^2 (t)  \, d\mathbf x \, d \mathbf v
\leq \left( \int_{I\!\!R^{N} \times I\!\!R^{N}} \hskip -9mm 
p^2 (0)  \, d\mathbf x \, d \mathbf v
+  \!\!  \int_0^T \hskip -2mm \int_{I\!\!R^{N} \times I\!\!R^{N}}
\hskip -9mm f^2  \, ds \, d\mathbf x \, d \mathbf v \right) 
e^{\beta t}, \label{energyL2}
\end{eqnarray}
\item $\nabla_{\mathbf v} p \in L^2((0,T)\times 
I\!\! R^{N} \times I\!\! R^{N})$ and
\begin{equation}
\sigma \!\! \int_0^T \!\!\!\! \int_{I\!\!R^{N} \times I\!\!R^{N}} 
\hskip -9mm |\nabla_{\mathbf v} p|^2 \, ds \, d\mathbf x \, d \mathbf v  \!\leq\!
\left( \!
\int_{I\!\!R^{N} \times I\!\!R^{N}} \hskip -9mm 
p^2 (0) \, d\mathbf x \, d \mathbf v 
+  \!\! \int_0^T \!\!\!\! \int_{I\!\!R^{N} \times I\!\!R^{N}} 
\hskip -9mm f^2 \, ds \, d\mathbf x \, d \mathbf v  \! \right) e^{\beta T}.
\label{energyH1}
\end{equation}
\end{itemize}
}

\section{Iterative scheme}
\label{sec:iterative}

Existence of solutions for the nonlinear angiogenesis model will be proved by means 
of an iterative scheme. For $m\geq 2$ and $t\in[0,T]$, we consider the system
\begin{eqnarray} 
\frac{\partial}{\partial t} p_m(t,\mathbf{x},\mathbf{v}) \!+\! 
\mathbf{v} \!\cdot\! \nabla_\mathbf{x}   p_m(t,\mathbf{x},\mathbf{v}) 
\!+\!  \nabla_\mathbf{v} \!\cdot\! [({\mathbf F}(c_{m-1}(t,\mathbf{x})) 
\!-\! k\mathbf{v} ) p_m(t,\mathbf{x},\mathbf{v}) ]  
 \label{eq:pm} \\
\!-\! \sigma \Delta_\mathbf{v} p_m(t,\mathbf{x},\mathbf{v})
\!+\! \gamma a_{m-1}(t,\mathbf{x})p_{m}(t,\mathbf{x},\mathbf{v}) 
\!=\! \alpha(c_{m-1}(t,\mathbf{x})) \rho(\mathbf{v}) p_{m}(t,\mathbf{x},\mathbf{v}),  
\nonumber  \\
a_{m-1}(t,\mathbf{x}) = \int_0^t \! \!ds \! \! \int_{I\!\! R^{N}} \! \! d{\bf v}' p_{m-1}(s,\mathbf{x},\mathbf{v}'),   \label{eq:am} \\
 \alpha(c_{m-1})=\alpha_1\frac{c_{m-1}}{c_R+ c_{m-1}}, \quad
 {\bf F}(c_{m-1})= \frac{d_1}{(1+\gamma_1c_{m-1})^{q_1}}\nabla_{\mathbf x} c_{m-1},
\label{eq:alphaFm} \\
p_m(0,\mathbf{x},\mathbf{v}) = p_0(\mathbf{x},\mathbf{v}),  \label{eq:pm0} \\
\frac{\partial}{\partial t}c_{m-1}(t,\mathbf{x}) = d \Delta_{\mathbf x} c_{m-1}(t,\mathbf{x})  - \eta c_{m-1}(t,\mathbf{x})
{ j_{m-1}(t,\mathbf{x})}, \label{eq:cm} \\
{ j_{m-1}(t,\mathbf{x})} \!=\! \int_{I\!\! R^{N}} 
{ { |\mathbf{v}'| }
\over  1+ e^{\delta (|{\mathbf v}'|^2 - v_{max}^2)} }
 p_{m-1}(t,\mathbf{x},\mathbf{v}') 
d{\bf v}'  \label{eq:jm}, \\
c_{m-1}(0,\mathbf{x}) = c_0(\mathbf{x}).   \label{eq:cm0}
\end{eqnarray}
Alternative schemes replacing (\ref{eq:pm}) with
\begin{eqnarray}
\frac{\partial}{\partial t} p_m(t,\mathbf{x},\mathbf{v}) \!+\! 
\mathbf{v} \!\cdot\! \nabla_\mathbf{x}   p_m(t,\mathbf{x},\mathbf{v}) 
\!+\!  \nabla_\mathbf{v} \!\cdot\! [({\mathbf F}(c_{m-1}(t,\mathbf{x})) 
\!-\! k\mathbf{v} ) p_m(t,\mathbf{x},\mathbf{v}) ]  
\label{eq:pmbis} \\
\!-\! \sigma \Delta_\mathbf{v} p_m(t,\mathbf{x},\mathbf{v})
\!+\! \gamma a_{m-1}(t,\mathbf{x})p_{m}(t,\mathbf{x},\mathbf{v}) 
\!=\! \alpha(c_{m-1}(t,\mathbf{x})) \rho(\mathbf{v}) p_{m-1}(t,\mathbf{x},\mathbf{v})   \nonumber
\end{eqnarray}
can also be used to generate solutions.

We initialize the scheme with $p_1=0$ and ${ j}_1=0$. 
$c_1$ is the solution of the corresponding heat equation with the appropriate behavior at infinity.  $p_2$ is a classical nonnegative solution of a  Fokker-Planck problem with smooth and bounded coefficient fields ${\mathbf F(c_1)}$ and $\alpha(c_1)$. We will show that the resulting sequence is well defined under our hypotheses on the data  and a subsequence converges to a solution of the original problem. \\

{\bf Theorem 5.1.} {\it Let us assume that:
\begin{eqnarray}
p_0 \geq 0, c_0 \geq 0, \label{signo} \\
c_0  \in L^{\infty}(I\!\!R^{N}), 
\nabla_{\mathbf x} c_0 \in
L^{\infty}(I\!\!R^{N})\cap L^2(I\!\!R^{N}),  \label{c0inf}  \\
(1+|{\mathbf v}|^2)^{\beta/2} p_0 \in L^{\infty}(I\!\!R^{N} \! \times \! I\!\!R^{N}), 
\quad   \beta >N,
\label{p0inf}  \\ 
(1+|{\mathbf v}|^2)^{\beta/2} p_0\in L^1(I\!\!R^{N}\times I\!\!R^{N}), \quad  \beta>N.  \label{p0mom}  
\end{eqnarray}
Then, there exists a nonnegative solution $(p, c)$ of (\ref{eq:p})-(\ref{eq:intpintvp}) satisfying:
\begin{eqnarray}
c \in L^{\infty}(0,T;L^{\infty}(I\!\!R^{N})), 
\nabla_{\mathbf x} c \in L^{\infty}(0,T;L^{\infty}\cap L^2(I\!\!R^{N})), \label{dcinf} \\
p \in L^{\infty}(0,T;L^{\infty}\cap L^1(I\!\!R^{N}\times I\!\!R^{N})),
\nabla_{\mathbf v }p \in L^{2}(0,T;L^2(I\!\!R^{N}\times I\!\!R^{N})),
\label{pinf} \\
(1+|{\mathbf v}|^2)^{\beta/2}p \in L^{\infty}(0,T;L^{\infty}(I\!\!R^{N} \!\times I\!\!R^{N})), 
\label{pdecay} \\
(1+|{\mathbf v}|^2)^{\beta/2} p \in L^{\infty}(0,T;L^1(I\!\!R^{N}\times I\!\!R^{N})), \label{pmom1}  \\
p   \in L^{\infty}(0,T;L^\infty_{\mathbf x}(I\!\!R^{N},
L^1_{\mathbf v}( I\!\!R^{N})), \label{pmominf}  
\end{eqnarray}
with norms bounded in terms of the norms of the data.  \\
If $\nabla_{\mathbf v} p_0 \in  L^\infty_{\mathbf x}(I\!\!R^N,
L^1_{\mathbf v}( I\!\!R^N)) $, then
$\nabla_{\mathbf v} p \in  L^{\infty}(0,T;L^\infty_{\mathbf x}(I\!\!R^N,
L^1_{\mathbf v}( I\!\!R^N)) $ and the solution is unique.
}\\

We prove this result below. The proof in organized in steps. First, we argue that the scheme is well defined using fundamental solutions of Fokker-Planck and parabolic operators. 
Next, we obtain the pertinent uniform estimates on the $L^q$ norms of solutions  
of the iterative scheme, relying mostly on comparison principles. 
The derivatives of $p_m$
with respect to ${\mathbf v}$ are bounded uniformly using energy inequalities.
Then, we use the velocity decay to obtain uniform estimates of the velocity integrals of $p_m$. In this way, we bound the norms of $\tilde p_m$, ${ j_m}$,  $\nabla_{\mathbf x} c_m$ and $\mathbf F(c_m)$. 
Finally, we pass to the limit in the equations using compactness results specific of Fokker-Planck operators. We obtain a nonnegative solution of the nonlinear problem with the stated regularity. To conclude, we establish the uniqueness result using expressions in terms of fundamental solutions and differential inequalities.  We work with scheme
(\ref{eq:pm})-(\ref{eq:cm0}) but a similar proof stands for  
(\ref{eq:pmbis}) and (\ref{eq:am})-(\ref{eq:cm0}).
\\

\noindent {\bf Proof.}\\
\leftline{\it Step 1: Existence of nonnegative solutions for the scheme}

Fundamental solutions of linear parabolic and Fokker-Planck
problems allow us to prove existence, nonnegativity and basic regularity of solutions for (\ref{eq:pm})-(\ref{eq:cm0}) in $[0,T]$.

Indeed, provided the coefficient ${ j_{m-1}} \in 
L^{\infty}((0,T) \times I\!\!R^{N})$, (\ref{eq:cm}) has a unique solution $c_{m-1}$ that can be expressed in terms of the fundamental solution of the linear parabolic operator \cite{aronson,kusuoka,gema}. In view of the
hypotheses on the initial datum, $c_{m-1}\geq 0$.
For initial data $c_0 \in L^{\infty}(I\!\!R^{N})$, 
the solution $c_{m-1} \in L^{\infty}((0,T) \times I\!\!R^{N})$. Therefore, $\alpha(c_{m-1})$ is a bounded function. Evaluating 
the field ${\mathbf F}(c_{m-1})$  involves $\nabla_{\mathbf x} c_{m-1}.$  $c_{m-1}$ can be seen as a solution of a heat equation with a source. 
Differentiating the resulting integral version of (\ref{eq:cm}) in terms of heat kernels,  $\nabla_{\mathbf x} c_{m-1}\in  L^{\infty}((0,T) \times I\!\!R^{N})$   if $\nabla_{\mathbf x} c_0 \in L^{\infty}(I\!\!R^{N})$ and  ${ j_{m-1}} \in L^{\infty}((0,T) \times I\!\!R^{N}).$ Therefore, ${\mathbf F}(c_{m-1})$ is a 
bounded function.

Provided that $a_{m-1}$ is bounded, (\ref{eq:pm}) has a unique nonnegative solution $p_m$, that  can be constructed in terms of fundamental solutions of the Fokker-Planck operator, as argued in Section \ref{sec:potential}.  
Since $p_0 \in L^q(I\!\!R^{N} \times I\!\!R^{N})$,
the solutions $p_m$ belong to 
$L^{\infty}(0,T;L^q(I\!\!R^{N} \times I\!\!R^{N}))$ 
for $1 \leq q \leq \infty$.  
Thanks to Proposition 4.4 and to Lemma 4.3, the coefficients 
$a_{m}$ and ${ j_{m}}$ are both bounded functions if  
$(1+|{\mathbf v}|^2)^{\beta/2}  p_{0} \in L^{\infty}(I\!\!R^{N} 
\times I\!\!R^{N})$ for $\beta >N$.

Therefore, we may go ahead and construct $c_m$ and $p_{m+1}$. We may consistently construct our sequence of iterates. \\

\leftline{\it Step 2: A priori estimates on the angiogenic factor $c_m$}

Uniform bounds on the $L^q$ norms of $c_m$ are obtained thanks to classical 
maximum principles for heat equations \cite{friedman,gema}.   Indeed, the 
solution of (\ref{eq:cm})-(\ref{eq:cm0}) is bounded from above by the solution of:
\begin{eqnarray}
\frac{\partial}{\partial t}{\cal C}(t,\mathbf{x}) = d \Delta_{\mathbf x} {\cal C}(t,\mathbf{x}),
\quad  
{\cal C}(0,\mathbf{x}) =c_0(\mathbf{x}),   \label{eq:Cheat}
\end{eqnarray}
since $c_{m-1}\geq 0$. Therefore, for $t \in [0,T]$
\begin{eqnarray}
0 \leq c_{m-1}(t,\mathbf{x}) \leq {\cal C}(t,\mathbf{x}) \leq \|c_0\|_{\infty}.
\label{cminf} 
\end{eqnarray}

For compactness, we will need estimates on the derivatives.
Notice that $c_{m-1}$ is the solution of a heat
equation with source $- \eta c_{m-1}(t,\mathbf{x})
{ j}_{m-1}(t,\mathbf{x})$.
To estimate the derivatives we can differentiate the integral expression:
\begin{eqnarray}
c_{m-1}(t) = K(t)*c_0 - \eta \int_0^t K(t-s)*c_{m-1}(s) { j}_{m-1}(s) ds, 
\label{cmint} 
\end{eqnarray}
to obtain:
\begin{eqnarray}
\partial_{x_i}c_{m-1}(t) = K(t)*\partial_{x_i}c_0 
- \eta \int_0^t \partial_{x_i} K(t-s)*c_{m-1}(s) { j}_{m-1}(s) ds,
\label{dcmint}
\end{eqnarray}
where $K$ is the heat kernel,   whose derivatives satisfy
$\| \partial_{x_i} K(t) \|_1 \leq M t^{-1/2}$, whereas $\|K(t)\|_1 \leq 1.$
This yields:
\begin{eqnarray}
\| \partial_{x_i}c_{m-1}(t) \|_{L^{\infty}_{\mathbf x}} \!\leq\!  
\|\partial_{x_i}c_0\|_{L^{\infty}_{\mathbf x}} \!+\! \eta M
\!\int_0^t \!\!\! (t-s)^{-1/2} \|c_{m-1}(s) { j}_{m-1}(s) \|_{L^{\infty}_{\mathbf x}} ds
\nonumber \\
\leq
\|\partial_{x_i}c_0\|_{L^\infty_{\mathbf x}} 
+ 2 \eta M t^{1/2} \|c_{m-1} { j}_{m-1}
\|_{L^{\infty}_t L^{\infty}_{\mathbf x}},
\quad t \in [0,T].
\label{dcminf}
\end{eqnarray}
This inequality would bound uniformly the spatial derivatives of $c_{m-1}$ 
provided the spatial derivatives of the initial datum belong to 
$L^{\infty}_{\mathbf x}$ and the sequence $\|  { j}_{m-1} \|_{L^{\infty}_t L^{\infty}_{\mathbf x}}$ is bounded. In this way we would ensure that $\| {\mathbf F}(c_{m-1}) \|_{L^{\infty}_t L^{\infty}_{\mathbf x}}$ is uniformly bounded.  We will come back to this point in  Step 4. 

Writing $c_{m-1}={\cal C}+\tilde{c}_{m-1}$, $\tilde{c}_{m-1}$ is a solution of:
\begin{eqnarray}
\frac{\partial}{\partial t}\tilde{c}_{m-1}(t,\mathbf{x}) 
- d \Delta_{\mathbf x} \tilde{c}_{m-1}(t,\mathbf{x}) 
= - \eta {c}_{m-1}(t,\mathbf{x}){ j}_{m-1}(t,\mathbf{x}), 
\label{eq:tildecm} \\
\tilde{c}_{m-1}(0,\mathbf{x}) =0. \label{eq:tildecm0} 
\end{eqnarray}
Inserting the inequality $\|K(t)\|_{1} \leq 1$
in the integral equation yields  the estimates: 
\begin{eqnarray}
\| \tilde{c}_{m-1}(t) \|_{L^{2}_{\mathbf x}} \leq   
\eta \int_0^t \|K(t-s)\|_{L^{1}_{\mathbf x}} \|c_{m-1}(s)\|_{L^{\infty}_{\mathbf x}} \|  { j}_{m-1}(s) \|_{L^{2}_{\mathbf x}} ds
\nonumber \\
\leq
 \eta t \|c_{m-1} \|_{L^\infty_t L^{\infty}_{\mathbf x}} 
 \|  { j}_{m-1} \|_{L^{\infty}_t L^2_{\mathbf x}}, \quad t \in [0,T],
\label{tildecmq}
\end{eqnarray}
that are uniform provided the sequence $\| { j}_{m-1} \|_{L^{\infty}_t L^2_x}$ is bounded. Additionally, 
\begin{eqnarray}
\| \partial_{x_i}\tilde c_{m-1}(t) \|_{L^{2}_{\mathbf x}} \!\leq\!  
2 \eta M t^{1/2} \|c_{m-1}\|_{L^{\infty}_t L^{\infty}_{\mathbf x}}
\|{ j}_{m-1}\|_{L^{\infty}_t L^2_{{\mathbf x}}}, \quad t \in [0,T],
\label{dcm2}
\end{eqnarray}
are uniformly bounded on the same condition. \\

\leftline{\it Step 3: A priori estimates on the vessel density $p_m$}

Uniform bounds on the $L^q$ norms of  $p_m$ are obtained thanks to the comparison Lemma 2.5, see also Proposition 3.1. 
Notice that $a_{m-1}p_m \geq 0$. Then, Lemma 2.5 implies that the solution 
$p_m$ of (\ref{eq:pm})-(\ref{eq:pm0}) is bounded from above by the solution of:
\begin{eqnarray} 
\frac{\partial}{\partial t} P_m(t,\mathbf{x},\mathbf{v}) \!+\! 
\mathbf{v} \!\cdot\! \nabla_\mathbf{x}   P_m(t,\mathbf{x},\mathbf{v}) 
\!+\!  \nabla_\mathbf{v} \!\cdot\! [({\mathbf F}(c_{m-1}(t,\mathbf{x})) 
\!-\! k\mathbf{v} ) P_m(t,\mathbf{x},\mathbf{v}) ] 
 \label{eq:Pm} \\  
\!-\! \sigma \Delta_\mathbf{v} P_m(t,\mathbf{x},\mathbf{v})  = 
\alpha(c_{m-1}(t,\mathbf{x})) \rho(\mathbf{v}) p_{m}(t,\mathbf{x},\mathbf{v}),  \nonumber  \\
P_m(0,\mathbf{x},\mathbf{v}) = p_0(\mathbf{x},\mathbf{v}).  
\label{eq:Pm0} 
\end{eqnarray}

Therefore, for $t \in [0,T]$
\begin{eqnarray}
0 \leq p_{m}(t,\mathbf{x}) \leq P_m(t,\mathbf{x}) \leq \|P_m\|_{L^{\infty}_{\mathbf x
\mathbf v}},
\label{pminf} \\
\|p_m(t)\|_{L^{q}_{\mathbf x \mathbf v}} \leq  \|P_m\|_{L^{q}_{\mathbf x
\mathbf v}}, \quad 1 \leq q \leq \infty. \nonumber 
\end{eqnarray}
The function $P_m$ admits the integral expression:
\begin{eqnarray}
P_m(t,{\mathbf x},{\mathbf v})= \int_{I\!\!R^{N} \times I\!\!R^{N}}
\hskip -8mm
\Gamma_{{\mathbf F}(c_{m-1})}(t, {\mathbf x}, {\mathbf v};  
0, {\boldsymbol \xi}, {\boldsymbol \nu})  
p_0({\boldsymbol \xi}, {\boldsymbol \nu}) d{\boldsymbol \xi} d{\boldsymbol \nu} 
+ \nonumber \\
\int_0^t  \hskip -2mm \int_{I\!\!R^{N} \times I\!\!R^{N}} \hskip -9mm
\Gamma_{{\mathbf F}(c_{m-1})}(t, {\mathbf x}, {\mathbf v};  \tau, {\boldsymbol \xi}, {\boldsymbol \nu})  \alpha(c_{m-1}(t,{\boldsymbol \xi})) \rho({\boldsymbol \nu}) p_{m}(t,{\boldsymbol \xi},{\boldsymbol \nu}) d{\boldsymbol \xi} d{\boldsymbol \nu} d \tau.
\label{solintPm}
\end{eqnarray}

Thanks to Lemma 2.3 in Section 2, we know that:
\begin{eqnarray}
\|p_m(t)\|_{L^{1}_{\mathbf x\mathbf v}} \leq 
\|P_m(t)\|_{L^{1}_{\mathbf x\mathbf v}} \leq 
\|p_0\|_{L^{1}_{\mathbf x\mathbf v}} \!+ \!
 \alpha_1 \|\rho\|_{L^{\infty}_{\mathbf v}} \!
\int_0^t \|p_{m}(\tau)\|_{L^{1}_{\mathbf x\mathbf v}} d\tau, 
\nonumber 
\\
\|p_m(t)\|_{L^{\infty}_{\mathbf x\mathbf v}} \!\leq\! 
\|P_m(t)\|_{L^{\infty}_{\mathbf x\mathbf v}} \!\leq\! 
e^{N kt}\|p_0\|_{L^{\infty}_{\mathbf x\mathbf v}}
\!+\! \alpha_1 \|\rho\|_{L^{\infty}_{\mathbf v}} \! \int_0^t \!\!\! e^{N k(t-\tau)}
\|p_{m}(\tau)\|_{L^{\infty}_{\mathbf x\mathbf v}} d\tau, 
\nonumber 
\\
\|p_m\|_{L^{q}_{\mathbf x\mathbf v}} \leq \|P_m\|_{L^{q}_{\mathbf x\mathbf v}} \leq \|P_m\|_{L^{1}_{\mathbf x\mathbf v}}^{1/q} \|P_m\|_{L^{\infty}_{\mathbf x\mathbf v}}^{1-1/q}, 
\quad 1 < q < \infty. \nonumber 
\end{eqnarray}

Let us set $\phi_q (t) =  \|p_m\|_{L^{q}_{\mathbf x\mathbf v}}$. The above
inequalities yield Gronwall type inequalities of the form:
\begin{eqnarray}
\phi_1 (t)  \leq  \|p_0\|_1 + \alpha_1 \|\rho\|_{\infty}
\int_0^t \phi_1 (\tau) d\tau, \nonumber 
\\
e^{- N kt}\phi_{\infty} (t) \leq  \|p_0\|_{\infty} + \alpha_1 \|\rho\|_{\infty} 
\int_0^t e^{-N k\tau} \phi_{\infty}(\tau) d\tau. 
\nonumber 
\end{eqnarray}

Applying Gronwall's inequality to the above relations, we get:
\begin{eqnarray}
\|p_m (t)\|_{L^{1}_{\mathbf x\mathbf v}}  \leq  
\|p_0\|_{L^{1}_{\mathbf x\mathbf v}}  e^{\alpha_1 
\|\rho\|_{L^{\infty}_{\mathbf v}} t}, \label{cotap1}\\
\|p_m (t)\|_{L^{\infty}_{\mathbf x\mathbf v}}   \leq  
\|p_0\|_{L^{\infty}_{\mathbf x\mathbf v}} e^{N kt} 
e^{\alpha_1 \|\rho\|_{L^{\infty}_{\mathbf v}}t},  \label{cotapinfty} \\
\|p_m (t)\|_{L^{q}_{\mathbf x\mathbf v}} \leq  
\|p_0\|_{L^{1}_{\mathbf x\mathbf v}}^{1/q}  
\|p_0\|_{L^{\infty}_{\mathbf x\mathbf v}}^{1-1/q} 
e^{\alpha_1 \|\rho\|_{L^{\infty}_{\mathbf v}} t} e^{N kt (1-1/q)}, 
\quad 1 < q < \infty, \label{cotapq}
\end{eqnarray}
for all $m$ and $t \in [0,T]$.
Notice that these bounds on $p_m$ do not depend on 
${\mathbf F}(c_{m-1})$. We use that these fields are bounded functions to construct $p_m$ using fundamental solutions, but we do not need uniform bounds on either ${\mathbf F}(c_{m-1})$ or $a_{m-1}$ to obtain these $L^q$ bounds.

Additionally, Lemma 4.7 provides a uniform estimate on 
$\| \nabla_{\mathbf v} p_m \|_{L^2_{t\mathbf x \mathbf v}}$, that
does not require uniform bounds on ${\mathbf F}(c_{m-1})$ or $a_{m-1}$
either. For all $m$:
\begin{eqnarray}
\| \nabla_{\mathbf v} p_m \|_{L^2((0,T)\times I\!\! R^{N} \times I\!\! R^{N})}^2
\leq \|p_0\|_{L^2( I\!\! R^{N} \times I\!\! R^{N})}^2
e^{(2k+2\alpha_1 \|\rho\|_{L^\infty_{\mathbf v}} +1)T}. 
\label{cotadpv}
\end{eqnarray}
\vskip 3mm

\leftline{\it Step 4: Uniform bounds on velocity integrals of $p_m$}

We have obtained uniform estimates on the density norms
$\|p_m\|_{L^\infty(0,T;L^{q}_{\mathbf x\mathbf v})}$ for $1 \leq q \leq \infty$
and the angiogenic factor norm
$\|c_m\|_{L^\infty(0,T;L^{\infty}_{\mathbf x})}$. 

Lemma 4.3 provides a uniform bound of the $L^q_{\mathbf x}$
norms of ${ j}_{m-1}$, $1 \leq q \leq \infty$, by interpolating
the estimates
\begin{eqnarray}
\|{ j}_{m-1} \|_{L^1_{\mathbf x}} \leq   \||\mathbf v|g\|_{L^\infty_{\mathbf v}}
\|p_{m-1}\|_{L^1_{\mathbf x \mathbf v}}, \quad
\| { j}_{m-1} \|_{L^\infty_{\mathbf x}} \leq   \||\mathbf v|g\|_{L^1_{\mathbf v}}
\|p_{m-1} \|_{L^\infty_{\mathbf x \mathbf v}}. \label{jinftym} 
\end{eqnarray}
Indeed, $\| { j}_{m-1} \|_{L^q_{\mathbf x}} \leq  
\| { j}_{m-1} \|_{L^1_{\mathbf x}} ^{1/q} 
\| { j}_{m-1} \|_{L^\infty_{\mathbf x}} ^{1-1/q}. $
To establish uniform bounds on $a_{m-1}$ we adapt Proposition 4.4 to equation (\ref{eq:pm}), 
setting $a= \gamma \int_0^t \tilde p_{m-1} \, ds -\alpha(c_{m-1}) \rho$ and
${\mathbf F}={\mathbf F}(c_{m-1})$, with ${ j}={ j}_{m-1}$ depending  on $p_{m-1}$, not $p_m$. We replace (iv) with
\begin{itemize}
\item[(iv)']  ${\bf F}(c_{m-1})$ is given by (\ref{eq:alphaFm}), $c_{m-1}$, $\nabla_{\mathbf x} c_{m-1}$ are given by (\ref{cmint})-(\ref{dcmint}).
$c_{m-1}\geq 0$ is coupled to equation  (\ref{eq:pm}) for $p_{m-1}$ by (\ref{eq:jm}), 
and $p_m\geq 0$ is coupled to equation (\ref{eq:cm}) for $c_{m-1}$ through (\ref{eq:alphaFm}).
\end{itemize}
We know by Step 1 that $a \in L^{\infty}$. Its negative part is 
$a^-(c_{m-1})=\alpha(c_{m-1})\rho$,  that satisfies 
$\|a^-\|_{\infty}\leq \alpha_1 \|\rho \|_{L^\infty_{\mathbf v}}.$   
Thanks to Step 2, Lemma 4.1 and the above uniform estimate on $\| j_{m-1}\|_{L^q_{\mathbf x}}$, $q>N$, we obtain a uniform bound on 
$\|\mathbf F (c_{m-1})\|_{\infty}$. The linear part of Proposition 4.4,
combined with these uniform estimates on $\|a^-(c_{m-1})\|_{\infty}$
and $\|\mathbf F (c_{m-1})\|_{\infty}$, yields a uniform bound on
$\|(1\!+\! |{\mathbf  v}|^2)^{\beta/2}  p_{m-1} \|_{L^\infty_{\mathbf x \mathbf v}}$.
Then, inequality (\ref{boundinf}) in Lemma 4.2  provides a uniform estimate on
the norms
$\|p_m\|_{L^\infty(0,T;L^{\infty}_{\mathbf x} L^1_{\mathbf v})} $.

Therefore,  the coefficients $a_{m-1}$, $\alpha(c_{m-1})$, ${ j}_{m-1}$ 
and $F(c_{m-1})$ appearing in the equations are uniformly bounded in 
$L^{\infty}(0,T; L^{\infty}_{\mathbf x})$. From Step 2, we also infer that  
the norms $\|\nabla_{\mathbf x} c_m\|_{L^\infty(0,T; L^\infty_{\mathbf x})}$,
$\|\nabla_{\mathbf x} \tilde c_m\|_{L^\infty(0,T; L^2_{\mathbf x})}$,
and $\|\tilde c_m\|_{L^\infty(0,T; L^2_{\mathbf x})}$,  are uniformly bounded too. 
\\
 
\leftline{\it Step 5: Compactness of the iterates}

The passage to the limit with a minimal set of uniform bounds is made possible 
by the following compactness result, taken from \cite{bouchut2}:\\

{\bf Lemma 5.2.}  {\it Let $\sigma > 0$, $k \geq 0$, $T >0$, $1 \leq q < \infty$,
$p_0 \in L^q(I \!\! R^{N})$, $f \in L^1(0,T;L^q(I \!\! R^{N}\times I \!\! R^{N}))$ and
consider the solution $p \in C([0,T], L^q( I\!\! R^{N} \times I\!\! R^{N}))$ of:
\begin{eqnarray}
{\partial p \over \partial t} + {\mathbf v} \nabla_{\mathbf x} p
- k \nabla_{\mathbf v} \cdot ({\mathbf v}  p) - \sigma \Delta_{\mathbf v}  p
&=& f \quad {\rm in} \,(0,T) \times I\!\! R^{N} \times I\!\! R^{N}, \label{eq:k0} \\
p(0) &=& p_0  \quad {\rm in} \; I\!\! R^{N} \times I\!\! R^{N}. \nonumber 
\end{eqnarray}
Assume that $p_0$ belongs to a bounded subset of $L^q( I\!\! R^{N} \times I\!\! R^{N})$ and $f$ belongs to a bounded subset of $L^r(0,T;L^q(I \!\! R^{N}\times I \!\! R^{N}))$ with $1 < r \leq \infty$. Then, for any $\eta >0$ and any bounded open subset
$\omega$ of $I \!\! R^{N}\times I \!\! R^{N}$, $p$ is compact in $C([\eta,T],L^q(\omega))$.
}\\

Previous compactness results by R.J. DiPerna and P.L. Lions \cite{diperna}  guarantee compactness in $L^1(0,T;L^q(\omega))$ for more general operators.

In our case, $p_0 \in L^{1}\cap L^{\infty}( I\!\! R^{N} \times I\!\! R^{N})$ is fixed.  The solutions $p_m$ of the iterative scheme satisfy (\ref{eq:k0}) with:
\begin{eqnarray}
f_m =\rho({\mathbf v}) \alpha(c_{m-1}(t,{\mathbf x})) p_{m}(t,{\mathbf x},{\mathbf v}) 
- \gamma a(p_{m-1}(t,\mathbf{x}))  p_m(t,\mathbf{x},\mathbf{v})\nonumber \\
- \mathbf{F}(c_{m-1}(t,\mathbf{x})) \nabla_\mathbf{v}p_m(t,\mathbf{x},\mathbf{v}), 
\nonumber \\
\alpha(c_{m-1})=\alpha_1\frac{\frac{c_{m-1}}{c_R}}{1+\frac{c_{m-1}}{c_R}}, \quad
{\bf F}(c_{m-1})= \frac{d_1}{(1+\gamma_1c_{m-1})^{q_1}}\nabla_{\mathbf x} c_{m-1}. \nonumber
\end{eqnarray}

The estimates established in Steps 2-4 and Lemma 4.5 ensure that:
\begin{itemize}
\item the coefficient functions $\alpha(c_{m-1})$, ${\bf F}(c_{m-1})$, $a(p_{m-1})$ 
are uniformly bounded in $L^{\infty}(0,T;L^{\infty}_{\mathbf x})$.
\item $p_m, \nabla_\mathbf{v}p_m$ are bounded in 
$L^2(0,T;L^2_{\mathbf x \mathbf v})$.
\end{itemize}
Lemma 5.2 guarantees that $p_m$ is compact in $C([\eta,\tau],L^2(\omega))$ 
for any $\omega$ and $\eta$. 

The integral inequalities (\ref{tildecmq})-(\ref{dcm2}) provide a uniform estimate on $\tilde{c}_m$ in $L^{2}(0,T;H^1(\omega))$, for any bounded $\omega \!\subset\! I\!\!R^{N}$. Using equation (\ref{eq:tildecm}), we bound ${\partial \tilde{c}_m \over \partial t}$   
in $L^2(0,T;H^{-1}(\omega))$. 
Standard compactness results in \cite{lions,simon} yield compactness in $L^2(0,T;L^2(\omega))$ for any bounded $\omega$.  

In view of the uniform bounds established in Sections 2, 3 and 4, we may extract
sequences $p_{m_i}$ and $\tilde c_{m_i}$ that converge weakly in all the pertinent spaces to limits $p$ and $c$. Compactness on bounded sets $\omega$ allows us to extract subsequences, again denoted $p_{m'}$ and $\tilde c_{m'}$, that converge strongly in $L^2_{loc}$ and pointwise almost everywhere to $p$ and $\tilde c$ by a diagonal extraction procedure. To exemplify the process, let us consider a family of balls $B_M$ of radius $M$ in $I\!\! R^{N}$. For $M=1$ we have a sequence $\tilde c_{m_i^1}$ that converges strongly in $L^2$ and pointwise a.e. to $\tilde c$ in $B_1$ as $i$ tends to infinity. For $M=2$, we extract from this sequence another sequence $\tilde c_{m_i^2}$ that converges strongly in $L^2$ and pointwise a.e. to $\tilde c$ in $B_2$ as $i$ tends to infinity.
By induction, for $M=j$ we may extract from $\tilde c_{m_i^{j-1}}$ a sequence 
$\tilde c_{m_i^{j}}$ that converges strongly in $L^2$ and pointwise a.e. to $\tilde c$ in $B_j$ as $i$ tends to infinity.
The diagonal subsequence $\tilde c_{m'}=\tilde c_{m_i^{i}}$ converges strongly in $L^2_{loc}$ and pointwise a.e. to $\tilde c$ as $i$ tends to infinity in any subset $B_j$, therefore in the whole space. Then, the sequence $c_{m'}$ tends to $c=\tilde{c}+{\cal C}$ pointwise. A similar argument constructs $p_{m'}.$

We have global pointwise and weak convergences, plus strong local convergences. Combined with uniform bounds in terms of integrable functions in the whole space, this will allow us to pass to the limit in the nonlocal terms and integrals appearing in Step 6. Let us construct the required controlling functions, exploiting representations in terms of fundamental solutions and the fact that the norms $\|\mathbf F(c_m)\|_{\infty}$ are uniformly bounded. 

Let us start with $p_m$.
We recall that $0 \leq p_{m}\leq P_m$, where $P_m$ is a solution of (\ref{eq:Pm}). Then, the right hand side satisfies
$\alpha(c_{m-1}) \rho\, p_{m} \leq  \alpha_1 \|\rho\|_{\infty} P_m.$
Therefore, $P_m \leq Q_m$, where $Q_m$ is the solution of 
\begin{eqnarray} 
\frac{\partial}{\partial t} Q_m \!+\! 
\mathbf{v} \!\cdot\! \nabla_\mathbf{x}   Q_m 
\!+\!  \nabla_\mathbf{v} \!\cdot\! [({\mathbf F}(c_{m-1}) 
\!-\! k\mathbf{v} ) Q_m ]   
\!-\! \sigma \Delta_\mathbf{v} Q_m = 
\alpha_1 \|\rho\|_{\infty} Q_m,   \label{eq:Qm} 
\end{eqnarray}
with initial datum $p_0$.  By Lemma 2.6, $Q_m \leq {\cal P}$ where  ${\cal P}=e^{\alpha_1 \|\rho\|_{\infty} t} P$ and $P$ is a solution of (\ref{upper2}) with  
$M(T, \|{\bf F}(c_{m-1})\|_{\infty})$ replaced by a constant $M(T)$ in the initial
datum.
Notice that $\|{\bf F}(c_{m-1})\|_{\infty}$ is uniformly bounded as a consequence of Step 4, therefore $M(T, \|{\bf F}(c_{m-1})\|_{\infty}) \leq M(T).$

By comparison principles for heat equations, the functions 
$|\tilde{c}_m|$ are uniformly bounded by the solution of
\begin{eqnarray}
{\partial \over \partial t} {\cal V} - d \Delta_{\mathbf x} {\cal V}
= \eta \|c_0\|_{\infty} { j}({\cal P}),
\quad {\cal V}(0)=0. \label{eq:U}
\end{eqnarray}
%
%
Setting $f=\eta \|c_0\|_{\infty} { j}({\cal P})$, we have ${\cal V}(t,\mathbf x) = \int_0^t ds \int_{I\!\!R^{N}} d \mathbf y \, K(t-s,\mathbf x-\mathbf y) f(s,\mathbf y),$ where $K$ is the heat kernel. 

Additionally, compactness of the sequences $p_m$ in  $L^2(0,T,L^2(I \!\! R^{N} \times I \!\! R^{N}))$ and $\tilde c_m$ in $L^2(0,T;L^2(I \!\! R^{N}))$ follows from compactness in bounded sets since  we  control uniformly  the decay of the sequences at infinity. 
This results from a standard strategy to establish compactness in spaces 
$L^q(I \!\! R^n)$, $1\leq q<\infty$, see reference \cite{brezis}. We must prove that
for any $\varepsilon >0$ we can find a radius $R>0$ such that the norms of the sequences outside that ball are smaller than $\varepsilon$.


Let us estimate the $L^2$ norm of $\tilde c_m$ outside a ball. Writing 
\[K(t\!-\!s,\mathbf x\!-\!\mathbf y) f(s,\mathbf y)=
K(t\!-\!s,\mathbf x\!-\!\mathbf y)^{1\over 2} K(t\!-\!s,\mathbf x\!-\!\mathbf y)^{1\over 2} f(s,\mathbf y),\] 
we obtain \cite{brezis}:
\[
|K(t-s)*_{\mathbf x} f(s)|^2 \leq \|K(t-s)\|_1 (K(t-s)*_{\mathbf x} 
|f(s)|^2).
\]
We split the norm of $(K(t-s)*_{\mathbf x} 
|f(s)|^2)$ outside a ball or radius $2R>0$ as:
\[
\int_{|\mathbf x| > 2R} \hskip -5mm d{\mathbf x} 
\int_{|\mathbf y| > R} \hskip -6mm d{\mathbf y} 
\, K(t-s,\mathbf x -\mathbf y)  |f(s, \! \mathbf y)|^2 \!+\!
\int_{|\mathbf x| > 2R} \hskip -5mm d{\mathbf x} 
\int_{|\mathbf y| < R} \hskip -6mm d{\mathbf y} 
\, K(t-s,\mathbf x -\mathbf y)  |f(s, \! \mathbf y)|^2.
\]
Notice that ${ j}(P) \in L^{\infty}(0,T;L^2_{\mathbf x})$ 
thanks to Lemma 4.3.
For the second integral $I_2$, $|\mathbf x -\mathbf y| >R$. Both integrals are bounded by:
\[
I_2 \leq {1 \over \pi^{N\over 2}} \| e^{-|\mathbf z|^2} \|_{L^1(|\mathbf z|>{R\over 2\sqrt{T}}  ) } \|f(s)\|_{L^2_{\mathbf x}}^2, \quad
I_1 \leq  {1 \over \pi^{N\over 2}} \| e^{-|\mathbf z|^2} \|_{1}
 \|f(s)\|_{L^2(|\mathbf x|>R)}^2.
\]
%
We can make $\|{\cal V}\|_{L^2(0,T,L^2(|\mathbf x| >R))}$ as small as wished  choosing $R$ large provided the same holds for
$ \|f\|_{L^\infty(0,T,L^2(|\mathbf x| >R))}$. With our definition of ${ j}$,
this leads to the study
of the norm $ \| \tilde P\|_{L^\infty(0,T,L^2(|\mathbf x| >R))}$. 
As stated in Lemma 2.6, $P$ admits the integral expression:
\begin{eqnarray}
P(t,{\mathbf x},{\mathbf v})=
M(T)  \!\! \int_{I\!\!R^{N} \times I\!\!R^{N}} 
\hskip -7mm 
G(t, {{\mathbf x}\over 2}, {{\mathbf v}\over 2};  0, {{\boldsymbol \xi}\over 2}, {{\boldsymbol \nu}\over 2}) p_0({\boldsymbol \xi}, {\boldsymbol \nu}) 
d{\boldsymbol \xi} d{\boldsymbol \nu},
\label{upperP}
\end{eqnarray}
where $G$ is the fundamental solution for the free field linear operator
given in (\ref{sf0}). Arguing in a similar way as we just did with the convolution
with the heat kernel, and taking into account that now $p_0$ is fixed and
does not depend on time, we find that
$ \|P\|_{L^2(0,T,L^2(|\mathbf x| >R, |\mathbf v|>R))}$,
$ \|\tilde P\|_{L^\infty(0,T,L^2(|\mathbf x| >R))}$ can be made as
small as wished choosing $R$ large enough. 
Notice that $ \|P\|_{L^2(0,\eta,L^2_{\mathbf x \mathbf y})}$
can be made arbitrarily small decreasing $\eta$.

Therefore, the sequence $p_m$ is compact in  $L^2(0,T,L^2(I \!\! R^{N} \times I \!\! R^{N}))$ and $\tilde c_m$ is compact
in $L^2(0,T;L^2(I \!\! R^{N})).$  
\\

\leftline{\it Step 6: Convergence to a solution}


The function ${\cal P}$ constructed in Lemma 2.6 and used to control uniformly
the sequence $p_m$ in Step 5, provides uniform integrable bounds to apply Lebesgue's dominated convergence theorem. For $t$ and ${\mathbf x}$ fixed, 
and for any ${\mathbf v}$, the subsequence $p_m'$ converging pointwise satisfies 
$0 \leq p_{m'} \leq {\cal P}$.
The functions ${\cal P}$ and $|{\mathbf v}|   g(|{\mathbf v}|) {\cal P}$ are integrable with respect to all the variables. The sequences $p_{m'}$ and ${|{\mathbf v}|}   g(|{\mathbf v}|)  p_{m'}$ converge pointwise to $p$ and ${|{\mathbf v}|}  g(|{\mathbf v}|)  p$, respectively. Therefore, $a(p_{m'-1})$ and ${ j}(p_{m'-1})$ tend to $a(p)$ and ${ j}(p)$. 

Let us pass to the limit in the weak formulation of the equations. For all $\phi \in C_c^2([0,T) \times I\!\!R^{N} \times I\!\!R^{N})$   
\begin{eqnarray}
\int_0^T \hskip -3mm \int_{I\!\!R^{N} \times I\!\!R^{N}} \hskip -10mm p_{m'}(t, \!{\mathbf x}, \!{\mathbf v}) \Big[
\frac{\partial}{\partial t} \phi(t,\!\mathbf{x},\!\mathbf{v}) \!+\! 
\mathbf{\!v} \!\!\cdot\!\! \nabla_\mathbf{\!x}   \phi(t,\!\mathbf{x},\!\mathbf{v}) 
\!+\!   ({\mathbf F}(c_{m'-1}(t,\!\mathbf{x})) \!-\! k\mathbf{v} ) \!\!\cdot\!\!
\nabla_\mathbf{\!v} \phi(t,\!\mathbf{x},\!\mathbf{v})  \nonumber \\  \!+ \sigma
\Delta_\mathbf{v} \phi(t,\!\mathbf{x},\!\mathbf{v}) \!-\! [a(p_{m'-1}(t,\!\mathbf x))\!-\! \alpha(c_{m'-1}(t,\!\mathbf x) \rho({\mathbf v})] \phi(t,\!\mathbf{x},\!\mathbf{v})  \Big] 
d{\mathbf x}  d{\mathbf v} dt  \nonumber \\
+ \int_{I\!\!R^{N} \times I\!\!R^{N}} \hskip -6mm \phi(0,\! \mathbf{x},\!\mathbf{v}) 
p_0(\mathbf{x},\!\mathbf{v}) d{\mathbf x}  d{\mathbf v} =0. 
\label{weakm}
\end{eqnarray}

Since ${\cal P}\in L^q_{t \mathbf x \mathbf v}$ for any $q \in [1,\infty]$, 
$p_{m'}$ converges to $p$ in $L^q_{t \mathbf x \mathbf v}$
for all finite $q$. The sequence $p_{m'} a(p_{m'-1})$ converges pointwise to 
$p a(p)$ and $p_{m'} a(p_{m'-1})$ is bounded by ${\cal P} a({\cal P})\in L^q_{t \mathbf x \mathbf v}$ for any $q \in [1,\infty]$. Thus, it also converges in $ L^q_{t \mathbf x \mathbf v}$ for finite $q$. 

Pointwise convergence of $c_{m'-1}$ together with continuity of the function $\alpha(s)$ imply pointwise convergence of $\alpha(c_{m'-1})$ to $\alpha(c)$. The sequence $p_{m'} \alpha(c_{m'-1})$ converges pointwise to $p \alpha(c)$ and $p_{m'} \alpha(c_{m'-1})$ is bounded by ${\cal P} \alpha_1\in L^q_{t \mathbf x \mathbf v}$ for any $q \in [1,\infty]$. Therefore, we have convergence in $ L^q_{t \mathbf x \mathbf v}$ for finite $q$.

These convergences allows to pass to the limit in all the terms present in identity
(\ref{weakm}), except in 
\begin{eqnarray}
\int_0^T \hskip -3mm \int_{I\!\!R^{N} \times I\!\!R^{N}} \hskip -8mm
p_{m'}(t, \!{\mathbf x}, \!{\mathbf v}) \frac{d_1}{(1+\gamma_1c_{m'-1})^{q_1}}
\nabla_{\mathbf x} c_{m'-1}(t,\!\mathbf{x}) \cdot \nabla_{\mathbf v} 
\phi(t,\!\mathbf{x},\!\mathbf{v}) d{\mathbf x}  d{\mathbf v}.
\label{extraterm}
\end{eqnarray}
The sequence $ p_{m'} \frac{d_1}{(1+\gamma_1c_{m'-1})^{q_1}} \nabla_{\mathbf v} \phi$ tends  pointwise to $ p \frac{d_1}{(1+\gamma_1c)^{q_1}} \nabla_{\mathbf v} \phi$ and is bounded by ${\cal P} d_1 \nabla_{\mathbf v} \phi \in L^q_{t \mathbf x \mathbf v}$ for any $q \in [1,\infty]$. Thus, we have strong convergence in $L^q_{t \mathbf x \mathbf v}$ for all finite $q$. 
The sequence $\nabla_{\mathbf x} c_{m'-1}$ is bounded in
$L^{2}_{t \mathbf x \mathbf v}$. Therefore, it tends
weakly to $\nabla_{\mathbf x} c$ in $L^{2}_{t \mathbf x \mathbf v}$ .  These convergences allow us to conclude the process. The limit $(p,c)$ satisfies:
\begin{eqnarray}
\int_0^T \hskip -3mm \int_{I\!\!R^{N} \times I\!\!R^{N}} \hskip -8mm p(t, \!{\mathbf x}, \!{\mathbf v}) \Big[
\frac{\partial}{\partial t} \phi(t,\!\mathbf{x},\!\mathbf{v}) \!+\! 
\mathbf{\!v} \!\!\cdot\!\! \nabla_\mathbf{\!x}   \phi(t,\!\mathbf{x},\!\mathbf{v}) 
\!+\!   ({\mathbf F}(c(t,\!\mathbf{x})) \!-\! k\mathbf{v} ) \!\!\cdot\!\! 
\nabla_\mathbf{\!v} \phi(t,\!\mathbf{x},\!\mathbf{v})  \nonumber \\  \!+ \sigma
\Delta_\mathbf{v} \phi(t,\!\mathbf{x},\!\mathbf{v}) \!-\! [a(p(t,\!\mathbf x))\!-\! \alpha(c(t,\!\mathbf x) \rho({\mathbf v})] \phi(t,\!\mathbf{x},\!\mathbf{v})  \Big] 
d{\mathbf x}  d{\mathbf v} dt  \nonumber \\
+ \int_{I\!\!R^{N} \times I\!\!R^{N}} \hskip -6mm \phi(0,\! \mathbf{x},\!\mathbf{v}) 
p_0(\mathbf{x},\!\mathbf{v}) d{\mathbf x}  d{\mathbf v} =0. 
\label{weakp}
\end{eqnarray}
This solution inherits all the bounds established for the converging sequences.

It remains to pass to the limit in the parabolic equations (\ref{eq:cm}).
The source term $c_{m'} { j}(p_{m'})$ tends pointwise to $c { j}(p)$
and is bounded by ${\cal C} { j}({\cal P}) 
\in L^q_{t \mathbf x }$ for any 
$q \in [1,\infty]$. Strong convergence of the source to its limit implies
that $c$ is a solution of (\ref{eq:c}).\\

\leftline{\it Step 7: Regularity}

Once existence of a solution of the nonlinear problem has been proved, it has
at least the regularity of solutions of linear Fokker-Planck equations 
with force field ${\bf F} \in L^{\infty}_t L^{\infty}_{\mathbf x}$ and either source 
$[\rho \alpha(c) - a(p) ] p \in L^{\infty}_tL^q_{\mathbf x \mathbf v}$ for all
$q$, or lower order term $[\rho \alpha(c) - a(p) ] p$ with coefficient
$ [\rho \alpha(c) - a(p) ]  \in L^{\infty}_tL^\infty_{\mathbf x \mathbf v}$.
Using the integral equation for derivatives with respect to $v$ and the 
estimate (F5) in Lemma 2.2 we obtain additional regularity adding 
hypotheses to the derivatives of the initial datum.\\
 
\leftline{\it Step 8: Uniqueness}

Set $\overline{p}=p_1-p_2$ and $\overline{c}=c_1-c_2$.  
These differences satisfy the equations:
\begin{eqnarray} 
\frac{\partial}{\partial t} \overline{p}\!+\! 
\mathbf{v} \!\cdot\! \nabla_\mathbf{x}   \overline{p}
\!+\!  \nabla_\mathbf{v} \!\cdot\! [({\mathbf F}(c_{1})
\!-\! k\mathbf{v} ) \overline{p} ]  \!-\! \sigma \Delta_\mathbf{v} \overline{p}
\!= \! [-\gamma a(p_{1}) \!+\! \alpha(c_{1}) \rho] \overline{p}  
\label{eq:pdif}  \\
 \!-\!   [{\mathbf F}(c_{1})-{\mathbf F}(c_{2})] \nabla_\mathbf{v}  
p_2
 \!+\! [-\gamma a(\overline{p}) \!+\! (\alpha(c_{1})-\alpha(c_2)) \rho] p_2,
\nonumber\\
\frac{\partial}{\partial t} \overline{c}\!-\! d
\Delta_{\mathbf x} \overline{c} + \eta { j} (p_1) \overline{c} \!= 
- \eta  [{ j}(p_1)- { j}(p_2)] c_2,
\label{eq:cdif} 
\end{eqnarray}
with $\overline{p}(0)=0$ and $\overline{c}(0)=0$. 

The mean value theorem applied to the definitions of ${\mathbf F}(c)$ 
and $\alpha(c)$ yields:
\begin{eqnarray}
|{\mathbf F}(c_1)-{\mathbf F}(c_2)|  
\leq {q_1 d_1 \gamma_1 |\nabla_{\mathbf x} c_1|
\over(1+\gamma_1 \xi)^{q_1+1}} |c_1-c_2| +
{d_1 \over (1+\gamma_1 c_2)^{q_1}} 
|\nabla_{\mathbf x} c_1- \nabla_{\mathbf x} c_2|,
\label{meanF}\\
|\alpha(c_1)\!-\!\alpha(c_2)| =  \alpha_1
|{c_1 \over c_R + c_1} \!-\! {c_2 \over c_R + c_2}|
=  {\alpha_1 c_R\over (c_R+ \chi)^2} |c_1\!-\!c_2|
\leq {\alpha_1 \over c_R} |c_1\!-\!c_2|,
\label{meanalpha}
\end{eqnarray}
where $\xi, \chi \in [c_1,c_2]$. Since $c_1$ and $c_2$ are
nonnegative, $\xi, \chi \geq 0$.
Particularizing (\ref{int1}) for (\ref{eq:pdif}) and making use
of (\ref{meanF})-(\ref{meanalpha}) we find:
\begin{eqnarray}
\|\overline {p}(t)\|_{L^1_{\mathbf x \mathbf v}} \leq 
[\gamma \|a(p_{1})\|_{\infty} \!+\! \alpha_1 \|\rho\|_{\infty}] 
\int_0^t ds \|\overline{p}(s) \|_{L^1_{\mathbf x \mathbf v}} 
+ \nonumber  \\
d_1 q_1 \gamma_1 \|\nabla_{\mathbf x} c_1\|_{\infty} 
\| \nabla_\mathbf{v}   p_2 \|_{L^{\infty}_t L^{\infty}_{\mathbf x}L^1_{\mathbf v}}
\int_0^t \!\! ds \|\overline{c}(s)\|_{L^1_{\mathbf x}} 
+ \nonumber \\
d_1 \| \nabla_\mathbf{v}   p_2 \|_{L^{\infty}_t L^{\infty}_{\mathbf x}L^1_{\mathbf v}}
\int_0^t \!\! ds \|\nabla_{\mathbf x}\overline{c}(s)\|_{L^1_{\mathbf x}}
\!+\!   \nonumber \\
\gamma \|p_2\|_{L^\infty_t L^{\infty}_{\mathbf x} L^1_{\mathbf v}} 
\int_0^t \!\! ds \!\! \int_0^s \!\! d\tau 
\|\overline{p}(\tau)\|_{L^{1}_{\mathbf x \mathbf v}}
\!+\! {\alpha_1 \|\rho\|_{\infty}\over c_R} \|p_2\|_{L^\infty_t L^{\infty}_{\mathbf x} L^1_{\mathbf v}} \int_0^t \!\! ds  \|\overline{c}(s)\|_{L^1_{\mathbf x}}. 
\label{intp1dif}
\end{eqnarray}

Expressing the solution of (\ref{eq:cdif}) in integral form in terms
of the fundamental solution of $c_t- d \Delta_{\mathbf x} c
+ \eta |{ j}(p_1)| c$ we obtain  \cite{gema}:
\begin{eqnarray}
\|\overline{c}(t)\|_{L^1_{\mathbf x}}  \leq   
C(\|{ j}(p_1)\|_{\infty}) \eta t \|c_2\|_{\infty}  
\| { j}(\overline{p}) \|_{L^{\infty}(0,t;L^1_{\mathbf x})}.
\label{intc1dif}
\end{eqnarray}
We have used that the fundamental solution is bounded by scaled
heat kernels \cite{aronson,gema,kusuoka} involving constants that depend
on the $L^{\infty}$ norm of the coefficient $|{ j}(p_1)| $.
Particularizing (\ref{intdc}) for this equation, using 
$\|\nabla K(t)\|_1\leq Mt^{-1/2}$ and estimate (\ref{intc1dif}), 
we find:
\begin{eqnarray}
\|\nabla_{\mathbf x}\overline{c}(t)\|_{L^1_{\mathbf x}} 
\leq 2 \eta M t^{1/2}  [ \|c_2\|_{\infty} 
\| { j}(\overline{p}) \|_{L^{\infty}(0,t;L^1_{\mathbf x})}
+ \|{ j}(p_1)\|_{\infty} 
\| \overline{c} \|_{L^{\infty}(0,t;L^1_{\mathbf x})} ]
\label{intdc1dif}  \\  \leq 2 \eta M [
\|c_2\|_{\infty} t^{1/2} 
+ C(\|{ j}(p_1)\|_{\infty}) 
\|c_2\|_{\infty} \eta   t^{3/2} ] 
\| { j}(\overline{p}) \|_{L^{\infty}(0,t;L^1_{\mathbf x})}.
\nonumber
\end{eqnarray}

By Lemma 4.3. we have:
\begin{eqnarray}
\|{ j}(\overline{p})\|_{L^1_{\mathbf x}} \leq 
\||\mathbf v|g\|_{L^\infty_{\mathbf v}}
\|\overline{p}\|_{L^1_{\mathbf x \mathbf v}}.
\label{uniqj}
\end{eqnarray}
Let us set:
\begin{eqnarray}
A= A_1 +A_2 = [\gamma \|a(p_{1})\|_{\infty} 
\!+\! \alpha_1 \|\rho\|_{\infty}]  
+ T \gamma \|p_2\|_{L^\infty_t L^{\infty}_{\mathbf x}L^1_{\mathbf v}},
\nonumber \\
B= B_1 +B_2 =  d_1 q_1 \gamma_1 \|\nabla_{\mathbf x} c_1\|_{\infty} 
\| \nabla_\mathbf{v}   p_2 \|_{L^{\infty}_t L^{\infty}_{\mathbf x}L^1_{\mathbf v}}
+ {\alpha_1 \|\rho\|_{\infty}\over c_R} \|p_2\|_{L^\infty_t L^{\infty}_{\mathbf x} L^1_{\mathbf v}}, 
\nonumber \\
D= d_1  \| \nabla_\mathbf{v}   p_2 \|_{L^{\infty}_t L^{\infty}_{\mathbf x}L^1_{\mathbf v}},
\nonumber \\
E = E_1 + 2 \eta M T^{1/2} E_2 =  2 \eta M \|c_2\|_{\infty} T^{1/2} 
+  2 \eta^2 M \|c_2\|_{\infty} C(\|{ j}(p_1)\|_{\infty}) T^{3/2}.
\nonumber
\end{eqnarray}
We define now 
\begin{eqnarray*}
U(t)= {\rm max}_{s\in[0,t]}  \|\overline{p}(s) \|_{L^1_{\mathbf x \mathbf v}}.
\label{controldif}
\end{eqnarray*}
Combining (\ref{intp1dif})-(\ref{uniqj}) we find: 
\begin{eqnarray}
\|\overline {p}(t)\|_{L^1_{\mathbf x \mathbf v}} \leq 
A \int_0^t ds \|\overline{p}(s) \|_{L^1_{\mathbf x \mathbf v}} 
+ ( B E_2 + D E)
\int_0^t ds \| { j}(\overline{p}) \|_{L^{\infty}(0,s;L^1_{\mathbf x})}
\label{control1} \\
\leq  \left(A + ( B E_2 + D E) 
 \||\mathbf v|g\|_{L^\infty_{\mathbf v}} 
\right) \int_0^t U(s) ds.
\end{eqnarray}
We deduce that $U(t)$ satisfies a Gronwall inequality of the form
\[ U(t) \leq G(T) \int_0^t U(s) ds\] for $t\in [0,T]$. Therefore, $U=0$
and $p_1=p_2$ in $[0,T]$ for any $T>0$. 
Estimate (\ref{intc1dif}) implies then $c_1=c_2$.




%
%

\section{Fluxes without velocity cut-offs}

In this section, we replace the flux defined in formula (\ref{eq:intpintvp}) by
\begin{eqnarray}
\mathbf{j}(t,\mathbf{x})= \int_{I\!\!R^N}  \hskip -2mm \mathbf{v}  p(t,\mathbf{x},\mathbf{v})\, d \mathbf{v}, \quad
{ j(t,\mathbf{x})= \int_{I\!\!R^N} \hskip -2mm
|\mathbf{v}|  p(t,\mathbf{x},\mathbf{v})\, d \mathbf{v}}. \label{flux}
\end{eqnarray}
We reconsider the system (\ref{eq:p})-(\ref{eq:alphaF}), (\ref{flux}) and the iterative
scheme  (\ref{eq:pm})-(\ref{eq:cm0}) with (\ref{eq:jm}) replaced by 
\begin{eqnarray}
\mathbf{j}_{m-1}(t,\mathbf{x})= \int_{I\!\!R^N} \hskip -2mm \mathbf{v}  p_{m-1}(t,\mathbf{x},\mathbf{v})\, d \mathbf{v},
\quad
{ j_{m-1}(t,\mathbf{x})= \int_{I\!\!R^N} \hskip -2mm |\mathbf{v}|  p_{m-1}(t,\mathbf{x},\mathbf{v})\, d \mathbf{v}}.
\label{fluxm}
\end{eqnarray}
In this case, local existence in time can be proven. This restriction comes
from the fact that we need uniform estimates on $\| \mathbf F(c_{m}) \|_{\infty}$
to be able to control uniformly the behavior of the sequences as $|\mathbf xÊ| \rightarrow \infty$ through a uniform estimate of the fundamental solutions
and pass to the limit in the weak formulation of the equations. Uniform
estimates on $\| \mathbf F(c_{m}) \|_{\infty}$ require uniform estimates
on $\| { j}(p_{m}) \|_{q}$ for $q$ large enough. 
In absence of Lemma 4.3,
we may obtain them via Lemma 4.2 from uniform estimates on velocity
moments whose order is high enough. At this time, the restriction appears,  
affecting also the counterpart of Proposition 4.4 for this choice of flux 
${ j}$. Let us first estimate the velocity moments. \\

{\bf Proposition 6.1.} {\it Let $p \geq 0$ be a solution of (\ref{eq:linpaj})-(\ref{eq:linp0aj}) under the hypotheses:
\begin{itemize}
\item [(i)] $a \in L^{\infty}((0,T)\times I\!\!R^{N} \times I\!\!R^{N})$, 
${\bf F} \in L^{\infty}((0,T)\times I\!\!R^{N})$,
\item[(ii)] $ (1+|{\mathbf  v}|^2)^{\beta/2} p_0 \in L^1(I\!\!R^{N} \times I\!\!R^{N})$, 
$ p_0 \in L^\infty(I\!\!R^{N} \times I\!\!R^{N})$, $p_0\geq 0$,
\item[(iii)] ${\bf F}(c(t,{\mathbf x}))$ is given by 
(\ref{eq:alphaF}), $c$, $\nabla_{\mathbf x} c$ are given by (\ref{intc})-(\ref{intdc}) with $c_0 \in L^{\infty}_{\mathbf x}$, 
$\nabla_{\mathbf x} c_0 \in L^{N+\beta}_{\mathbf x}$, $c\geq 0$,
$c$ is coupled to (\ref{eq:linpaj}) by (\ref{eq:intpintvp}), with
$a= \gamma \int_0^t \tilde p \, ds -\alpha(c) \rho$,
\end{itemize}
for $\beta > N^2 - N \geq 2$. Then,  there exists $\tau_{\beta}$
such that for any $\tau < {\rm min}(\tau_\beta,T)$ the norms
$\| |{\mathbf v}|^2 p \|_{L^{\infty}(0,\tau;L^1_{\mathbf x \mathbf v})}$, 
$\| |{\mathbf v}| p \|_{L^{\infty}(0,\tau;L^1_{\mathbf x \mathbf v})}$,
$\|{ j}\|_{L^{\infty}(0,\tau;L^q_{\mathbf x})}$, 
$q \in [1, {N+\beta \over N+1}]$,
and $\|\tilde p \|_{L^{\infty}(0,\tau;L^q_{\mathbf x})}$, 
$q \in [1,  {N+\beta \over N}]$,
are bounded by constants depending  on the parameters 
 $\sigma$, $k$, $d$, $\eta$,  $d_1$,  $\alpha_1$, $T$, $\beta$, $N$, 
and the norms $\|c_0\|_{L^{\infty}_{\mathbf x}}$,  
$\|\nabla_{\mathbf x}c_0\|_{L^{N+\beta}_{\mathbf x}}$
$\| |{\mathbf  v}|^2 p_0 \|_{L^1_{\mathbf x \mathbf v}}$,
$\|\rho \|_{L^\infty_{\mathbf v}}$,
$\| p \|_{L^{\infty}(0,T;L^1_{\mathbf x \mathbf v} \cap L^\infty_{\mathbf x \mathbf v})}$.\\
}

{\bf Proof.} Using fundamental solutions  we get $p \in C([0,T];L^{1}_{\mathbf x \mathbf v}\cap L^{\infty}_{\mathbf x \mathbf v} ),$ see references \cite{bouchut3,bouchut2,gema}. Moreover, $(1+|\mathbf v |^2)^{\beta/2} p \in C([0,T];L^{1}_{\mathbf x \mathbf v})$. 

The conservation law for moments of order $\beta$ reads \cite{carrillo}:
\begin{eqnarray}
{d \over dt} \int_{I\!\!R^{N} \times I\!\!R^{N}} \hskip -8mm |{\mathbf  v}|^{\beta} 
p \, d{\mathbf  x} d{\mathbf  v}
= \beta(\beta-2+N) \sigma \int_{I\!\!R^{N} \times I\!\!R^{N}} \hskip -8mm
|{\mathbf  v}|^{\beta-2}  p \, d{\mathbf  x} d{\mathbf  v}
 \nonumber \\
 - \int_{I\!\!R^{N} \times I\!\!R^{N}} \hskip -8mm (\beta k+a) |{\mathbf  v}|^{\beta} 
 p \, d{\mathbf  x} d{\mathbf  v}
+ \beta \int_{I\!\!R^{N} \times I\!\!R^{N}} \hskip -8mm 
{\mathbf F} \cdot {\mathbf v} |{\mathbf v} |^{\beta-2} p \, d{\mathbf  x} d{\mathbf  v}
=I_1+I_2+I_3.
\label{bounddiffbeta3N}
\end{eqnarray}
Applying H\"older's inequality:
\begin{eqnarray}
\left| \int_{I\!\!R^{N}} 
{\mathbf F} \cdot \left[\int_{I\!\!R^{N}} {\mathbf v}  |{\mathbf v}|^{\beta-2} p  \, 
 d{\mathbf  v}\right] d{\mathbf  x} \right| \leq 
\|\mathbf F \|_{L^r_{\mathbf x}} 
\| \int_{I\!\!R^{N}} |{\mathbf v} |^{\beta-1} p d{\mathbf  v} \|_{L^{r'}_{\mathbf x}}, \label{IFbeta}
\end{eqnarray}
with ${1 \over r} + {1 \over r'}= 1.$ 
Thanks to (\ref{intdc}), (\ref{convolution}) and (\ref{decayheat}):
\begin{eqnarray}
\| {\mathbf F} (t) \|_{L^r_{\mathbf x}} \leq d_1 \|\nabla_{\mathbf x} c_0 \|_{L^r_{\mathbf x}}
\!+\!  d_1 \eta C_{N,q}\| c_0\|_{L^\infty_{\mathbf x}} 
\hskip -2mm \int_0^t \hskip -2mm (t-s)^{-{1\over 2} - {N\over 2}(1 - {1\over q})} 
\| { j}(s) \|_{L^{q'}_{\mathbf x}} ds, \nonumber \\
\leq d_1 \|\nabla_{\mathbf x} c_0 \|_{L^r_{\mathbf x}}
+  d_1 \eta C_{N,r,q'}\| c_0\|_{L^\infty_{\mathbf x}} t^{{1\over 2} - {N\over 2}({1\over q'} 
- {1\over r})} 
\| { j} \|_{L^\infty(0,t;L^{q'}_{\mathbf x})},
\label{Fbeta}
\end{eqnarray}
where $1+{1\over r}={1\over q}+{1\over q'},$ provided $1/N > 1 - 1/q = 1/q'-1/r$,
or equivalently, $1+1/N > 1/q'+1/r'.$ 

By Lemma 4.2, we have the estimates:
\begin{eqnarray}
\| |\mathbf v|^{\beta-2}  p\|_{L^1_{\mathbf x \mathbf v}} \leq
\|p\|_{L^1_{\mathbf x \mathbf v}}^{2 \over \beta}
\||\mathbf v|^{\beta}   p\|_{L^1_{\mathbf x \mathbf v}}^{\beta -2 \over \beta},
\quad \beta\! \geq  2,
\label{mbeta2} \\
\| \int_{I\!\!R^{N}} |{\mathbf v} |^{\beta-1} p d{\mathbf  v} \|_{L^{r'}_{\mathbf x}}
\leq C_{N,\beta}  \|p\|_{\infty}^{1 \over N+\beta}
\||\mathbf v|^{\beta}   p\|_{L^1_{\mathbf x \mathbf v}}^{N+\beta-1 \over N+\beta},
\quad r'\!=\!{N\!+\!\beta \over N\!+\!\beta\!-\!1},
\label{mbeta1} \\
\|{ j} \|_{L^{q'}_{\mathbf x}} \leq 
\| \int_{I\!\!R^{N}} |{\mathbf v} | p d{\mathbf  v} \|_{L^{q'}_{\mathbf x}}
\leq C_{N,\beta} \|p\|_{\infty}^{\beta -1 \over N+\beta}
\||\mathbf v|^{\beta}  p\|_{L^1_{\mathbf x \mathbf v}}^{N+1\over N+\beta},
\quad q'\!=\!{N\!+\!\beta \over N\!+\!1}.
\label{m1}  
\end{eqnarray}
In view of (\ref{mbeta1}) and (\ref{m1}), inequality (\ref{Fbeta}) holds provided:
\begin{eqnarray*}
{1\over q'} +{1\over r'} = {N+1\over N+\beta} + {N+\beta-1\over N+\beta} 
 = 1 + {N \over N+\beta}<1 + {1\over N},
\end{eqnarray*}
that is, $N^2-N<\beta$.
Let us analyze each term in the right hand side of identity (\ref{bounddiffbeta3N}),
keeping track of the dependence on $\beta$. Using Young's inequality 
\cite{brezis},
\begin{eqnarray}
|I_1| \leq  \varepsilon_1^{2-\beta \over 2}
\left[ \beta(\beta-2+N) \sigma \right]^{\beta \over 2}
\|p\|_{L^1_{\mathbf x \mathbf v}} 
+ \varepsilon_1 \||\mathbf v|^{\beta}   p\|_{L^1_{\mathbf x \mathbf v}}.
\label{momentoI1}
\end{eqnarray}
We choose $\varepsilon_1={\beta k \over 2}$ to cancel the moment contribution
with the negative terms in equation (\ref{bounddiffbeta3N}). The integral
$I_3$ is bounded gathering the contributions from $\nabla_{\mathbf x}
c_0$ and ${ j}$. For the $\nabla_{\mathbf x} c_0$ part, we have the 
upper estimate: 
\begin{eqnarray} \varepsilon_2^{1-N-\beta}
\left[ \beta d_1 \|\nabla_{\mathbf x} c_0 \|_{L^{N+\beta}_{\mathbf x}} C_{N,\beta} \right]^{N+\beta}
\|p\|_{\infty} + \varepsilon_2
\||\mathbf v|^{\beta}   p\|_{L^1_{\mathbf x \mathbf v}},
\label{momentoI31}
\end{eqnarray}
using again Young's inequality with $\varepsilon_2={\beta k \over 4}$ to cancel
the moment contribution with the negative terms. For the positive part
of $I_2$, Young's inequality implies:
\begin{eqnarray}
\|a^-Ê\|_{\infty} \||\mathbf v|^{\beta}  p(t)\|_{L^1_{\mathbf x \mathbf v}} \leq 
{N \over 2N + \beta} \|a^-Ê\|_{\infty}^{2N+\beta \over N}  + {N+\beta \over 2N + \beta}
\||\mathbf v|^{\beta}  p(t)\|_{L^1_{\mathbf x \mathbf v}}^{2N +\beta \over N+\beta},  
\label{split}
\end{eqnarray}
with $\| a^- \|_{\infty}= \alpha_1 \|\rho \|_{\infty}$.
Inserting (\ref{IFbeta})-(\ref{split}) in (\ref{bounddiffbeta3N}), we find:
\begin{eqnarray}
{d \over dt} \||\mathbf v|^{\beta}  p(t)\|_{L^1_{\mathbf x \mathbf v}}
 \leq \left[\beta k / 2\right]^{2-\beta \over 2}
\left[ \beta(\beta-2+N) \sigma \right]^{\beta \over 2}
\|p(t)\|_{L^1_{\mathbf x \mathbf v}}  + \nonumber \\
\left[\beta k / 4\right]^{1-N-\beta }
\left[ \beta d_1 \|\nabla_{\mathbf x} c_0 \|_{L^{N+\beta}_{\mathbf x}} C_{N,\beta} 
\right]^{N+\beta} \|p(t)\|_{\infty}  + 
{N \over 2N + \beta} \|a^-Ê\|_{\infty}^{2N+\beta \over N} +
\nonumber \\
\left[  \beta d_1 \eta \tilde C_{N,r,q'}
  \| c_0\|_{L^\infty_{\mathbf x}} t^{{1\over 2} - {N\over 2}{N \over
N+\beta}} \tilde C_{N,\beta} \|p(t) \|_{\infty}^{\beta \over N + \beta} + 1 \right]
\||\mathbf v|^{\beta}  p(t)\|_{L^1_{\mathbf x \mathbf v}}^{2N +\beta \over N+\beta}.
\label{bounddiffbeta3N2}
\end{eqnarray}
We will get a superlinear Gronwall inequality with exponent
$1 + {N \over N+\beta}>1$, unless we have
external information on ${ j}$. Let us set
\begin{eqnarray}
A_{\beta}= \||\mathbf v|^{\beta}  p_0\|_{L^1_{\mathbf x \mathbf v}} +
\left[\beta k / 2\right]^{2-\beta \over 2}
\left[ \beta(\beta-2+N) \sigma \right]^{\beta \over 2}
T\|p\|_{L^{\infty}_t L^1_{\mathbf x \mathbf v}}  + \nonumber \\
\left[\beta k / 4\right]^{1-N-\beta \over 2}
\left[ \beta d_1 \|\nabla_{\mathbf x} c_0 \|_{L^{N+\beta}_{\mathbf x}} C_{N,\beta} 
\right]^{N+\beta} \hskip -6mm
T\|p\|_{L^\infty_t L^\infty_{\mathbf x \mathbf v}}  \!+ \!
{N T \over 2N \!+ \!\beta} \left[ \alpha_1\| \rho\|_{\infty} \right]^{2N + \beta \over N} ,
\label{momentoA} \\
B_{\beta}=   
\beta d_1 \eta \tilde C_{N,r,q'}
\| c_0\|_{L^\infty_{\mathbf x}} T^{{1\over 2} - {N\over 2}
{N\over N+ \beta}} \tilde C_{N,\beta} 
\|p\|_{L^{\infty}_t L^\infty_{\mathbf x \mathbf v}}^{\beta \over N + \beta} 
\!+\! 1.
\label{momentoB}
\end{eqnarray}
Integrating (\ref{bounddiffbeta3N2})  in time we have the inequality:
\begin{eqnarray}
\||\mathbf v|^{\beta}  p(t)\|_{L^1_{\mathbf x \mathbf v}}
\leq A_{\beta} + B_{\beta} \int_0^t \||\mathbf v|^{\beta}  p(s)
\|_{L^1_{\mathbf x \mathbf v}}^{2N +\beta \over N+\beta} ds.
\label{diffmoment}
\end{eqnarray}
Set $\phi(t) = A_{\beta} + B_{\beta} \int_0^t F(s)^{1+\delta} ds \geq 0$,
with $F(t)= \||\mathbf v|^{\beta}  p(t)\|_{L^1_{\mathbf x \mathbf v}}$
and $\delta = {N \over N + \beta}.$ Then,
$\phi'(t)= B_\beta F(t)^{1+\delta}$, and 
$F(t)= \left({\phi'(t) \over B_\beta}\right)^{1\over 1+\delta}$.
Therefore, $\phi'(t) \leq B_\beta \phi(t)^{1+\delta}, \; \phi(0)=A_\beta,$ and
$F(t) \leq \phi(t)$ is bounded from above by the solution of the equation,
${1\over (- B_\beta \delta t + A_\beta^{-\delta})^{1\over \delta}},$ 
until it blows up a time $\tau_\beta= (A_\beta^{\delta} B_\beta \delta )^{-1}$.
This yields:
\begin{eqnarray}
\||\mathbf v|^{\beta}  p(t)\|_{L^1_{\mathbf x \mathbf v}}
\leq {1\over (- {N\over N+\beta}  B_{\beta} t + 
(A_{\beta})^{-{N\over N+\beta}})^{N+\beta\over N}},
\label{intmoment}
\end{eqnarray}
for $t < \tau_{\beta}= {N\!+\! \beta \over (A_{\beta})^{N\over N+\beta} B_{\beta} N}.$

Then, Lemma 4.2 combined with interpolation inequalities in $L^q$ spaces
 \cite{brezis}, provides the required estimates on  the norms
 $\||{\mathbf v}| p\|_{L^\infty(0,\tau;L^1_{\mathbf x \mathbf v})}$,  
 $\|{ j}\|_{L^\infty(0,\tau;L^q_{\mathbf x})}$, 
 $1 \leq q \leq {N+\beta \over N+1},$ and 
 $\|\tilde p\|_{L^\infty(0,\tau;L^q_{\mathbf x})}$, 
 $1 \leq q \leq {N+\beta \over N},$ for $\tau < {\rm min}(\tau_\beta,T).$ \\
 
This variation on the available estimates of 
$\| { j} \|_{L^q_{\mathbf x}}$
affects Proposition 4.4, that takes now the following form. \\
 
{\bf Proposition 6.2.} {\it Let $p\geq 0$ be a solution of (\ref{eq:linpaj})-(\ref{eq:linp0aj}). Assume that conditions (i)-(iii) in Proposition 6.1 hold. 
Assume further that:
\begin{itemize}
\item[(iv)]
$(1+|{\mathbf v}|^2)^{\beta/2} p_0({\mathbf x},{\mathbf v}) \in
L^{\infty}(I\!\!R^{N} \times I\!\!R^{N}),  \quad
\beta > {max}(N+1,N^2-N).$
\end{itemize}
Then,  
$\| (1+|{\mathbf v}|^{2})^{\beta/2} p \|_{L^\infty(0,T;L^{\infty}_{\mathbf x \mathbf v})
}$, $\||{\mathbf v}| p\|_{L^\infty(0,T;L^\infty_{\mathbf x}L^1_{\mathbf v})}$ 
and $\| p\|_{L^\infty(0,T;L^\infty_{\mathbf x}L^1_{\mathbf v})}$ 
are bounded up to any time $\tau < {min}(\tau_\beta,T)$,
by constants depending on the  parameters $\sigma$, $k$, $d$, $\eta$,  
$d_1$, $\alpha_1$, $T$, $\beta$, $N$, and the norms
$\| (1+|{\mathbf v}|^{2})^{\beta/2} p_0 \|_{L^{\infty}_{\mathbf x \mathbf v
}}$, $\|c_0\|_{L^{\infty}_{\mathbf x}}$, 
$\|\nabla_{\mathbf x}c_0\|_{L^\infty_{\mathbf x} \cap L^{N+\beta}_{\mathbf x}}$,
$\|\rho \|_{L^{\infty}_{\mathbf v}}$,
$\|p\|_{L^\infty(0,T;L^{\infty}_{\mathbf x \mathbf v}\cap
L^{1}_{\mathbf x \mathbf v})}$.}\\


{\bf Proof.} We revisit the proof of Proposition 4.4 to find
\begin{eqnarray}
\| Y(t) \|_{L^\infty_{\mathbf x \mathbf v}} 
\!\!\leq \!\! \| Y(0) \|_{L^\infty_{\mathbf x \mathbf v}} 
\!\!+\!\!\! \int_0^t \hskip -2mm\Big( \!\!
(N k\!+\!\|a^-\|_{\infty} \!) \| Y \|_{L^\infty_{\mathbf x \mathbf v}} 
\!\!+\!\! \|R_1 \|_{L^\infty_{\mathbf x \mathbf v}} 
\!\!+\!\! \|R_2 \|_{L^\infty_{\mathbf x \mathbf v}} 
\!\!+\!\! \|R_3 \|_{L^\infty_{\mathbf x \mathbf v}} \!\! \Big) ds, \nonumber
\end{eqnarray}
with  $\|R_2\|_{L^\infty_{\mathbf x \mathbf v}} \leq   
  k \beta \|Y \|_{L^\infty_{\mathbf x \mathbf v}}, $
$ \|R_3\|_{L^\infty_{\mathbf x \mathbf v}} \leq   \sigma \beta (\beta + 2 + N) 
 \|Y \|_{L^\infty_{\mathbf x \mathbf v}}, $ and
\begin{eqnarray*}
 \|R_1\|_{L^\infty_{\mathbf x \mathbf v}} \leq  
 \beta \|(1+|{\mathbf v}|^2)^{\beta/2 -1}  {\mathbf F} \cdot {\mathbf v} p 
 \|_{L^\infty_{\mathbf x \mathbf v}}. \nonumber
\end{eqnarray*}

We resort to the interpolation inequality (\ref{boundinterp}) 
in Lemma 4.2. For $\beta >1$
\begin{eqnarray}
\|(1+|{\mathbf v}|^2)^{\beta/2 -1}  | {\mathbf v} |p \|_{L^\infty_{\mathbf x \mathbf v}}
\leq \|(1+|{\mathbf v}|^2)^{(\beta - 1)/2}  p \|_{L^\infty_{\mathbf x \mathbf v}}
\nonumber \\
\leq C_1(\beta)  \| p \|_{L^\infty_{\mathbf x \mathbf v}}^{1/\beta} \|(1+|{\mathbf v}|^2)^{\beta/2}  p \|_{L^\infty_{\mathbf x \mathbf v}}^{1-1/\beta}.
\label{boundF1}
\end{eqnarray}
We must now relate $\| {\mathbf F}  \|_{\infty}$ to some norm of 
${ j}$, where
${\mathbf F}$ is given by formulas (\ref{eq:alphaF}) and (\ref{intc})-(\ref{intdc}).  
Lemma 4.1 yields inequality (\ref{boundFinf}) in terms of 
$\| { j} \|_{L^q_{\mathbf x}}$ for $q>N,$ $t \in [0,T]$.  
Proposition 6.1 yields
$\| { j} \|_{L^{\infty}(0,\tau;L^{q_\beta}_{\mathbf x})} 
\leq C(\|p \|_{L^{\infty}(0,T;L^\infty_{\mathbf x \mathbf v}
\cap L^1_{\mathbf x \mathbf v} )},I)$, where $I$ stands for the remaining list of parameters and data norms detailed in that proposition, $q_\beta={N+\beta
\over N+1}$ and $\tau < {\rm min}(\tau_\beta,T)$. 
Then, $L^q$ norms can be estimated by means 
of $L^{q_\beta}$ and $L^{\infty}$ norms for $t\leq \tau $:
\begin{eqnarray*}
\| { j} \|_{L^q_{\mathbf x}}^q  \leq 
\int d{\mathbf x} \left( \int d{\mathbf v} \, |{\mathbf v}| p \right)^{q_\beta}
\| \int d{\mathbf v} \, |{\mathbf v}| p \|_{L^\infty_{\mathbf x}}^{q-q_\beta}   
\nonumber \\
\leq \||\mathbf v | p \|^{q_\beta}_{L^{q_\beta}_{\mathbf x} L^1_{\mathbf v}} 
\| \int d{\mathbf v} \, |{\mathbf v}| p \|_{L^\infty_{\mathbf x}}^{q-q_\beta}.
\label{Lq2inf}
\end{eqnarray*}
Thanks to inequality (\ref{boundvinf}) in Lemma 4.2:
\begin{eqnarray}
\| { j} \|_{L^q_{\mathbf x}}  \leq 
C(\|p \|_{L^{\infty}(0,T;L^\infty_{\mathbf x \mathbf v}
\cap L^1_{\mathbf x \mathbf v} )},I)^{q_\beta \over q}
\| \int d{\mathbf v} \, |{\mathbf v}| p \|_{L^\infty_{\mathbf x}}^{1-
{q_\beta\over q}} \leq \nonumber \\
C_{\beta} \tilde C(\|p \|_{L^{\infty}(0,T;L^\infty_{\mathbf x \mathbf v}
\cap L^1_{\mathbf x \mathbf v} )}, I, q)
\| (1+|{\mathbf v}|^2)^{\beta/2} p\|_{L^\infty_{\mathbf x \mathbf v}}^{{N+1\over \beta}
(1-{q_\beta \over q})},
\label{Lq2infbis}
\end{eqnarray}
for $\beta > N+1$. 
To ensure ${N+1\over \beta}(1-{q_\beta\over q})\leq 
{1\over \beta}$, we may select $q$ in the interval $N<q \leq q_\beta {N+1 \over N}.$  We choose $q=q_\beta {N+1 \over N}$. Notice that $q_\beta ={N+\beta
\over N+1} > {N^2 \over N+1}$ because $\beta >N^2-N$.

We set 
\begin{eqnarray}
A'\!=\! \| Y(0) \|_{\infty}, \; 
B'\!=\! \sigma \beta (\beta \!+\! 2 \!+\! N)\!+\! 
(N\!+\! \beta) k\!+\!\|a^-\|_{\infty}
\!+\! \beta d_1 N \| \nabla_{\mathbf x} c_0 \|_{L^\infty_{\mathbf x}},   \label{ABY} \\
C'\!=\! d_1 \beta \eta \|c_0\|_{L^\infty_{\mathbf x}}  C(T,\beta,N)
\tilde C(\|p \|_{L^{\infty}(0,T;L^\infty_{\mathbf x \mathbf v}
\cap L^1_{\mathbf x \mathbf v} ) }, I). \label{CY}
\end{eqnarray}
Combining (\ref{boundF1}), (\ref{boundFinf}) and (\ref{Lq2infbis}), we find:
\begin{eqnarray}
\| Y(t) \|_{L^\infty_{\mathbf x \mathbf v}} \!\leq \!   A' \!+\!\!  
\int_0^t  \hskip -1.5mm B' \| Y (s) \|_{L^\infty_{\mathbf x \mathbf v}} ds \!+\!\!
\int_0^t  \hskip -1.5mm C' \| Y (s)\|_{L^\infty_{\mathbf x \mathbf v}}^{1-{1\over\beta}+{N+1\over \beta}(1-{q_\beta\over q})} ds. \label{Yfinal}
\end{eqnarray}
Since we have chosen $q=q_\beta {N+1 \over N}$, we have ${N+1\over \beta}(1-{q_\beta\over q}) + 1-{1\over \beta}=1$, and the right hand side grows linearly. Gronwall's inequality implies then:
\begin{eqnarray}
\| Y(t) \|_{\infty} \leq A' e^{(B'+C') t}, \quad t \in [0,T],
\nonumber 
\end{eqnarray}
with constants depending on the values specified in (\ref{ABY})-(\ref{CY}).

Once the velocity decay has been established, the $L^{\infty}$ bounds on 
$\int_{I\!\! R^{N}} p d {\mathbf v}$ and  $\int_{I\!\! R^{N}} 
{|{\mathbf v }| }p d {\mathbf v}$ follow from Lemma 4.2.
\\

We can now state the existence result for the flux without velocity cut-off. \\

{\bf Theorem 6.3.} {\it Let us assume that
\begin{eqnarray}
p_0 \geq 0, c_0 \geq 0, \label{signo2} \\
c_0  \in L^{\infty}(I\!\!R^{N}), \nabla_{\mathbf x} c_0 \in L^{\infty}\cap L^{N+\beta_2}(I\!\!R^{N}),  \label{c0inf2}  \\
(1+|{\mathbf v}|^2)^{\beta_1/2} p_0 \in L^{\infty}(I\!\!R^{N}\times I\!\!R^{N}), \quad  \beta_1 >{max}(N\!+\!1, N^2-N),
\label{p0inf2}  \\ 
(1+|{\mathbf v}|^2)^{\beta_2/2} p_0\in L^1(I\!\!R^{N}\times I\!\!R^{N}), \quad \beta_2>{max}(N\!+\!2,N^2\!-\!N). \label{p0mom2}  
\end{eqnarray}
Then, there exists a nonnegative solution $(p, c)$ of (\ref{eq:p})-(\ref{eq:alphaF}),
(\ref{flux}), satisfying:
\begin{eqnarray}
c \in L^{\infty}(0,\tau;L^{\infty}(I\!\!R^{N})), 
\nabla_{\mathbf x} c \in L^{\infty}(0,\tau;L^{\infty}\cap L^{N+\beta_2}(I\!\!R^{N})), \label{dcinf2} \\
p \in L^{\infty}(0,\tau;L^{\infty}\cap L^1(I\!\!R^{N}\times I\!\!R^{N})),
\nabla_{\mathbf v }p \in L^{2}(0,\tau;L^2(I\!\!R^{N}\times I\!\!R^{N})), \label{pinf2} \\
(1+|{\mathbf v}|^2)^{\beta_1/2}p \in L^{\infty}(0,\tau;L^{\infty}(I\!\!R^{N} \!\times I\!\!R^{N})),  \label{pdecay2} \\
(1+|{\mathbf v}|^2)^{\beta_2/2} p \in L^{\infty}(0,\tau;L^1(I\!\!R^{N}\times I\!\!R^{N})), \label{pmom2}  \\
p, |{\mathbf v}| p  \in L^{\infty}(0,\tau;L^\infty_{\mathbf x}(I\!\!R^{N},
L^1_{\mathbf v}( I\!\!R^{N})), \label{pmominf2}  
\end{eqnarray}
for any $\tau < {\rm min}(T,\tau_{\beta_2})$, with norms bounded in term of the norms of the data.\\
Additionally, { if
$(1+|{\mathbf v}|) \nabla_{\mathbf v} p \in L^{\infty}(0,\tau;L^\infty_{\mathbf x}(I\!\!R^{N}, L^1_{\mathbf v}( I\!\!R^{N}))$, then the solution is unique. In particular, uniqueness holds when $(1+|{\mathbf v}|)\nabla_{\mathbf v} p_0 \in  L^\infty_{\mathbf x}(I\!\!R^N,L^1_{\mathbf v}( I\!\!R^N)) $.}
}\\

{\bf Proof.} The proof proceeds with the same steps as the proof of Theorem 5.1.

Step 1 remains the same, except that $a_m$ and ${ j_m}$ are both bounded functions as a consequence of Proposition 4.4 when $\beta_1 >N+1$. Steps 2 and 3 proceed in the same way.  

The main difference arises in Step 4, when obtaining uniform bounds on
$\|{ j_{m-1}}\|_{\infty}$, $\|\mathbf F(c_{m-1})\|_{\infty}$, 
$\|a_{m-1}\|_{\infty}$. The previous steps provide
uniform estimates on the density norms
$\|p_m\|_{L^\infty(0,T;L^{q}_{\mathbf x\mathbf v})}$ for $1 \leq q \leq \infty$
and the angiogenic factor norm
$\|c_m\|_{L^\infty(0,T;L^{\infty}_{\mathbf x})}$. 
To establish uniform estimates on $a_{m-1}$, ${\mathbf F}_{m-1}$ and 
${ j}_{m-1}$ we adapt Propositions 6.1 and 6.2  to equation (\ref{eq:pm}), 
setting $a= \gamma \int_0^t \tilde p_{m-1} \, ds -\alpha(c_{m-1}) \rho$ and
${\mathbf F}={\mathbf F}(c_{m-1})$, with ${ j}={ j}_{m-1}$ depending 
on $p_{m-1}$, not $p_m$. We replace (iii) with
\begin{itemize}
\item[(iii)']  ${\bf F}(c_{m-1})$ is given by (\ref{eq:alphaFm}), $c_{m-1}$, $\nabla_{\mathbf x} c_{m-1}$ are given by (\ref{cmint})-(\ref{dcmint}).
$c_{m-1}\geq 0$ is coupled to equation  (\ref{eq:pm}) for $p_{m-1}$ by 
(\ref{fluxm}),
and $p_m\geq 0$ is coupled to equation (\ref{eq:cm}) for $c_{m-1}$ through (\ref{eq:alphaFm}).
\end{itemize}
We know by Step 1 that $a \in L^{\infty}$. Its negative part is $a^-=\alpha(c_{m-1})\rho$, 
that satisfies $\|a^-\|_{\infty}\leq \alpha_1 \|\rho \|_{L^\infty_{\mathbf v}}.$  Step 1
also ensures that ${\mathbf F}(c_{m-1})$ and ${ j_{m-1}}$ are bounded functions.
By Step 2, we have $0 \leq c_{m-1}\leq \|c_0\|_{\infty}$.
By Step 3, $\|p_m (t)\|_{L^{1}_{\mathbf x\mathbf v}}  \leq  \|p_0\|_{L^{1}_{\mathbf x\mathbf v}}  e^{\alpha_1 \|\rho\|_{L^\infty_{\mathbf v}} T}$ and
$\|p_m (t)\|_{L^{\infty}_{\mathbf x\mathbf v}}  \leq  \|p_0\|_{L^{1}_{\mathbf x\mathbf v}}  e^{NkT} e^{\alpha_1 \|\rho\|_{L^\infty_{\mathbf v}} T}$
for all $m$ and $t\in[0,T]$.

The previous estimates on the iterates hold in $[0,T]$, for any $T>0$. Uniform estimates on velocity integrals are only guaranteed for finite time $\tau_{\beta_2}$.
Revisiting the proof of Proposition 6.1, we have identity (\ref{bounddiffbeta3N}) and inequality (\ref{Fbeta}) for $p_m$ and $\mathbf F(c_{m-1}).$ We recover inequality
(\ref{bounddiffbeta3N2}) for $p_m$ with $\||\mathbf v|^{\beta_2}  p(t)\|_{L^1_{\mathbf x \mathbf v}}^{(2N +\beta_2) / (N+\beta_2)}$ replaced by
$\||\mathbf v|^{\beta_2}  p_{m-1}(t)\|_{L^1_{\mathbf x \mathbf v}}^{(N+1) / (N+\beta_2)}
\||\mathbf v|^{\beta_2}  p_{m}(t)\|_{L^1_{\mathbf x \mathbf v}}^{(N+\beta_2-1) / (N+\beta_2)} $. Defining $A_{\beta_2}$ and $B_{\beta_2}$ as in (\ref{momentoA})-(\ref{momentoB}), but in terms of upper bounds of the norms of $p_m$, and integrating in time, we find the analogous of inequality (\ref{diffmoment}):
\[
\||\mathbf v|^{\beta_2}  p_m(t)\|_{L^1_{\mathbf x \mathbf v}}
\leq A_{\beta_2} + B_{\beta_2} \int_0^t \||\mathbf v|^{\beta_2}  p_{m-1}(s)\|_{L^1_{\mathbf x \mathbf v}}^{N+1\over N+\beta_2}
\||\mathbf v|^{\beta_2}  p_{m}(s)\|_{L^1_{\mathbf x \mathbf v}}^{N+\beta_2-1\over N+\beta_2} ds.
\]
For any natural number $M$, set
$\Phi_M= {\rm Max}_{1\leq m \leq M} {F_m},$
$F_m = \||\mathbf v|^{\beta_2}  p_m(t)\|_{L^1_{\mathbf x \mathbf v}}.$
When $m\leq M$, we have
$ F_m(t) \leq   A _{\beta_2} +   B_{\beta_2}  \int_0^t \Phi_M(s)^{1 \over q} 
\Phi_M(s)^{1 \over q'} ds =
  A_{\beta_2}  +  B_{\beta_2}  \int_0^t \Phi_M(s)^{1+\delta} ds, $
  $1/q+1/q'=1+\delta$.
Thus,
$ \Phi_M(t) \leq    A _{\beta_2} +   B_{\beta_2}  \int_0^t \Phi_M(s)^{1+\delta} ds. $
Arguing as we did in Proposition 6.1 to find inequality (\ref{intmoment}) we obtain:
$F_m(t) \leq \Phi_M(t) \leq (- B_{\beta_2} \delta t + A_{\beta_2}^{-\delta})^{-1\over \delta}$
for $t \leq \tau = (A_{\beta_2}^{\delta} B_{\beta_2} \delta)^{-1}$ and $m\leq M$, for
any $M$. Therefore,
\begin{eqnarray}
\||\mathbf v|^{\beta_2}  p_m(t)\|_{L^1_{\mathbf x \mathbf v}}
\leq {1\over (- {N\over N+\beta_2}  B_{\beta_2} t + 
(A_{\beta_2})^{-{N\over N+\beta_2}})^{N+\beta_2\over N}},
\label{momemtbeta2}
\end{eqnarray}
for $t <\tau_{\beta_2}= {N\!+\! \beta_2 \over (A_{\beta_2})^{N\over N+\beta_2} B_{\beta_2} N},$
and for all $m$.

Then, Lemma 4.2 combined with interpolation inequalities in $L^q$ spaces, and
the uniform bound on $\|p_m\|_{L^\infty_t L^\infty_{\mathbf x \mathbf v}}$, provides uniform estimates on  
 $\||{\mathbf v}| p_m(t)\|_{L^1_{\mathbf x \mathbf v}}$,  
 $\|{ j}(p_m(t))\|_{L^q_{\mathbf x}}$, $1 \leq q \leq {N+\beta_2 \over N+1},$ 
 and $\|\tilde p_m(t)\|_{L^q_{\mathbf x}}$, $1 \leq q \leq {N+\beta_2 \over N}$,
 for $t\in [0,\tau]$, $\tau< \tau_{\beta_2}$.
 Using (\ref{tildecmq}), we obtain uniform estimates on 
 $\|\tilde{c}_{m-1}(t)\|_{L^{2}_{\mathbf x}}$  for $t\in [0,\tau]$ provided
  $\|{ j}(p_m(t))\|_{L^2_{\mathbf x}}$ is bounded. Taking
  $\beta_2 \geq  N+2$, we have ${N+\beta_2 \over N+1} \geq 2$ and this is
  guaranteed.

Let us consider now Proposition 6.2. The force field ${\mathbf F}(c_{m-1})$ depends now on ${ j}_{m-1}$, not ${ j}_{m}$. The Gronwall type inequality (\ref{Yfinal})  takes  the form
\[
\| Y_m(t) \|_{L^\infty_{\mathbf x \mathbf v}} \!\leq \!   A' \!+\!\!  
\int_0^t  \hskip -1.5mm B'\| Y_m (s) \|_{L^\infty_{\mathbf x \mathbf v}} ds \!+\!\!
\int_0^t  \hskip -1.5mm C' \| Y_m (s)\|_{L^\infty_{\mathbf x \mathbf v}}^{1-1/\beta} 
\| Y_{m-1} (s)\|_{L^\infty_{\mathbf x \mathbf v}}^{1/\beta} 
ds,  
\]
where
$Y_m(t)=\| (1+|{\mathbf v}|^2)^{\beta_1/2} p_m \|_{L^\infty_{\mathbf x \mathbf v}}$,
for $ t \in[0,\tau]$.
The constants $A'$, $B'$, $C'$ are defined as in formulas (\ref{ABY})-(\ref{CY}), but
in terms of upper bounds of the norms of $p_m$. Reasoning as in reference \cite{degond}, or arguing as we did above to obtain (\ref{momemtbeta2}) we find  
\begin{eqnarray}
\|(1+|{\mathbf v}|^2)^{\beta_1/2} p_m (t) \|_{L^\infty_{\mathbf x \mathbf v}} \leq A' e^{(B'+C')t},  \quad t\in [0,\tau],
\label{momemtbeta1}
\end{eqnarray} 
for all $m$. Combining this estimate with Lemma 4.2, we conclude that
\[
\||{\mathbf v}| p_m\|_{L^\infty(0,\tau;L^{\infty}_{\mathbf x} L^1_{\mathbf v})}, 
\|p_m\|_{L^\infty(0,\tau;L^{\infty}_{\mathbf x} L^1_{\mathbf v})}, 
\]
are uniformly bounded by constants depending on fixed parameter values and on $\|c_0\|_{W^{1,\infty}_{\mathbf x}}$, 
$\| \nabla_{\mathbf x} c_0\|_{L^{N+\beta_2}_{\mathbf x}}$
$, \|\rho\|_{L^\infty_{\mathbf v}}$,  
$\|(1\!+\! |{\mathbf  v}|^2)^{\beta_2/2}  p_0 \|_{L^1_{\mathbf x \mathbf v}}$ \!\!, 
$\|(1\!+\! |{\mathbf  v}|^2)^{\beta_1/2}  p_0 \|_{L^\infty_{\mathbf x \mathbf v}}$.

Therefore,  the coefficients $a_{m-1}$, $\alpha(c_{m-1})$, ${ j}_{m-1}$ and $F(c_{m-1})$ appearing in the equations are uniformly bounded in 
$L^{\infty}(0,\tau; L^{\infty}_{\mathbf x})$. From Step 2, we also infer that  
the norms  $\| \nabla_{\mathbf x} c_m\|_{L^\infty(0,\tau; L^\infty_{\mathbf x})}$,
$\| \nabla_{\mathbf x} \tilde c_m\|_{L^\infty(0,\tau; L^{N+\beta_2}_{\mathbf x})}$, 
$\| \nabla_{\mathbf x} \tilde c_m\|_{L^\infty(0,\tau; L^{2}_{\mathbf x})}$ and
 $\| \tilde c_m\|_{L^\infty(0,\tau; L^2_{\mathbf x})}$
are uniformly bounded too. 

Steps 5 and 6 proceed as in Theorem 5.1, providing similar compactness results up to the time $\tau$, that allow to pass to the limit in the weak formulation of the equations yielding a solution $(p,c)$ of the system up to time $\tau$. The same regularity considerations hold as in Step 7 in Theorem 5.1.

Uniqueness follows the same argument as in
Step 8  in Theorem 5.1.
The proof of uniqueness holds as long as solutions satisfying the required 
regularity conditions exist.
Let us assume that we have two solutions $(p_1,c_1)$ and $(p_2,c_2)$,
satisfying such conditions up to a time $T$.
We recover identities and inequalities (\ref{eq:pdif})-(\ref{intdc1dif}).
{
Let us first notice that $\| j(\overline{p}) \|_{L^1_{\mathbf x}} \leq 
\| |\mathbf v| \overline{p}\|_{L^1_{\mathbf x \mathbf v}}$.
Instead of Lemma 4.3, we must resort to integral expressions of solutions
of (\ref{eq:pdif}) in terms of fundamental solutions to relate 
$\| |\mathbf v| \overline{p}\|_{L^1_{\mathbf x \mathbf v}}$ and 
$\| \overline{p}\|_{L^1_{\mathbf x \mathbf v}}$. 
Let us set
\begin{eqnarray}
f\!=\![-\gamma a(p_{1}) \!+\! \alpha(c_{1}) \rho] \overline{p}  
 \!-\!   [{\mathbf F}(c_{1})\!-\!{\mathbf F}(c_{2})] \nabla_\mathbf{v}  
p_2 \!+\! [-\gamma a(\overline{p}) \!+\! (\alpha(c_{1})\!-\!\alpha(c_2)) \rho] p_2.
\nonumber
\end{eqnarray}
Using the fundamental solution $\Gamma=\Gamma_{\mathbf F(c_1)}$
and identity (\ref{solint})
\begin{eqnarray} 
\int_{I\!\!R^N \times I\!\!R^N}  \hskip -4mm
|\mathbf v| |\overline{p}(t, {\mathbf x}, {\mathbf v})| d \mathbf v
d \mathbf x \leq  \hskip -2mm
\int_0^t \hskip -2mm
\int_{I\!\!R^N\times I\!\!R^N\times I\!\!R^N \times I\!\!R^N}  
\hskip - 2.5cm |\mathbf v|
\Gamma(t, {\mathbf x}, {\mathbf v};  \tau, {\boldsymbol \xi}, {\boldsymbol \nu})  
|f(\tau, {\boldsymbol \xi}, {\boldsymbol \nu})| d{\boldsymbol \xi} 
d{\boldsymbol \nu} d \tau d \mathbf x d \mathbf v = I. \nonumber
\end{eqnarray}
Thanks to estimate (F4) in Lemmas 2.2 and 2.4, $\Gamma$ is bounded from
above by the fundamental solution $G$ of the free field linear Fokker-Planck
operator, rescaled in space and velocity, and multiplied by a constant.
By assertion (e) in Lemma 1.3 in reference \cite{rein}
\[
\int_{I\!\!R^N \times I\!\!R^N} |\mathbf v| G(t, {\mathbf x}, {\mathbf v};  \tau, {\boldsymbol \xi}, {\boldsymbol \nu}) d \mathbf x d \mathbf v 
\leq C_1 (1+|{\boldsymbol \nu}|),
\]
for ${\mathbf x},$ $ {\mathbf v}, $ ${\boldsymbol \xi}, $ ${\boldsymbol \nu},$
and $t > \tau \geq 0$. This is also easily checked using expression (12)-(15)
for the fundamental solution as in references \cite{bouchut3,gema}. Therefore,
$I \leq I_1 + I_2 $ where
\begin{eqnarray}
I_1 = C_2  \int_0^t \hskip -2mm \int_{I\!\!R^N\times I\!\!R^N}  \hskip - 1cm 
|f(\tau, {\boldsymbol \xi}, {\boldsymbol \nu})| d{\boldsymbol \xi} 
d{\boldsymbol \nu} d \tau, \quad
I_2 = C_3
\int_0^t \int_{I\!\!R^N\times I\!\!R^N}  \hskip - 1cm
|{\boldsymbol \nu}| |f(\tau, {\boldsymbol \xi}, {\boldsymbol \nu})| d{\boldsymbol \xi} 
d{\boldsymbol \nu} d \tau. \nonumber 
\end{eqnarray}
Both integrals can be estimated in a similar way to the right hand side in (\ref{intp1dif}), making use of inequalities (\ref{meanF})-(\ref{meanalpha}). This
yields
\begin{eqnarray}
\||\mathbf v|\overline{p}(s) \|_{L^1_{\mathbf x \mathbf v}}\! \leq  \!
[\gamma \|a(p_{1})\|_{\infty} \!+\! \alpha_1 \|\rho\|_{\infty}]  \!\!
\int_0^t \! \!\! ds [C_2\|\overline{p}(s) \|_{L^1_{\mathbf x \mathbf v}} 
\! \! +\!  C_3 \||\mathbf v|\overline{p}(s) \|_{L^1_{\mathbf x \mathbf v}}]
\! +\!  \nonumber  \\
d_1 q_1 \gamma_1 \|\nabla_{\mathbf x} c_1\|_{\infty} 
[C_2 \| \nabla_\mathbf{v} p_2 \|_{L^{\infty}_t L^{\infty}_{\mathbf x}L^1_{\mathbf v}} \! + \! 
C_3 \| |\mathbf v|\nabla_\mathbf{v} p_2 \|_{L^{\infty}_t L^{\infty}_{\mathbf x}L^1_{\mathbf v}} ] \! \! 
\int_0^t \!\!\!  ds \|\overline{c}(s)\|_{L^1_{\mathbf x}} 
\! +\!  \nonumber \\
d_1 
[C_2 \| \nabla_\mathbf{v} p_2 \|_{L^{\infty}_t L^{\infty}_{\mathbf x}L^1_{\mathbf v}} + 
C_3 \| |\mathbf v|\nabla_\mathbf{v} p_2 \|_{L^{\infty}_t L^{\infty}_{\mathbf x}L^1_{\mathbf v}} ]
\int_0^t \!\! ds \|\nabla_{\mathbf x}\overline{c}(s)\|_{L^1_{\mathbf x}}
\!+\!   \nonumber \\
\gamma [C_2 \|p_2\|_{L^\infty_t L^{\infty}_{\mathbf x} L^1_{\mathbf v}} 
+ C_3 \||\mathbf v|p_2\|_{L^\infty_t L^{\infty}_{\mathbf x} L^1_{\mathbf v}} ]
\int_0^t \!\! ds \!\! \int_0^s \!\! d\tau 
\|\overline{p}(\tau)\|_{L^{1}_{\mathbf x \mathbf v}} \!+\! 
\nonumber \\
{\alpha_1 \|\rho\|_{\infty}\over c_R} 
[C_2 \|p_2\|_{L^\infty_t L^{\infty}_{\mathbf x} L^1_{\mathbf v}} 
+ C_3 \||\mathbf v|p_2\|_{L^\infty_t L^{\infty}_{\mathbf x} L^1_{\mathbf v}} ]
\int_0^t \!\! ds  \|\overline{c}(s)\|_{L^1_{\mathbf x}}. 
\label{intp2dif}
\end{eqnarray}
}

Let us set:
\begin{eqnarray}
A= A_1 +A_2 = [\gamma \|a(p_{1})\|_{\infty} 
\!+\! \alpha_1 \|\rho\|_{\infty}]  
+ T \gamma \|p_2\|_{L^\infty_t L^{\infty}_{\mathbf x}L^1_{\mathbf v}},
\nonumber \\
B= B_1 +B_2 =  d_1 q_1 \gamma_1 \|\nabla_{\mathbf x} c_1\|_{\infty} 
\| \nabla_\mathbf{v}   p_2 \|_{L^{\infty}_t L^{\infty}_{\mathbf x}L^1_{\mathbf v}}
+ {\alpha_1 \|\rho\|_{\infty}\over c_R} \|p_2\|_{L^\infty_t L^{\infty}_{\mathbf x} L^1_{\mathbf v}}, 
\nonumber \\
D= d_1  \| \nabla_\mathbf{v}   p_2 \|_{L^{\infty}_t L^{\infty}_{\mathbf x}L^1_{\mathbf v}},
\nonumber \\
E = E_1 + 2 M \eta T^{1/2} E_2 =  2 \eta M \|c_2\|_{\infty} T^{1/2} 
+  2 \eta^2 M \|c_2\|_{\infty} C(\|{ j}(p_1)\|_{\infty}) T^{3/2}.
\nonumber
\end{eqnarray}
We denote by $\tilde A, \tilde A_i, \tilde B, \tilde B_i, \tilde D$ the same constants replacing the norms
$\| \nabla_\mathbf{v}   p_2 \|_{L^{\infty}_{\mathbf x}L^1_{\mathbf v}}$,
$\|   p_2 \|_{L^{\infty}_{\mathbf x}L^1_{\mathbf v}}$
with
$\|{ |{\mathbf v}| }\nabla_\mathbf{v}   p_2 \|_{L^{\infty}_{\mathbf x}L^1_{\mathbf v}}$,
$\|{ |{\mathbf v}|}  p_2 \|_{L^{\infty}_{\mathbf x}L^1_{\mathbf v}}$.
We define now 
\begin{eqnarray}
U(t)= {\rm max}_{s\in[0,t]} \left\{ \|\overline{p}(s) \|_{L^1_{\mathbf x \mathbf v}},  
 \|{|{\mathbf v}|} \overline{p}(s) \|_{L^1_{\mathbf x \mathbf v}}
\right\}.
\label{controldif2}
\end{eqnarray}
Combining (\ref{intp1dif})-(\ref{intc1dif}) with (\ref{intp2dif}) we find: 
\begin{eqnarray}
\|\overline {p}(t)\|_{L^1_{\mathbf x \mathbf v}} \leq 
A \int_0^t ds \|\overline{p}(s) \|_{L^1_{\mathbf x \mathbf v}} 
+ ( B E_2 + D E)
\int_0^t ds \| {j}(\overline{p}) \|_{L^{\infty}(0,s;L^1_{\mathbf x})}
\label{control2} \\
\leq  \left(A + ( B E_2 + D E)  \right) \int_0^t U(s) ds,
\nonumber\\
\| {|{\mathbf v}|} \tilde {p}(t)\|_{L^1_{\mathbf x \mathbf v}} \leq 
{ \left(A + A_1+\tilde A_2 + (B+\tilde B)E_2+
(D+\tilde D) E \right) }
\int_0^t U(s) ds.
\end{eqnarray}
We deduce that $U(t)$ satisfies a Gronwall inequality of the form
\[ U(t) \leq G(T) \int_0^t U(s) ds\] for $t\in [0,T]$. Therefore, $U=0$
and $p_1=p_2$ in $[0,T]$ for any $T>0$. 
{ Estimate (\ref{intc1dif}) implies then $c_1=c_2$.}



\vskip 5mm

{\bf Acknowledgements.} This work has been supported by 
MINECO grants No. FIS2011-28838-C02-02 
and No. MTM2014-56948-C2-1-P. 
The authors thank LL Bonilla and V. Capasso for suggesting the problem and for insight on the modeling.


\vskip 5mm

\begin{thebibliography}{9}
\bibitem{anderson} A.R.A. Anderson, M.A.J. Chaplain, Continuous and discrete mathematical models of tumor induced angiogenesis, Bull. Math. Biol. 
60, 857-900, 1998.
\bibitem{aronson} D.G. Aronson, Nonegative solutions of linear parabolic equations, Ann. Sci. Norm. Sup. Pisa 22, 607-694, 1968.  
\bibitem{capasso} LL Bonilla, V Capasso, M. Alvaro, M. Carretero, Hybrid modeling of tumor induced angiogenesis, Physical Review E 90, 062716, 2014.
\bibitem{bouchut} F. Bouchut, Existence and uniqueness of a global smooth solution for the Vlasov-Poisson-Fokker-Planck system in three dimensions, J. Funct. Anal. 111, 239-258, 1993.
\bibitem{bouchut2} F. Bouchut, J. Dolbeault, On long time asymptotics of the
Vlasov-Fokker-Planck equation and of the Vlasov-Poisson-Fokker-Planck system with coulombic and newtonian potentials, Diff. Int. Eqs. 8(3), 487-514, 1995.
\bibitem{bouchut3} F. Bouchut, Smoothing effect for the nonlinear Vlasov-Poisson-Fokker-Planck system, J. Diff. Eqs 122, 225-238, 1995.
\bibitem{brezis} H. Br\'ezis, Functional analysis, Sobolev spaces and partial differential equations, Springer,  2011.
\bibitem{angiogenesis3} V. Capasso and D. Morale, Stochastic modelling of tumour-induced angiogenesis, J. Math. Biol. 58, 219-233, 2009.
\bibitem{angiogenesis1} P. F. Carmeliet, Angiogenesis in life, disease and medicine, Nature  438, 932-936, 2005.
\bibitem{angiogenesis2} P. Carmeliet and R. K. Jain, Molecular mechanisms and clinical applications of angiogenesis, Nature 473, 298-307, 2011.
\bibitem{carpio} A. Carpio, Long time behavior of solutions of the Vlasov-Poisson-Fokker-Planck equation, Mathematical Methods in the Applied Sciences 21, 985-1014, 1998.
\bibitem{gema} A. Carpio, G. Duro, Well posedness of an integrodifferential diffusion model related to angiogenesis, Applied Mathematical Modelling 40, 5560-5575, 2016
\bibitem{carrillo} J.A. Carrillo, Global solutions for the initial boundary  value problem
to the Vlasov-Poisson-Fokker-Planck system, Math. Meth. Appl. Sc., 21, 907-938, 
1998.
\bibitem{angiogenesis6} S. L. Cotter, V. Klika, L. Kimpton, S. Collins, and A. E. P.
Heazell, A stochastic model for early placental development, J. R. Soc. Interface 11, 20140149, 2014.
\bibitem{chandrasehkar} S. Chandrasehkar, Stochastic problems in physics and astronomy, Rev. Mod. Phys. 15,1-89, 1943.
\bibitem{chen} J.C. Chen, C. He, On local existence of the Vlasov-Fokker-Planck equations in a 2D anisotropic space, Boundary value problems, 2013:233, 2013
\bibitem{diperna} R.J. DiPerna and P.L. Lions, On the Fokker-Planck-Boltzmann equation, Comm. Math. Phys., 120, 1-23, 1988.
\bibitem{degond} P. Degond, Global existence of smooth solutions for the Vlasov-Fokker-Planck equation in 1 and 2 space dimensions, Ann. Sci. Ec. Norm. Super., 19(4), 519-542, 1986.
\bibitem{friedman} A. Friedman, Partial differential equations of parabolic type, Prentice-Hall, 1964.
\bibitem{giga} M.H. Giga, Y. Giga, J. Saal,  Nonlinear partial differential equations, Birkhauser, 2010.
\bibitem{ilin} A.M. Il'in, R.Z. Khasminski, On equations of Brownian motion, Theor. Prob. Appl. IX, 421-444, 1964.
\bibitem{kusuoka} S. Kusuoka, H\"older continuity of the fundamental solutions to parabolic equations with irregular coefficients, arXiv:1310.4600v2 [math.PR], 2014.
\bibitem{lions} J.L. Lions, Quelques m\'ethodes pour les probl\`emes aux limites nonlin\'eaires, Gauthier-Villards, 1969.
\bibitem{perthame} P.L. Lions, B. Perthame, Propagation of moments and regularity
for the three dimensional Vlasov-Poisson system, Invent. Math. 105, 415-430,
1991.
\bibitem{rein} G. Rein, J. Weckler,  Generic global classical solutions of the Vlasov-Fokker-Planck-Poisson system in three dimensions. J. Differential Equations 99, 59-77, 1992.
\bibitem{angiogenesis5} M. Scianna, L. Munaron, and L. Preziosi, A multiscale hybrid approach for vasculogenesis and related potential blocking therapies, Prog. Biophys. Mol. Biol. 106, 450-462, 2011.
\bibitem{simon} J. Simon, Compact sets in the space $L^p(0,T;B)$, Ann. Mat. Pura ed Applicata (IV) CXLVI, 65-96, 1987.
\bibitem{angiogenesis4} K. R. Swanson, R. C. Rockne, J. Claridge, M. A. Chaplain, E.
C. Alvord Jr., and A. R. A. Anderson, Quantifying the role of angiogenesis in malignant progression of gliomas: in silico modeling integrates imaging and histology, Cancer Res. 71, 7366-7375, 2011.
\bibitem{victoryweak} H.D. Victory,  On the existence of global weak solutions for Vlasov-Poisson-Fokker-Planck systems, J. Math. Anal. Appl., 160, 525-555, 1991.
\bibitem{victoryclassical} H.D. Victory, B.P. O'Dwyer, On classical solutions of Vlasov-Poisson-Fokker-Planck systems, Ind. Univ. Math. Math. J., 3 (1), 105-155, 1990.
\end{thebibliography}
\end{document}